\documentclass[final,3p,times,authoryear,11pt]{elsarticle}

\usepackage{amsmath}
\usepackage{graphicx}
\usepackage{enumerate}
\usepackage{natbib} 
\usepackage{url} 
\usepackage{xcolor}
\usepackage{amsfonts, amsthm, latexsym, amssymb}
\usepackage{algorithm}
\usepackage{algpseudocode}
\usepackage{booktabs}
\usepackage{multirow}
\usepackage{makecell}
\usepackage{subcaption}
\usepackage{longtable}
\usepackage{enumitem}
\usepackage{setspace}
\usepackage{lineno}
\usepackage{color, colortbl}
\usepackage{bbm}

\makeatletter
\newenvironment{breakablealgorithm}
  {
   \begin{center}
     \refstepcounter{algorithm}
     \hrule height.8pt depth0pt \kern2pt
     \renewcommand{\caption}[2][\relax]{
       {\raggedright\textbf{\ALG@name~\thealgorithm} ##2\par}%
       \ifx\relax##1\relax 
         \addcontentsline{loa}{algorithm}{\protect\numberline{\thealgorithm}##2}%
       \else 
         \addcontentsline{loa}{algorithm}{\protect\numberline{\thealgorithm}##1}%
       \fi
       \kern2pt\hrule\kern2pt
     }
  }{
     \kern2pt\hrule\relax
   \end{center}
  }
\makeatother

\newtheorem{Theorem}{Theorem}
\newtheorem{Lemma}{Lemma}
\newtheorem{Proposition}{Proposition}
\newtheorem{Assumption}{Assumption}
\newtheorem{Remark}{Remark}

\journal{European Journal of Operational Research}

\begin{document}

\makeatletter
\def\ps@pprintTitle{%
   \let\@oddhead\@empty
   \let\@evenhead\@empty
   \let\@oddfoot\@empty
   \let\@evenfoot\@empty
}
\makeatother

\onehalfspacing


\begin{frontmatter}



\title{A Budget-Adaptive Allocation Rule for Optimal Computing Budget Allocation} 


\author[label1]{Zirui Cao} 
\author[label1]{Haowei Wang\corref{cor1}}
\author[label1]{Ek Peng Chew}
\author[label1]{Haobin Li}
\author[label2]{Kok Choon Tan}

\cortext[cor1]{Corresponding author.}

\affiliation[label1]{organization={Department of Industrial Systems Engineering and Management, National University of Singapore},
            addressline={1 Engineering Drive 2}, 
            city={Singapore},
            postcode={117576}, 
            country={SINGAPORE}}

\affiliation[label2]{organization={Department of Analytics and Operations, National University of Singapore},
            addressline={15 Kent Ridge Drive}, 
            city={Singapore},
            postcode={119245}, 
            country={SINGAPORE}}

\begin{abstract}
Simulation-based ranking and selection (R\&S) is a popular technique for optimizing discrete-event systems (DESs). It evaluates the mean performance of system designs by simulation outputs and aims to identify the best system design from a set of alternatives by intelligently allocating a limited simulation budget. In R\&S, the optimal computing budget allocation (OCBA) is an efficient budget allocation rule that asymptotically maximizes the probability of correct selection (PCS). In this paper, we first show the asymptotic OCBA rule can be recovered by considering a large-scale problem with a specific large budget. Considering a sufficiently large budget can greatly simplify computations, but it also causes the asymptotic OCBA rule ignoring the impact of budget. To address this, we then derive a budget-adaptive rule under the setting where budget is not large enough to simplify computations. The proposed budget-adaptive rule determines the ratio of total budget allocated to designs based on the budget size, and its budget-adaptive property highlights the significant impact of budget on allocation strategy. Based on the proposed budget-adaptive rule, two heuristic algorithms are developed. In the numerical experiments, the superior efficiency of our proposed allocation rule is shown.
\end{abstract}




\end{frontmatter}



\section{Introduction}
\label{sec:introduction}
Discrete-event systems (DESs) are a widely-used technical abstraction for complex systems \citep{zeigler2000theory}, such as traffic control systems, manufacturing systems, and communication systems. When the complexity of DESs is high and analytical models are unavailable, a powerful tool for evaluating the performance of DESs is discrete-event system simulation \citep{banks2005discrete}. In this paper, we consider a simulation optimization problem of identifying the best system design from a set of competing alternatives, where “best” is defined with respect to the smallest mean performance. The performance of each design is unknown and can be learnt by samples, i.e., by the sample mean of simulation outputs. Such problem is often called statistical ranking and selection (R\&S) problem \citep{kim2006selecting,chen2011stochastic,fu2015handbook} or Ordinal Optimization (OO) problem \citep{ho1992ordinal,ho2008ordinal}.

In R\&S, sampling efficiency is of significant concern as the simulation budget is often limited. First, simulation can be expensive. For example, the running time for a single simulation replication of the 24-hour dynamics of the scheduling process in a transportation network with 20 intersections can be about 2 hours, and it can take 30 minutes to obtain an accurate estimate of a maintenance strategy's average cost by running 1,000 independent simulation replications for a re-manufacturing system \citep{ho2008ordinal}. Although the recent increasing computing power has considerably alleviated concerns regarding computational cost, sampling efficiency remains one of the fundamental limitations in the application of simulation optimization methods to optimize complex systems \citep{xu2016simulation}. Second, the number of alternative designs can be very large. In practical implementations, the problems we face may involve over ten thousands of designs and relative to the scale of problems, simulation budget can be limited \citep{hong2021review}. And to address such problems, high-performance computing clusters and parallel computing environments may even be required \citep{hong2022solving}. With limited simulation replications, it is impossible to guarantee a correct selection of the best design occurs with probability 1. This nature of the problem motivates the need of implementing R\&S techniques to intelligently allocate simulation replications to designs for efficiently identifying the best design. In our problem, we consider a fixed-budget setting, and the probability of correct selection (PCS), a primary criterion in R\&S literature, is used to measure the quality of budget allocation rules. The goal is to derive a budget allocation rule that can maximize the PCS subjecting to a constraint simulation budget.

Although the simulation budget is limited and of vital importance, many R\&S algorithms allocate a simulation budget either by asymptotically optimal or by one-step-ahead optimal allocation rules, both of which can not adapt to the simulation budget. Intuitively, we argue that a desirable budget allocation rule should be adaptive to the simulation budget. This intuition is consistent with theoretical analyses on the optimal budget allocation rules under different simulation budgets. The optimal computing budget allocation (OCBA) asymptotically maximizes the PCS, and it tends to allocate large budget allocation ratios to competitive designs, where competitive designs include the best design and non-best designs that are hard to distinguish from the best \citep{chen2000simulation}. However, when the simulation budget is small, assigning large budget allocation ratios to competitive designs, according to the OCBA allocation rule, may decrease the PCS \citep{peng2015non}. Such scenario is referred to as the low-confidence scenario \citep{peng2017gradient} and also takes place in the expected value of information (EVI) in \cite{chick2010sequential} and knowledge gradient (KG) polices in \cite{frazier2008knowledge} and \cite{ryzhov2012knowledge}. To avoid the decrease of PCS, the budget allocation ratios of competitive designs should be discounted and the budget allocation ratios of non-competitive designs should be increased \citep{peng2017gradient,shi2022dynamic}. This counter-intuitive result emphasizes the significant impact of simulation budget on the budget allocation rule. It motivates the need of deriving a desirable budget allocation rule that considers and adapts to the simulation budget. 

In this paper, we consider a fixed-budget R\&S problem and develop a budget-adaptive allocation rule. Our approach follows the OCBA framework, however, our derivation does not force the simulation budget to be large enough to simplify computations. First, we approximate PCS by its lower-bound derived from Bonferroni inequality; formulate the approximated budget allocation problem; and characterize its optimality conditions under a general budget. The Bonferroni bound serves as a cheap approximation of PCS and enables us to evaluate PCS in a computationally efficient way. Second, we show that the asymptotic OCBA rule in \cite{chen2000simulation} can be recovered by considering a large-scale problem with a specific large budget, i.e., the number of designs $k\to \infty$ and the simulation budget $T = \omega(k \log k)$. Third, we derive a budget-adaptive rule under the setting where $k \to \infty$ but $T = \mathcal{O}(k\log k)$. Under such setting, the simulation budget is different from and smaller than what is considered in recovering the asymptotic OCBA rule. To derive a closed-form rule, we use first-order Taylor series expansion at the asymptotic OCBA rule to linearize a set of optimality conditions. Compared with letting $T=\omega(k \log k)$, Taylor series expansion keeps more features (e.g., $T$) of the optimality conditions from being ignored and therefore is more promising \citep{shi2022dynamic}. In contrast to the asymptotic OCBA rule, the budget-adaptive rule can discount the budget allocation ratios of non-best designs that are hard to distinguish from the best design while increase the budget allocation ratios of those that are easy to distinguish. These adjustments to budget allocation ratios are based on $T$. When $T$ is sufficiently large, i.e, when $T = \omega(k \log k)$, the proposed budget-adaptive rule reduces to the asymptotic OCBA rule. The budget-adaptive property possessed by our proposed rule illustrates the significant impact of budget on the optimal budget allocation strategy \citep{peng2015non,shi2022dynamic}. Finally, based on the budget-adaptive rule, we develop two fully sequential algorithms, called final-budget anchorage allocation (FAA) and dynamic anchorage allocation (DAA). It is shown that both FAA and DAA are consistent, that is, as iteration increases, PCS converges to 1, and allocations to designs made by both algorithms converge to the asymptotic OCBA rule. In the numerical experiments, various problem settings are examined to illustrate the superior efficiency of both algorithms over benchmarks due to their budget-adaptive property. 

To summarize, the main contributions of this paper are as follows: 
\begin{itemize}
    \item[1.] We derive a budget-adaptive rule and explicitly demonstrate its distinct behavior under different budgets. Unlike the OCBA rule, which does not adjust to budget, our resulting rule depends on $T$. This budget-adaptive property is the major novelty of our proposed rule, as budget is often limited in practice. Our derivation is based on the setting where $k \to \infty$ and $T = \mathcal{O}(k \log k)$, which differs from the asymptotic OCBA rule that requires a stronger condition $k \to \infty$ and $T = \mathcal{\omega}(k \log k)$. 
    \item[2.] We develop two heuristic algorithms: FAA and DAA, that implement the proposed budget-adaptive rule based on different approaches. Additionally, we demonstrate that both algorithms will converge to the asymptotic OCBA rule. 
    \item[3.] We conduct comprehensive numerical experiments, including both small- and large-scale problems, to demonstrate the desirable budget-adaptive property of FAA and DAA. These experiments show that both algorithms significantly improve performance under small-budget conditions compared with the asymptotic OCBA rule.
\end{itemize}

\section{Related literature}
\label{related literature}
In R\&S, the objective is to identify the best design among a set of alternatives with respect to a performance metric, e.g., PCS. Due to the simulation noise, it is impossible to surely identify the best design within finite observations. Therefore, a strategy that intends to intelligently allocate simulation replications among designs is supposed to be developed. 

There are two branches of problem settings in R\&S literature. One is fixed-confidence setting and the other is fixed-budget setting. Fixed-confidence R\&S primarily focuses on the indifference zone (IZ) formulation and tries to guarantee a pre-specified level of PCS while using as little simulation budget as possible \citep{kim2001fully}. In later work, the IZ formulation is implemented to develop Frequentist procedures that can adapt to fully sequential setting \citep{hong2005tradeoff,batur2006fully,hong2007selecting}. Indifference-zone-free procedure that does not require an IZ parameter is proposed in \cite{fan2016indifference}. More recently, IZ procedures for large-scale R\&S problems in parallel computing environment are developed \citep{luo2015fully,zhong2022knockout,hong2022solving}. The fixed-budget R\&S procedures are designed to optimize a certain kind of performance metric by efficiently allocating a fixed simulation budget. In the fixed-budget setting, there are procedures that allocate a simulation budget according to an asymptotically optimal allocation rule, such as OCBA \citep{chen2000simulation}, the large deviation allocation  \citep{glynn2004large}, and the optimal expected opportunity cost allocation (OEA) \citep{gao2017new}; and procedures that myopically maximize the expected one-step-ahead improvement, such as the expected value of information (EVI) \citep{chick2010sequential}, the knowledge gradient \citep{frazier2008knowledge}, and the approximately optimal myopic allocation policy (AOMAP) \citep{peng2016myopic}. In particular, the approximately optimal allocation policy (AOAP) achieves both one-step-ahead optimality and asymptotic optimality \citep{peng2018ranking}. In most cases, fixed-budget R\&S procedures require less simulation budget than fixed-confidence R\&S procedures to achieve the same level of PCS due to their better adaptiveness to the observed simulation outputs, however, they can not provide a statistical guarantee as fixed-confidence R\&S procedures do.

There is a stream of literature in R\&S taking a global optimization perspective and focusing on the asymptotic behavior of allocation rules. The premise of these allocations is that if such allocations perform optimally when the simulation budget is sufficiently large, then they should also have satisfactory performances when the simulation budget is small. OCBA is such a typical method that allocates simulation budget according to an asymptotically optimal allocation rule when sampling distributions are normal \citep{chen2000simulation}. In later work, \cite{glynn2004large} applies the Large Deviation theory and extend the analyses to a more general setting where sampling distributions are non-Gaussian. More recently, \cite{chen2023balancing} follow this line and propose a new budget allocation rule that adaptively learns the optmality conditions obtained by using the Large Deviation theory.  \cite{gao2017new} present a budget allocation rule that uses the expected opportunity cost (EOC) as the quality measure of their procedure and is shown to be asymptotically optimal. \cite{peng2016dynamic} formulate the problem in a stochastic dynamic program framework and derive an approximately optimal design selection policy as well as an asymptotically optimal budget allocation policy. Furthermore, this stream of methods, which explore the asymptotic behavior of allocations, are extended to solving many variants of R\&S problem and their applications, such as the subset selection problem \citep{chen2008efficient,zhang2015simulation,gao2016new,zhang2022efficient}, ranking and selection with input uncertainty \citep{gao2017robust,xiao2018simulation,xiao2020optimal}, ranking and selection with multiple objectives \citep{lee2010finding}, stochastically constrained ranking and selection problem \citep{hunter2013optimal,pasupathy2014stochastically}, contextual ranking and selection problem \citep{li2022efficient}, efficient estimation of the risk measure problem \citep{wang2023efficient}, and preventive maintenance
optimization problem \citep{shi2021joint}. The most common simplification made by such methods is to consider an asymptotically large budget, which leads to solving a simplified budget allocation problem. However, this simplification results in derived allocation rules ignoring the impact of budget size on the budget allocation strategy.

While a huge number of works contribute to developing asymptotically optimal rules, few works investigate the impact of simulation budget on the budget allocation strategy. Typical myopic allocation rules \citep{frazier2008knowledge,chick2010sequential,ryzhov2016convergence,peng2018ranking,wang2023efficient} optimize one-step-ahead improvement. In particular, \cite{peng2017gradient} consider a low-confidence scenario and propose a gradient-based myopic allocation rule, which takes the induced correlations into account and performs well in such scenario. In later work, a myopic allocation rule that possesses both one-step-ahead optimality and asymptotic optimality is developed in \cite{peng2018ranking}. However, existing myopic allocation rules can not adapt to the simulation budget and some of them are not asymptotically optimal, even though they have excellent performances especially when the simulation budget is small. More recently, \cite{qin2022non} formulate the budget allocation problem as a dynamic program (DP) problem and develop a non-myopic knowledge gradient (KG) policy, which can look multiple steps ahead. \cite{shi2022dynamic} propose a dynamic budget-adaptive allocation rule for feasibility determination (FD) problem, a variant of R\&S problem, and show their allocation rule possesses both finite-budget properties and asymptotic optimality. \cite{cheng2023finite} propose an adaptive budget allocation rule that intends to deliver statistical guarantee on PCS while consuming as little budget as possible. None of existing works consider and develop a budget allocation rule that can not only adapt to the simulation budget but also achieve asymptotic optimality, for solving R\&S problems under a fixed budget setting.

The rest of the paper is organized as follows. In Section \ref{problem formulation}, we formally formulate the budget allocation problem. In Section \ref{optimal budget allocation strategy}, (1) we show the asymptotic OCBA rule can be recovered under the asymptotic regime where $k \to \infty$ and $T = \omega(k \log k)$; (2) we derive a budget-adaptive rule under the setting where $k \to \infty$ but $T = \mathcal{O}(k \log k)$ and explicitly present its distinct behaviour under different budgets; and (3) we develop two heuristic algorithms implementing the proposed budget-adaptive rule. In Section \ref{numerical experiments}, numerical experiments on synthetic examples and a case study are conducted. In the end, Section \ref{conclusion} concludes the paper.

\section{Problem formulation}
\label{problem formulation}
We formally introduce the following notations in our paper.
\begin{longtable}{ll}
\normalsize
$k$                     & Total number of designs;   \\
$\mathcal{K}$           & Set of designs, i.e., $\mathcal{K} = \{1,2,\dots,k\}$;   \\
$T$                     & Simulation budget;   \\
$X_{i,j}$               & The $j$-th simulation output sample of design $i$;   \\
$\mu_i$                 & Mean of the performance of design $i$, i.e., $ \mu_i = \mathbb{E}[X_{i,j}]$; \\
$\sigma^2_i$            & Variance of the performance of design $i$, i.e., $\sigma^2_i = \text{Var}[X_{i,j}]$;  \\
$b$                     & Real best design, i.e., $b = \arg \min_{i \in \mathcal{K}} \mu_i$;   \\
$\mathcal{K}^{\prime}$  & Set of non-best designs, i.e., $\mathcal{K}^{\prime} = \mathcal{K} \backslash b$;   \\
$N_i$                   & The number of simulation replications allocated to design $i$;   \\
$w_i$                   & Proportion of simulation budget allocated to design $i$, i.e., $w_i = N_i/T$;   \\
$\hat{\mu}_i$           & Sample mean of the performance of design $i$, i.e., $\hat{\mu}_i = (1/N_i) \sum_{j=1}^{N_i} X_{i,j}$;    \\
$\hat{b}$               & Observed best design, i.e., $\hat{b} = \arg \min_{i \in \mathcal{K}} \hat{\mu}_i$.  
\end{longtable}

Suppose that there are $k \geq 3$ designs in contention. For each design $i \in \mathcal{K}$, its mean performance $\mu_i$ is unknown and can only be estimated by sampling replications via a stochastic simulation model. The goal of R\&S is to identify the real best design $b$, where “best” is defined with respective to the smallest mean. Assume the best design $b$ is unique, i.e., $\mu_b < \mu_i$, for $ i \in \mathcal{K}^{\prime}$. This assumption basically requires the best design can be distinguished from the others. Common random numbers and correlated sampling are not considered in the paper, and we assume the simulation output samples are independent across different designs and replications, i.e., $X_{i,j}$ is independent for all $i$ and $j$. The most common assumption on the sampling distribution is that the simulation observations of each design $i$ are i.i.d. normally distributed with mean $\mu_i$ and variance $\sigma_i^2$, i.e., $X_{i,j} \sim N(\mu_i,\sigma_i^2)$, for $i \in \mathcal{K}$ and $j \in \mathbb{Z}^+$. In our analysis, we focus on the normal sampling distribution. For non-Gaussian distributions, batches of independent samples of each design are considered as a single sample.  Though batching does not impact the convergence rate of PCS from the large deviation perspective \citep{glynn2004large}, using batching and assuming normality is implemented for practical convenience. 

After the simulation budget $T$ is depleted, the observed best design $\hat{b}$ (with the smallest sample mean) is selected. The event of correct selection occurs when the selected design, design $\hat{b}$, is the real best design, design $b$. Thus, we define the probability of correct selection (PCS) as
\begin{align*}
    \text{PCS} &= \Pr ( \hat{b}  = b  )\\
    &= \Pr \left( \bigcap_{i \in \mathcal{K}^{\prime}} \left\{ \hat{\mu}_b < \hat{\mu}_{i} \right\}  \right).
\end{align*}
The problem of interest is to determine $N_1,N_2,\dots,N_k$, such that by the time the simulation budget is exhausted and we select the observed best design, PCS is maximized. We model the budget allocation problem as follows:  
\begin{equation*}
\begin{split}
    Problem \ \mathcal{P}: \  &\max \ \text{PCS}   \\
    &\text{s.t.} \ \sum\limits_{i \in \mathcal{K}} N_i = T,\\
    & \quad \quad N_i \geq 0, \quad  i \in \mathcal{K}.
\end{split}
\end{equation*}
For simplicity, we ignore the integer constraints on $N_i = w_i T$, for $i \in \mathcal{K}$, and we assume that $N_i$ are continuous variables (see \citet[p.71]{chen2011stochastic} for a comprehensive discussion). In practical implementations, the simulation replication numbers $N_i$, for $i \in \mathcal{K}$, derived by solving Problem $\mathcal{P}$ can be rounded up to the largest integer smaller than $N_i$, i.e., $ \lfloor N_i \rfloor $, where $ \lfloor \cdot \rfloor $ is the flooring function. 

Under general settings, the major difficulty in solving Problem $\mathcal{P}$ is that there is no closed-form expression for PCS. To evaluate PCS, one can use Monte Carlo simulation, but its computational cost is usually unaffordable, especially when the simulated systems are of high complexity. More recently, \cite{eckman2022posterior} propose to use the integral function provided by MATLAB to calculate the posterior PCS under a Bayesian setting. However, computational burden of this method is still of concern if one intends to calculate the posterior PCS repeatedly during the allocation procedure. In this paper, we follow the path of OCBA method and approximate PCS by a computationally cheap lower bound. One of the most commonly used lower bounds of PCS can be derived by applying the Bonferroni inequality \citep{galambos1977bonferroni,chen2000simulation,zhang2015simulation},
\begingroup
\allowdisplaybreaks
\begin{align}
\label{Bonferroni inequality}
    \text{PCS} &= \Pr \left( \bigcap_{i \in \mathcal{K}^{\prime}} \left\{ \hat{\mu}_b < \hat{\mu}_{i} \right\}  \right) \notag \\
    & \geq 1 - \sum_{i \in \mathcal{K}^{\prime}} \Pr \left( \hat{\mu}_i \leq \hat{\mu}_{b} \right) \notag \\
    &= 1 - \sum_{i \in \mathcal{K}^{\prime}}  \Pr \left( Z \leq - \frac{\mu_i -  \mu_b}{ \sqrt{ \sigma_i^2/N_i + \sigma_b^2/N_b } }  \right)  \notag \\
    &= 1 - \sum_{i \in \mathcal{K}^{\prime}}  \Phi \left( - \frac{ \delta_{i,b} }{ \sigma_{i,b} }  \right) = \text{APCS},
\end{align}
\endgroup
where $Z$ is a random variable follows the standard normal distribution, $\Phi$ denotes the cumulative distribution function (c.d.f.) of the standard normal random variable, $\delta_{i,b} = \mu_i -  \mu_b $, and $\sigma_{i,b} = \sqrt{\sigma_i^2/N_i + \sigma_b^2/N_b} $. Another commonly used lower bound of PCS can be derived by applying the Slepian's inequality \citep{slepian1962one}. However, the OCBA framework requires the objective to be concave such that convex programming theory applies, nonetheless, it can be checked that the Slepian bound is neither convex nor concave. The tightness of the APCS bound in \eqref{Bonferroni inequality} grows with $T$. And when $T$ is sufficiently large, APCS is tight as it will converge to 1. The APCS's closed-form expression enables us to efficiently evaluate PCS, and its concavity (as will be shown in Lemma 1) makes it possible for us to characterize its optimality conditions in the context of a general budget. Therefore, instead of solving Problem $\mathcal{P}$, with new objective $\text{APCS}$, we consider the following optimization problem:
\begin{equation*}
\begin{split}
    Problem \ \mathcal{P}1: \  &\max \ \text{APCS}   \\
    &\text{s.t.} \ \sum\limits_{i \in \mathcal{K}} w_i = 1,\\
    & \quad \quad w_i \geq 0, \quad  i \in \mathcal{K}.
\end{split}
\end{equation*}

\section{Budget allocation strategy}
\label{optimal budget allocation strategy}
In Section \ref{subsection: optimal computing budget allocation}, we recover the asymptotic OCBA rule in \cite{chen2000simulation}. In Section \ref{subsection: budget-adaptive allocation rule}, we propose a budget-adaptive rule and analyze how it differs from the asymptotic OCBA rule due to its budget-adaptive property. Based on the proposed budget-adaptive rule, two heuristic algorithms are developed in Section \ref{subsection: budget allocation procedure}. 

Of note, we highlight the most important implication of this section: simulation budget significantly impacts the budget allocation strategy. To enhance readability, all the proofs, and the definitions of   Big-$\mathcal{O}$, Big-$\Theta$, and Little-$\omega$ notations are relegated to the Online Appendix.

\subsection{Optimal computing budget allocation}
\label{subsection: optimal computing budget allocation}
In this subsection, we first extend some asymptotic results in the development of the asymptotic OCBA rule to the case of a general budget. Then, we show the asymptotic OCBA rule can be recovered under the asymptotic regime where (1) $k \to \infty$; and (2) $T = \omega(k\log k)$.

In the OCBA paradigm, the derivation of optimality conditions essentially requires Problem $\mathcal{P}$1 to be a convex optimization problem. \cite{zhang2015simulation} show APCS is concave when $T \to \infty$. We generalize this result and rigorously show in Lemma 1 that the concavity of APCS in Problem $\mathcal{P}1$ indeed holds for any $T \in \mathbb{Z}^+$, thereby establishing Problem $\mathcal{P}$1 as a convex optimization problem.

\begin{Lemma}
\label{lemma: convexity of APCS}
    APCS is concave and therefore Problem $\mathcal{P}1$ is a convex optimization problem.
\end{Lemma}

\textit{Proof:} See Online Appendix B.1. 

With Lemma 1, the theory of convex programming applies, and the solution that satisfies the Karush-Kuhn-Tucker (KKT) conditions is optimal to Problem $\mathcal{P}1$ \citep{boyd2004convex}. According to \cite{chen2000simulation}, we give in Lemma 2 without proof the optimality conditions of Problem $\mathcal{P}1$.

\begin{Lemma}
\label{lemma: optimality conditions}
If the solution $w = (w_1,w_2,\dots,w_k)$ maximizes the $\text{APCS}$ in Problem $\mathcal{P}1$, it satisfies the following optimality conditions $\mathcal{C}_1,\mathcal{C}_2,\mathcal{C}_3$, and $\mathcal{C}_4$
\begin{itemize}
    \item[$\mathcal{C}_1$:] $w_b = \sigma_b \sqrt{\sum_{i \in \mathcal{K}^{\prime}} w_i^{2}/\sigma_i^2}$,
    \item[$\mathcal{C}_2$:] $  - \frac{\delta_{i,b}^2}{ 2(\sigma_i^2/w_i + \sigma_b^2/w_b) } T  +  \log \frac{\delta_{i,b}\sigma_i^2}{(\sigma_i^2/w_i + \sigma_b^2/w_b)^{\frac{3}{2}}} - 2 \log w_i = \lambda, \quad  i \in \mathcal{K}^{\prime}$,
    \item[$\mathcal{C}_3$:] $\sum_{i \in \mathcal{K}} w_{i}=1$,
    \item[$\mathcal{C}_4$:] $w_i \geq 0, \quad  i \in \mathcal{K}$,
\end{itemize}
where $\lambda$ is a constant.
\end{Lemma}

\begin{Remark}
\normalfont{Another approximation for PCS, which resembles the Bonferroni bound in \eqref{Bonferroni inequality}, is introduced in \cite{gao2017new},
\begin{equation}
\label{Bonferroni-like bound}
    \text{PCS} \geq 1 - \sum_{i \in \mathcal{K}^\prime} \exp \left( -\frac{T \delta_{i,b}^2}{2( \sigma_i^2/w_i + \sigma_b^2/w_b )} \right) = \text{APCS}_{0}.
\end{equation}
Compared with the APCS in \eqref{Bonferroni inequality}, APCS$_{0}$ approximates the probability of each pairwise incorrect selection $\Pr(\hat{\mu}_i \leq \hat{\mu}_b)$ by $\exp ( - T \delta_{i,b}^2/2( \sigma_i^2/w_i + \sigma_b^2/w_b ))$. APCS and APCS$_{0}$ are different approximations of PCS, even though they have similar forms.}
\end{Remark}

Notice that, $\mathcal{C}_1$ and $\mathcal{C}_2$ are highly nonlinear, and determining the exact solution to Problem $\mathcal{P}1$ requires a numerical solver. Since sampling efficiency is of significant concern in R\&S, one would expect to derive a solution in analytical form. To do this, one can make mild modifications to optimality conditions, i.e., letting $w_b \gg w_i$ for $i \in \mathcal{K}^{\prime}$, and assume a large enough budget, i.e., letting $T \to \infty$, to simplify computations. This is the basic idea behind the derivation of the asymptotic OCBA rule in \cite{chen2000simulation}. \cite{zhang2015simulation} show that when both $T \to \infty$ and $k \to \infty$, the approximation $w_b \gg w_i$ can be justified. However, when $T$ is small, this result does not straightforwardly apply.  Proposition \ref{proposition:verify w_b is far larger than w_i} formally generalizes the analyses in \cite{zhang2015simulation} to a general budget, i.e., for any $T \in \mathbb{Z}^+$, if $k \rightarrow \infty$, we have $w_b \gg w_i$, for $i \in \mathcal{K}^{\prime}$. To prove Proposition 1, we need the following assumption, which requires the problem we aim to solve is neither too hard nor too easy.

\begin{Assumption}
\label{assumption: problem structure}
    The variances of all designs are lower bounded by $\underline{\sigma}^2 > 0$ and upper bounded by $\bar{\sigma}^2 > 0$, i.e., $0 < \underline{\sigma}^2 \leq \sigma^2_i \leq \bar{\sigma}^2$, for $i \in \mathcal{K}$, and the mean differences between the best design and any other designs are lower bounded by $\underline{\delta} > 0$ and upper bounded by $\bar{\delta} > 0$, i.e., $0 < \underline{\delta} \leq \delta_{i,b} \leq \bar{\delta}$, for $i \in \mathcal{K}^{\prime}$.
\end{Assumption}

\begin{Proposition}
\label{proposition:verify w_b is far larger than w_i}
Let Assumption \ref{assumption: problem structure} always hold. Suppose that as $k$ increases, the identity of the best design is fixed, and newly added designs are sub-optimal. If the solution $w = (w_1,w_2,\dots,w_k)$ satisfies Lemma 2, for any $T\in \mathbb{Z}^+$, we have
\begin{equation*}
    \frac{w_i}{w_b} = \mathcal{O} \left( \frac{1}{\sqrt{k}} \right), \quad i \in \mathcal{K}^{\prime}.
\end{equation*}
Therefore, if $k \rightarrow \infty$, we have $w_b \gg w_i$, for $i \in \mathcal{K}^{\prime}$ and $T \in \mathbb{Z}^+$.
\end{Proposition}

\textit{Proof:} See Online Appendix C.1. 

We give an example to illustrate Proposition \ref{proposition:verify w_b is far larger than w_i}. Consider a slippage configuration of means, i.e., $\mu_i - \mu_b = \delta$ for $i  \in \mathcal{K}^\prime$ for some $\delta > 0$, with a common variance $\sigma^2$, i.e., $ \sigma_i^2 = \sigma^2$ for $i \in \mathcal{K}$. For this problem structure, it can be checked that $w$ satisfies $w_b = 1/(\sqrt{k-1} + 1)$ and $w_i = 1/(k - 1 + \sqrt{k - 1)}$, for $i \in \mathcal{K}^\prime$. Then, we have $w_i/w_b = 1/\sqrt{k-1} = \mathcal{O}(1/\sqrt{k})$. As $k \rightarrow \infty$, $w_i/w_b \rightarrow 0$, and therefore $w_b \gg w_i$, for $i \in \mathcal{K}^\prime$. In practical implementations, the problem we face may involve a large number of alternative designs, i.e., $k$ can be very large. For example, in a scheduling problem \citep{hong2021review}, the manager wants to simultaneously determine the jobs to be scheduled, the values assigned to the jobs, and the time when the scheduling happens. If there are 50 choices for each component, the total number of alternative designs will be over $10^5$. More recently, \cite{hong2022solving} has considered and managed to solve problems with over $10^6$ alternative designs. With Proposition \ref{proposition:verify w_b is far larger than w_i}, the optimality conditions can be greatly simplified so that it is possible to derive a solution in analytical form.

We recover the asymptotic OCBA rule. First, we consider a large-scale problem, i.e., $k \to \infty$, such that $w_b \gg w_i$, for $i \in \mathcal{K}^{\prime}$ (because Proposition 1 applies), and $\mathcal{C}_2$ in Lemma \ref{lemma: optimality conditions} becomes
\begin{equation}
\label{simplified condition C2-1}
     \log I_i + \log w_i + \frac{T}{I_i} w_i = \lambda, \quad  i \in \mathcal{K}^{\prime},
\end{equation}
where $I_i = \sigma_i^2/\delta_{i,b}^2$, for $i \in \mathcal{K}^\prime$. For notational simplicity, we further let $I_b = \sigma_b \sqrt{\sum_{i \in \mathcal{K}^{\prime}}I_i^{2} /\sigma_i^2}$. The following proposition shows if $T = \omega(k\log k)$, the two terms $\log I_i$ and $\log w_i$ in \eqref{simplified condition C2-1} are negligibly small to $Tw_i/I_i$.

\begin{Proposition}
\label{proposition: conditions on ignoring terms}
Let Assumption \ref{assumption: problem structure} always hold. Suppose that as $k$ increases, the identity of the best design is fixed, and newly added designs are sub-optimal. If the solution $w = (w_1,w_2,\dots,w_k)$ satisfies $\mathcal{C}_1$, $\mathcal{C}_3$, $\mathcal{C}_4$, and \eqref{simplified condition C2-1}, we have $(\log I_i) / (w_i/I_i) = \Theta(k)$ and $(\log w_i) / (w_i/I_i) = \Theta(k \log k)$. Therefore, as $k \to \infty$, if $T = \omega(k \log k)$, both $\log I_i$ and $\log w_i$ are negligibly small compared with $Tw_i/I_i$, i.e., $(\log I_i) / (Tw_i/I_i) \to 0$ and $(\log w_i) / (Tw_i/I_i) \to 0$ for $i \in \mathcal{K}^\prime$. 
\end{Proposition}

\textit{Proof:} See Online Appendix C.2.

Second, with Proposition \ref{proposition: conditions on ignoring terms}, we let $T = \omega(k \log k)$ such that both $\log I_i$ and $\log w_i$ can be ignored relative to $Tw_i/I_i$. By doing so, we obtain the asymptotic OCBA rule
\begin{equation}
\label{OCBA allocation rule}
    w^*_i = \frac{I_i}{\sum_{i \in \mathcal{K}} I_i}, \quad i \in \mathcal{K},
\end{equation}
Clearly, the asymptotic OCBA rule can be recovered by considering a large-scale problem (i.e., $k \to \infty$) with a specific large budget setting (i.e., $T = \omega(k \log k)$). 

\begin{Remark}
\normalfont{If one considers a small-scale problem, i.e., if  $k$ is small (so that Proposition 1 does not apply), by simply letting $T \to \infty$, it can be checked that $w$ satisfies
\begingroup
\allowdisplaybreaks
\begin{equation}
\label{large deviation ratio}
\begin{split}
    w_b &= \sigma_b \sqrt{\sum\limits_{i \in \mathcal{K}^{\prime}} \frac{w_i^{2}}{\sigma_i^2}}, \\
    \frac{\delta_{i,b}^2}{(\sigma_i^2/w_i + \sigma_b^2/w_b) }  &= \frac{\delta_{j,b}^2}{(\sigma_j^2/w_j + \sigma_b^2/w_b) }, \ i,j \in \mathcal{K}^{\prime} \ \text{and} \  i \neq j ,
\end{split}
\end{equation}
\endgroup
which corresponds to the optimality conditions in \cite{glynn2004large} derived by maximizing the asymptotic convergence rate of probability of false selection (PFS) under normal sampling distributions. Equation \eqref{large deviation ratio} is a great simplification but still requires a numerical solver to determine $w$.}
\end{Remark}

The asymptotic OCBA rule $w^*$ tends to assign high budget allocation ratios to non-best designs with large $I_i$, while it tends to assign low budget allocation ratios to those with small $I_i$. Since $I_i = \sigma_i^2/(\mu_b-\mu_i)^2$ are the inverse signal-to-noise ratios, the behavior of $w^*$ indicates that more simulation budget should be allocated to non-best designs that are hard to distinguish from the best, while less simulation budget should be allocated to those that are easy to distinguish.

In practical implementations, the asymptotic OCBA rule is attractive  due to its closed-form expression and impressive performance. However, it ignores the impact of budget on it. This observation motivates the need to derive a desirable allocation rule which is adaptive to the simulation budget.

\subsection{Budget-adaptive allocation rule}
\label{subsection: budget-adaptive allocation rule}
In this subsection, we develop a budget-adaptive rule that incorporates the budget and analyze its distinct behavior under different budgets. Our derivations consider the setting where $k \to \infty$ but $T = \mathcal{O}(k \log k)$. Under such setting, $T$ is not large enough to simplify computations. 
To derive a closed-form solution, we use the first-order Taylor series expansion at the asymptotic OCBA rule to linearize a set of optimality conditions. Doing so not only keeps $T$ from being ignored, but also enables us to understand how $T$ impacts the allocation strategy. 

First, we let $w_b \gg w_i$, for $i \in \mathcal{K}^{\prime}$ (because Proposition \ref{proposition:verify w_b is far larger than w_i} applies), and condition $\mathcal{C}_2$ in Lemma \ref{lemma: optimality conditions} becomes
\begin{equation}
\label{simplified condition C2-1 (section 4.2)}
     \log I_i + \log w_i + \frac{T}{I_i} w_i = \lambda, \quad  i \in \mathcal{K}^{\prime}.
\end{equation}
According to Proposition \ref{proposition: conditions on ignoring terms}, when $T = \mathcal{O}(k \log k)$, the terms $\log I_i$ or $\log w_i$ for $i \in \mathcal{K}^\prime$ cannot be ignored to simplify computations. Notice that, even with $w_b \gg w_i$, conditions $\mathcal{C}_1$, $\mathcal{C}_3$, $\mathcal{C}_4$, and \eqref{simplified condition C2-1 (section 4.2)} are still a system of highly nonlinear equations and require a numerical solver to determine $w$. To derive a solution in analytical form, we further approximate the term $\log w_i$ in \eqref{simplified condition C2-1 (section 4.2)} by its first-order Taylor series expansion at point $w^*_i$
\begin{equation}
\label{first order taylor expansion}
    \log w_i \approx \log w^*_i + \left(w_i - w^*_i \right)/w^*_i, \quad  i \in \mathcal{K}^{\prime}.
\end{equation}
In \eqref{first order taylor expansion}, the asymptotically optimal solution $w^* = (w^*_1,w^*_2,\allowbreak\dots,w^*_k)$ is regarded as a “good” approximation of the optimal solution $w = (w_1,w_2,\allowbreak\dots,w_k)$ to Problem $\mathcal{P}1$. This approximation tends to be accurate when $T$ increases so that the gap between $w^*$ and $w$ would decrease. Then, we substitute the term $\log w_i$ with its approximation provided by \eqref{first order taylor expansion} and obtain the approximated optimality conditions for Problem $\mathcal{P}1$:
\begin{itemize}
    \item[$\mathcal{C}_1$:] $w_b = \sigma_b \sqrt{\sum_{i \in \mathcal{K}^{\prime}} w_i^{2}/\sigma_i^2}$,
    \item[$\widehat{\mathcal{C}}_2$:] $  2 \log I_i + \left(\frac{T}{I_i} + \frac{1}{w^*_i} \right) w_i = \lambda, \quad  i \in \mathcal{K}^{\prime}$,
    \item[$\mathcal{C}_3$:] $\sum_{i \in \mathcal{K}} w_{i}=1$,
    \item[$\mathcal{C}_4$:] $w_i \geq 0, \quad i \in \mathcal{K}$,
\end{itemize}
where $\lambda$ is a constant. In contrast to the sufficiently large budget approximation (i.e., letting $T = \omega(k \log k)$), the first-order Taylor series expansion used in \eqref{first order taylor expansion} keeps more features (including $T$) of the
optimality conditions from being ignored, and it yields $\widehat{\mathcal{C}}_2$ that incorporates $T$. Additionally, compared with $\mathcal{C}_2$, condition $\widehat{\mathcal{C}}_2$ is a linear equation of $w_i$, and therefore solving $\mathcal{C}_1$, $\widehat{\mathcal{C}}_2$, $\mathcal{C}_3$, and $\mathcal{C}_4$ leads to a solution in analytical form. We temporarily omit the non-negativity constraints (i.e., $\mathcal{C}_4$) and consider $\mathcal{C}_1$, $\widehat{\mathcal{C}}_2$, and $\mathcal{C}_3$ in Lemma \ref{lemma: budget-adaptive allocation rule}. In Lemma \ref{lemma: feasibility of allocation rule}, the non-negativity constraints are discussed to guarantee the feasibility of the solution obtained in Lemma \ref{lemma: budget-adaptive allocation rule}.

\begin{Lemma}
\label{lemma: budget-adaptive allocation rule}
If the solution $W(T) = (W_1(T),W_2(T),\dots,\allowbreak W_k(T))$ solves conditions $\mathcal{C}_1$, $\widehat{\mathcal{C}}_2$, and $\mathcal{C}_3$, it satisfies
\begin{equation}
    W_i(T) =\left\{\begin{array}{ll}
w^*_i \alpha_i(T) & \quad \text{if\ } i \in \mathcal{K}^{\prime} \\
\sigma_b \sqrt{\sum_{i \in \mathcal{K}^{\prime}}\frac{(W_i(T))^{2} }{\sigma_i^2}} & \quad \text{if\ } i = b
\end{array}\right.
\end{equation}
\textit{where}
\begingroup
\allowdisplaybreaks
\begin{align}
	\alpha_i(T) &= \frac{ (\lambda - 2 \log I_i)}{1 + T/S}, \notag \\
    S &= \sum_{i \in \mathcal{K}} I_i, \notag \\
    \lambda &=\left\{\begin{array}{ll}
\frac{- q + \sqrt{q^2 - 4 p r}}{2 p}  & \quad \text{if} \  w^*_b \neq \frac{1}{2}  \\
\frac{4\sum_{i \in \mathcal{K}^{\prime}} I_i \log I_i + T + S}{2\sum_{i \in \mathcal{K}^{\prime}} I_i} & \quad \text{if} \ w^*_b = \frac{1}{2} \\
\end{array}\right. \notag \\
    p  &=  S(2 I_b - S), \notag \\
    q  &= - 4 \sigma_b^2  \sum\limits_{i \in \mathcal{K}^{\prime}} \frac{I_i^2 \log I_i }{\sigma_i^2 } \notag + 2 (S- I_b)  \left(2 \sum\limits_{i \in \mathcal{K}^{\prime}}  I_i \log I_i  + T+S \right), \notag\\
    r  &=  4 \sigma_b^2  \sum\limits_{i \in \mathcal{K}^{\prime}} \frac{I_i^2 \log^2 I_i }{\sigma_i^2 } -  \left( 2 \sum\limits_{i \in \mathcal{K}^{\prime}}  I_i \log I_i + T+S \right)^2. \notag
\end{align}
\endgroup
\end{Lemma}

\textit{Proof:} See Online Appendix B.2. 

The solution $W(T)$ is an analytical function of $T$, and it reduces to the asymptotic OCBA rule $w^*$ when $T = \omega(k \log k)$. For a certain budget $T$, a non-best design $i$, for $i \in \mathcal{K}^{\prime}$, tends to be allocated more  budget by $W_i(T)$ than by $w^*_i$ if $\alpha_i(T) > 1$ and be allocated less budget by $W_i(T)$ than by $w^*_i$ if $\alpha_i(T) < 1$. The “balance” between the proportional allocations to the best design and non-best designs remains unchanged. In particular, when $T = \omega(k \log k)$, we have $\alpha_i(T) \rightarrow 1 $, then $W_i(T) \rightarrow w_i^*$, for $i \in \mathcal{K}^{\prime}$, and consequently, $W_b(T) \rightarrow w_b^*$. This implies that $W(T)$ achieves asymptotic optimality, and it is always feasible. However, some of $W_i(T)$, for $i \in \mathcal{K}^\prime$, defined in Lemma 3 may violate the non-negativity constraints when $T = \mathcal{O}(k \log k)$. Because $W(T)$ is derived by temporarily omitting the non-negativity constraints. For non-best designs $i \in \mathcal{K}^{\prime}$, let $\langle j \rangle$, $j = 1,2,\dots,k - 1$, be the ascending order statistics of $I_i$, i.e., $I_{\langle 1 \rangle} \leq I_{\langle 2 \rangle} \leq \dots \leq I_{\langle k-1 \rangle}$. Lemma 4 gives a sufficient condition for $W(T)$ being feasible. 

\begin{Lemma}
\label{lemma: feasibility of allocation rule}
Suppose that the solution $W(T) = (W_1(T),\allowbreak W_2(T),\dots, W_k(T))$ solves conditions $\mathcal{C}_1$, $\widehat{\mathcal{C}}_2$, and $\mathcal{C}_3$. If $T \geq T_0$, $W(T)$ is always feasible, i.e., $W_i(T) \geq 0$, for $i \in \mathcal{K}$, where
\begingroup
\allowdisplaybreaks
\begin{equation}
\label{T0}
    T_0 = 
    \left\{\begin{array}{ll}
\max\{0,T_1,T_2\}  & \text{if} \  w^*_b \neq \frac{1}{2}  \\
\max \{ 0, 4 \sum_{i \in \mathcal{K}^{\prime}} I_i \log\frac{I_{\langle k-1 \rangle}}{I_i} - S\} &  \text{if} \ w^*_b = \frac{1}{2} \\
\end{array}\right.
\end{equation}
and
\begin{align*}
    T_1 &= 2 \sum\limits_{i \in \mathcal{K}^{\prime}} \left[\frac{\sigma^2_b I_i^2}{\sigma_i^2 (S - I_b)} - I_i \right] \log \frac{I_{\langle k - 1 \rangle}}{I_i}  - S,\\
    T_2 &= 2 \sum\limits_{i \in \mathcal{K}^{\prime}}  I_i \log \frac{I_{\langle k - 1 \rangle}}{I_i} + 2 \sigma_b \sqrt{  \sum\limits_{i \in \mathcal{K}^{\prime}} \frac{  I_i^2 }{\sigma_i^2 } \left(  \log \frac{I_{\langle k - 1 \rangle}}{I_i}   \right)^2 } - S.
\end{align*}
\endgroup
\end{Lemma}

\textit{Proof:} See Online Appendix B.3. 

\begin{Remark}
\normalfont{Notice that the value of $T_0$ depends on the variance of designs and differences in means (between the best design and non-best designs). Additionally, the calculation of $T_0$ requires true means and variances of designs, which can not be known in practical implementations and must be estimated from samples (e.g., by plugging in sample means and variances). Therefore, besides the problem structure, the value of $T_0$ also depends on the allocation status.}
\end{Remark}

Lemma 4 shows that $W(T)$ is always feasible when $T \geq T_0$. However, when $T < T_0$, there exists a factor $\alpha_i(T)$ could become negative such that $W_i(T) = w_i^* \alpha_i(T) < 0$. This implies that $W_i(T)$ could be discounted too heavily to be feasible due to the effect of $\alpha_i(T)$. To address this issue, when $T < T_0$, we allocate $T$ simulation budget according to $W(T_0)$, which is always a feasible solution to Problem $\mathcal{P}1$. Let $\lceil T_0 \rceil$ denote the smallest integer that is larger than or equal to $T_0$. For non-best designs $i \in \mathcal{K}^{\prime}$, we define
\begin{equation}
\label{new allocation ratio for design i}
    \widetilde{W}_i(T) =\left\{\begin{array}{ll}
W_i(T)  & \quad \text{if\ }   T \geq T_0  \\
W_i(\lceil T_0 \rceil) & \quad \text{if\ } T \leq  T_0 \\
\end{array}\right.
\end{equation}
and for the best design  $b$, we define
\begin{equation}
\label{new allocation ratio for design b}
    \widetilde{W}_b(T) = \sigma_b \sqrt{\sum_{i \in \mathcal{K}^{\prime}}\frac{( \widetilde{W}_i(T))^{2} }{\sigma_i^2}}.
\end{equation}

\begin{Theorem}
\label{theorem: always feasible solution}
Let Assumption 1 always hold. Suppose that as $k$ increases, the identity of the best design is fixed, and newly added designs are sub-optimal. When $k \to \infty$, the solution $\widetilde{W}(T) = (\widetilde{W}_1(T),\widetilde{W}_2(T),\dots,\allowbreak\widetilde{W}_k(T))$ defined in \eqref{new allocation ratio for design i} and \eqref{new allocation ratio for design b} solves Problem $\mathcal{P}1$ and approximately maximizes the $\text{APCS}$.
\end{Theorem}

The budget-adaptive property and closed-form expression of $\widetilde{W}(T)$ not only greatly facilitate its implementation, but generates insights on how $T$ impacts the allocation strategy for identifying the best design. To understand how $T$ influences the allocation strategy, we need Proposition \ref{proposition: properties of alpha}. This proposition demonstrates the behavior of $\alpha_i(T)$, which plays key roles in influencing the budget allocation strategy.

\begin{Proposition}
\label{proposition: properties of alpha}
For $ T \in \mathbb{Z}^+$, we have $\alpha_{\langle 1 \rangle}(T) \geq 1 $, $\alpha_{\langle k - 1 \rangle}(T) \leq 1$, and  $\alpha_{\langle 1 \rangle}(T) \geq \alpha_{\langle 2 \rangle}(T) \geq \dots \geq \alpha_{\langle k - 1 \rangle}(T)$, where all the equalities hold if and only if $I_i$, for $i \in \mathcal{K}^{\prime}$, are all equal.
\end{Proposition}

\textit{Proof:} See Online Appendix C.3.

Notice that a large $I_i$ indicates that non-best design $i$ is hard to distinguish from the best design, and it corresponds to a small $\alpha_i(T)$. Compared with the asymptotic OCBA rule $w^*$, $\widetilde{W}(T)$ discounts the budget allocation ratios of non-best designs with large $I_i$ (e.g., design $\langle k - 1 \rangle$) and increases the budget allocation ratios of non-best designs with small $I_i$ (e.g., design $\langle 1 \rangle$) due to the effects of $\alpha_i(T)$, for $i \in \mathcal{K}^{\prime}$. These adjustments to budget allocation ratios are based on $T$. More specifically, relative to $w^*$, $\widetilde{W}(T)$ discounts the budget allocated to non-best designs that are hard to distinguish from the best, while it increases the budget allocated to those that are easy to distinguish. In particular, if all non-best designs are equally hard or easy to distinguish from the best, $\widetilde{W}(T)$ is identical to $w^*$. As for the best design, condition $\mathcal{C}_1$ and \eqref{new allocation ratio for design b} show that there exists a general “balance” between the proportions of total budget allocated to the best design and non-best designs. And the same “balance” exists in both $\widetilde{W}(T)$ and $w^*$. Additionally, we have
\begin{equation*}
    \frac{\widetilde{W}_b(T)}{w^*_b} = \sqrt{\frac{ \sum_{i \in \mathcal{K}^\prime} \widetilde{W}^2_i(T) / \sigma_i^2}{  \sum_{i \in \mathcal{K}^\prime} (w_i^*)^2 / \sigma_i^2 }}.
\end{equation*}
Besides $T$, whether the value of $\widetilde{W}_b(T)/w^*_b$ is larger or smaller than 1 also depends on the specific problem structure and allocation status. Therefore, how $T$ impacts $\widetilde{W}_b(T)$, especially whether it is larger or smaller than $w_b^*$, is in a complex manner and not involved in this discussion.

Based on the preceding analyses, it is clear that $\widetilde{W}(T)$ and $w^*$ admit different allocations when $T = \mathcal{O}(k \log k)$, even though they are equivalent when $T = \omega(k \log k)$. One major contribution of this paper is that we develop a budget-adaptive rule, and our derivations do not let $T$ be sufficiently large, i.e., let $T= \omega(k\log k)$, to simplify computations. In practical implementations, the budget-adaptive property is significant as simulation budget is often limited by its high-expense.

\subsection{Budget allocation algorithms}
\label{subsection: budget allocation procedure}
In this subsection, we develop two heuristic algorithms based on two approaches implementing the proposed budget-adaptive rule. Without loss of generality, we consider a fully sequential setting, i.e., only one sample is allocated in an iteration. To facilitate presentation, we introduce some additional notations.

\begin{longtable}{ll}
\normalsize
$t$                     & Consumed budget (the number of samples allocated);         \\
$A_{t}$                 & Design received the $t$-th sample;    \\
$\hat{b}^{(t)}$         & Design with the smallest sample mean with consumed budget $t$; \\
$\hat{\mu}^{(t)}_i$     & Sample mean of design $i$ with consumed budget $t$;     \\
$(\hat{\sigma}^{(t)}_i)^2$  & Sample variance of design $i$ with consumed budget $t$;       \\
$N^{(t)}_i$             & The number of samples received by design $i$ with consumed budget $t$; \\
$w^{(t)}_i$             & Proportion of samples allocated to design $i$ with consumed budget $t$; \\
$\hat{T}_0^{(t)}$ & Estimated value of $T_0$ with consumed budget $t$; \\
$\hat{w}^{*,(t)}_i$ & Estimated OCBA ratio of design $i$ with consumed budget $t$; \\
$\widetilde{W}^{(t)}_i(T)$ & Estimated budget-adaptive ratio of design $i$ with total budget $T$ and consumed budget $t$   \\
&(for notation simplicity, we omit the “hat” on $\widetilde{W}^{(t)}_i(T)$). \\
\end{longtable}
To calculate the budget allocation ratios, we use every design's sample mean and sample variance as plug-in estimates for its true mean and true variance, respectively. \cite{chen1996lower} and \cite{chick2001new} describe the main superiority of fully sequential procedures is that it can improve each stage's sampling efficiency by incorporating information from all earlier stages.

\subsubsection{Final-budget anchorage allocation}

We develop a fully sequential algorithm called final-budget anchorage allocation (FAA). At the beginning, the total budge $T$ is specified, and each design is sampled $n_0$ replications to get initial sample estimates. In each iteration, we run additional one replication according to $\widetilde{W}^{(t)}(T)$; update sample estimates; and repeat till the budget is exhausted. During the allocation procedure, the final budget $T$ is anchored, and the goal is to maximize the PCS after $T$ budget is depleted. The “most starving” technique introduced in \cite{chen2011stochastic} can be applied to define an allocation policy
\begin{equation}
\label{allocation policy: sequential}
    A_{t+1}^{FAA} = \arg \max_{i \in \mathcal{K}} \left\{ (t+1) \times \widetilde{W}^{(t)}_i(T) - N^{(t)}_i  \right\},
\end{equation}
which allocates the $(t+1)$-th sample to a design that is the most starving for it. After the simulation budget is exhausted, the design with the smallest sample mean is selected as the best.
The FAA algorithm is described in Algorithm \ref{alg:FAA}.  
\begin{breakablealgorithm}
\caption{FAA}\label{alg:FAA}
\renewcommand{\algorithmicrequire}{\textbf{Input:}}
\renewcommand{\algorithmicensure}{\textbf{Output:}}
\begin{algorithmic}[1]
\Require Set of designs $\mathcal{K}$, initial sample size $n_0$, and simulation budget $T$.
\State \textbf{Initialization}: Set $t = n_0 \times k$, and $N^{(t)}_i= n_0$, for $i \in \mathcal{K}$. Perform $n_0$ replications for each design. 
\While{$t < T$}
\State Update $\hat{\mu}^{(t)}_i$, $(\hat{\sigma}^{(t)}_i)^2$, for $i \in \mathcal{K}$, $\hat{b}^{(t)} = \arg \min_{i \in \mathcal{K}} \hat{\mu}_i^{(t)}$.
\State Calculate $\hat{w}^{*(t)}_i$, for $i \in \mathcal{K}$, according to (\ref{OCBA allocation rule}). 
\State Calculate $\hat{T}^{(t)}_0$ according to \eqref{T0}.
\State Calculate $\widetilde{W}^{(t)}_i(T)$, for $i \in \mathcal{K}$, according to (\ref{new allocation ratio for design i}) and (\ref{new allocation ratio for design b}).
\State Find $A_{t+1}^{FAA}$ according to (\ref{allocation policy: sequential}).
\State Perform additional one replication for design $A_{t}^{FAA}$.
\State Set $N^{(t+1)}_{A^{FAA}_{t+1}} = N^{(t)}_{A^{FAA}_{t+1}} + 1$, $N^{(t+1)}_i = N^{(t)}_i $, for $ i \in \mathcal{K}$ and $i \neq A^{FAA}_{t+1}$, and $t = t + 1$.
\EndWhile
\Ensure $\hat{b}^{(T)} = \arg \min_{i \in \mathcal{K}} \hat{\mu}_i^{(T)}$
\end{algorithmic}
\end{breakablealgorithm}

On the one hand, when the specified total budget $T = \mathcal{O}(k \log k)$, allocations to designs made by FAA may not exactly converge to the proposed budget-adaptive rule. Because FAA cannot iterate more than $T$ times, and there are estimation errors in sample means and variances of designs. On the other hand, when $T = \omega(k \log k)$ such that $\widetilde{W}(T)$ reduces to $w^*$, the budget allocation policy of FAA defined in \eqref{allocation policy: sequential} becomes
\begin{equation*}
    A_{t+1}^{FAA} = \arg \max_{i \in \mathcal{K}} \left\{ (t+1) \times 
    \hat{w}^{*,(t)}_i - N^{(t)}_i  \right\},
\end{equation*}
which is identical to what is defined in the OCBA algorithm (with the “most starving” technique) in \cite{chen2011stochastic}. Thus, the convergence of FAA to the asymptotic OCBA rule is guaranteed by the consistency of the OCBA algorithm proved by \cite{li2023convergence}. 

\begin{Proposition}
\label{proposition: consistency of FAA}
If the specified simulation budget $T = \omega(k \log k)$ and $k \to \infty$, for the FAA algorithm, the following statements hold (all limits hold almost surely):
\begin{itemize}
    \item $\lim_{t \rightarrow T} \hat{b}^{(t)} = b$;
    \item $\lim_{t \rightarrow T} w_i^{(t)} = w^*_i$, for $i \in \mathcal{K}.$
\end{itemize}
\end{Proposition}
Proposition \ref{proposition: consistency of FAA} shows that when  $T = \omega(k \log k)$, as iteration increases, PCS converges to 1, i.e., FAA is guaranteed to learn the best design with probability 1; and the proportional allocations to designs made by FAA converge to the asymptotic OCBA allocation ratios.

\subsubsection{Dynamic anchorage allocation}

We extend FAA to a more flexible variant, named as dynamic anchorage allocation (DAA), by allowing dynamically changing the anchored final budget instead of fixing it during the procedure. Specifically, in each iteration, the next simulation replication is anchored by DAA, and the goal becomes maximizing the PCS after the additional allocation. Again, the “most starving” technique in \cite{chen2011stochastic} can be used to obtain a budget allocation policy, which is defined as
\begin{equation}
\label{allocation policy: dynamic}
    A_{t+1}^{DAA} = \arg \max_{i \in \mathcal{K}} \left\{ (t+1) \times \widetilde{W}^{(t)}_i(t+1) - N^{(t)}_i  \right\}.
\end{equation}
The DAA algorithm can be implemented by Algorithm \ref{alg:DAA}. 
\begin{breakablealgorithm}
\caption{DAA}
\label{alg:DAA}
\renewcommand{\algorithmicrequire}{\textbf{Input:}}
\renewcommand{\algorithmicensure}{\textbf{Output:}}
\begin{algorithmic}[1]
\Require Set of designs $\mathcal{K}$, initial sample size $n_0$, and simulation budget $T$.
\State \textbf{Initialization}: Set $t = n_0 \times k$, and $N^{(t)}_i= n_0$, for $i \in \mathcal{K}$. Perform $n_0$ replications for each design.
\State Complete step 2-5 in Algorithm \ref{alg:FAA}.
\State Calculate $\widetilde{W}^{(t)}_i(t+1)$, for $i \in \mathcal{K}$, according to (\ref{new allocation ratio for design i}) and (\ref{new allocation ratio for design b}).
\State Find $A_{t+1}^{DAA}$ according to (\ref{allocation policy: dynamic}).
\State Complete step 8-10 in Algorithm \ref{alg:FAA}.
\Ensure $\hat{b}^{(T)} = \arg \min_{i \in \mathcal{K}} \hat{\mu}_i^{(T)}$.
\end{algorithmic}
\end{breakablealgorithm}

When $T = \mathcal{O}(k \log k)$, allocations to designs made by DAA may not exactly achieve the proposed budget-adaptive rule due to similar issues to FAA. Proposition \ref{proposition: consistency of DAA} demonstrates that when $T = \omega(k \log k)$, as iteration increases, PCS converges to 1, and allocations made by DAA converge to the asymptotic OCBA rule. 

\begin{Proposition}
\label{proposition: consistency of DAA}
If the specified simulation budget $T = \omega(k \log k)$ and $k \to \infty$, for the DAA algorithm, the following statements hold (all limits hold almost surely):
\begin{itemize}
    \item $\lim_{t \rightarrow T} \hat{b}^{(t)} = b$;
    \item $\lim_{t \rightarrow T} w^{(t)} = w^*_i$, for $i \in \mathcal{K}.$
\end{itemize}
\end{Proposition}

\textit{Proof:} See Online Appendix C.4.

\begin{Remark}
\normalfont{In Proposition \ref{proposition: consistency of FAA}\&\ref{proposition: consistency of DAA}, the convergence is shown by considering a large-scale problem with a specific large budget, i.e., $k \to \infty$ and $T = \omega(k \log k)$. When one considers a small-scale problem, all claims in Proposition \ref{proposition: consistency of FAA}\&\ref{proposition: consistency of DAA} are valid if one specifies a budget $T \to \infty$.}
\end{Remark}

\begin{Remark}
\normalfont{FAA anchors the final budget while DAA dynamically anchors the next replication. Clearly, DAA is more flexible than FAA as the calculation of its allocation policy in \eqref{allocation policy: dynamic} does not rely on $T$. The design motivation behind DAA is to approach optimality conditions across various budgets simultaneously. It achieves this by anchoring the next replication and maintaining an adaptive allocation policy that evolves with each iteration. In contrast, the FAA is primarily designed to reach the optimality conditions only when the total budget is fully utilized, as it ties its allocation to the final budget and does not fully “mature” until budget is exhausted. As will be shown in Section \ref{numerical experiments}, the difference in performance between FAA and DAA is subtle for both small- and large-scale problems; and the budget-adaptive property improves both algorithms' performance under small-budget conditions. In practical implementations, when the total budget is specified and the users only cares about the PCS eventually achieved, both FAA and DAA are preferred. Otherwise, DAA is recommended for implementation.}
\end{Remark}

\section{Numerical experiments}
\label{numerical experiments}
In this section, the objectives of numerical experiments are twofold: (1) to test the accuracy of APCS and the first-order Taylor series expansion; and (2) to demonstrate the desirable budget-adaptive property of both FAA and DAA improves their small-budget performances compared with the asymptotic OCBA rule. All experiments are conducted in MATLAB R2022b on a computer with Intel Core i5-10400 CPU with 2.90 GHz, 16 GB memory, a 64-bit operating system, and 6 cores with 12 logical processors. Source code for all experiments is available at \url{https://github.com/Haowei-Wang/Budget-adaptive-rule-for-OCBA}.

\subsection{Accuracy test for APCS and Taylor series expansion}

\begin{figure}[t]
     \centering
     \begin{subfigure}[b]{0.495\textwidth}
         \centering
         \includegraphics[width=\textwidth]{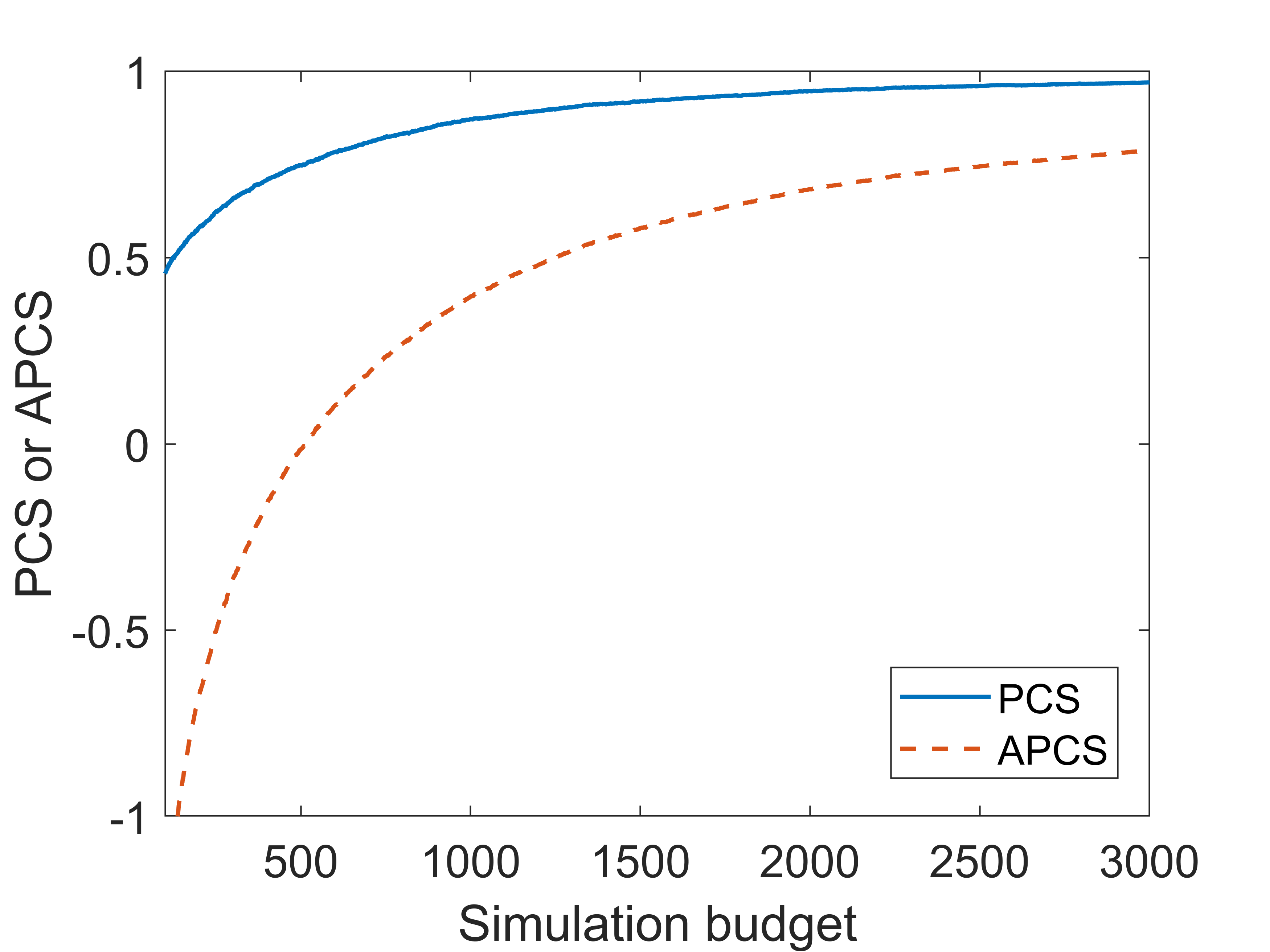}
         \caption{$k=10$}
         \label{sfig:Comparison of PCS and APCS with k = 10}
     \end{subfigure}
     \hfill
     \begin{subfigure}[b]{0.495\textwidth}
         \centering
         \includegraphics[width=\textwidth]{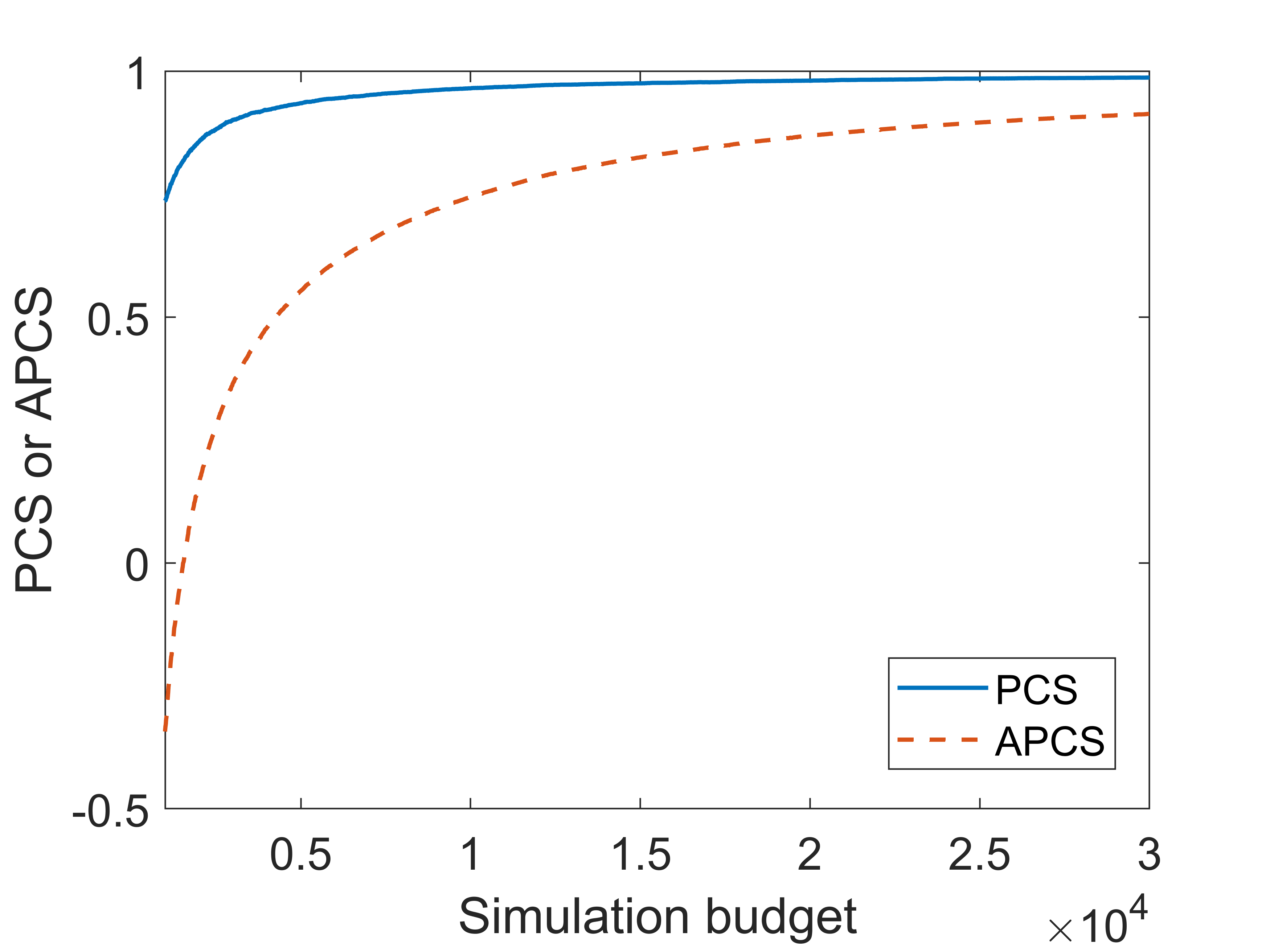}
         \caption{$k = 10^2$}
         \label{sfig:Comparison of PCS and APCS with k = 10^2}
     \end{subfigure}
     \hfill
     \begin{subfigure}[b]{0.495\textwidth}
         \centering
         \includegraphics[width=\textwidth]{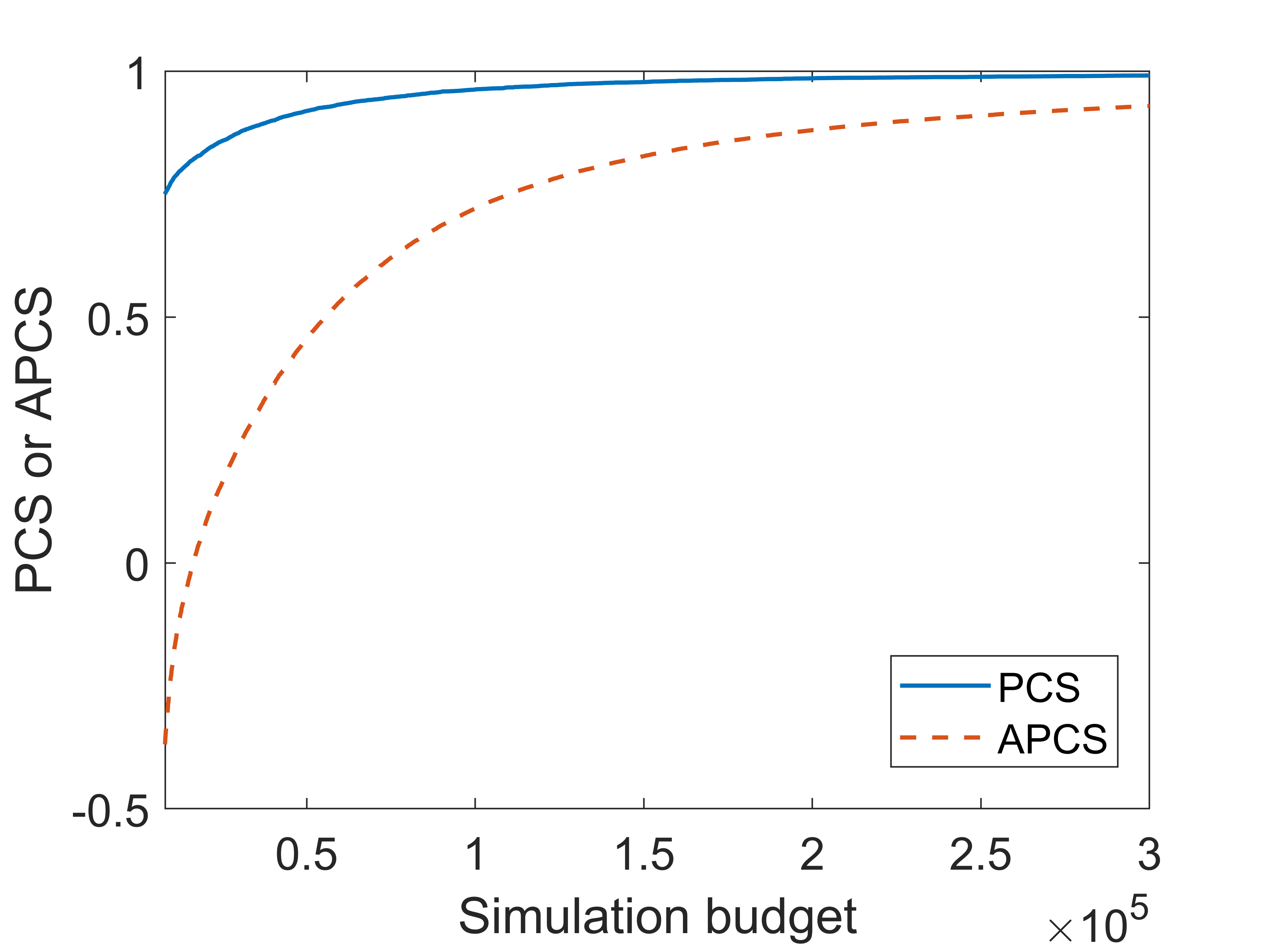}
         \caption{$k = 10^3$}
         \label{sfig:Comparison of PCS and APCS with k = 10^3}
     \end{subfigure}
     \hfill
     \begin{subfigure}[b]{0.495\textwidth}
         \centering
         \includegraphics[width=\textwidth]{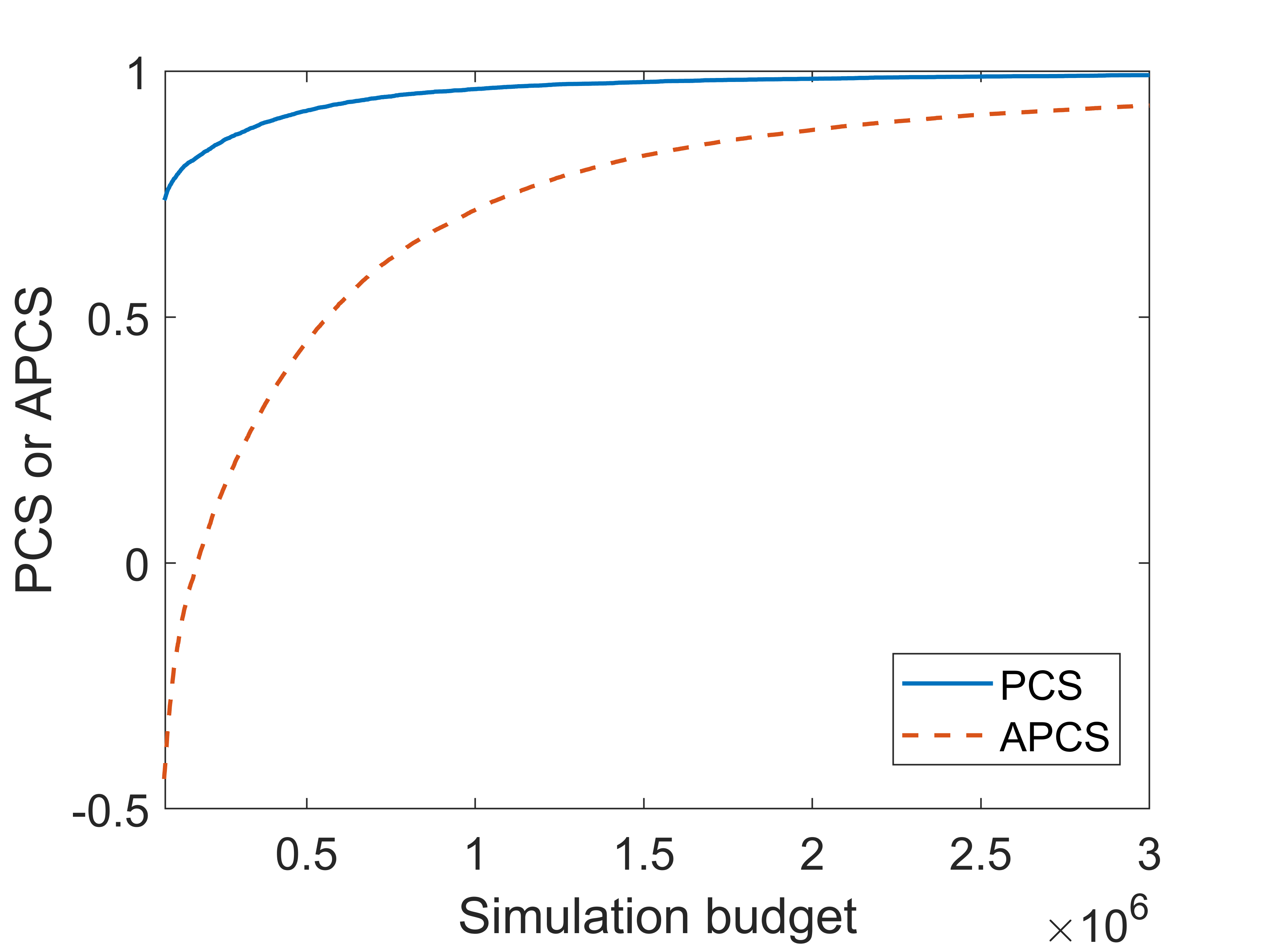}
         \caption{$k = 10^4$}
         \label{sfig:Comparison of PCS and APCS with k = 10^4}
     \end{subfigure}
     \caption{(Color online) Comparison of PCS and APCS}
    \label{fig:Comparison of PCS and APCS}
\end{figure}

To demonstrate the accuracy of APCS and the first-order Taylor series expansion, we  consider a synthetic problem setting where there are $k$ alternative designs with sampling distributions $N(i,10^2)$, for $i = 1,2,\dots,k$. The identity of the best design is fixed, i.e., $b=1$.

The APCS value is calculated based on known distribution parameters, and it is used to evaluate the gap between PCS and APCS. While calculating PCS and APCS, the simulation budget is allocated by DAA without any prior knowledge of underlying distributions. The initial number of replications is set as $3$, i.e., $n_0 = 3$. Figure \ref{fig:Comparison of PCS and APCS} illustrates the value of PCS and APCS based on the average of $10^5$ independent replications with $k$ ranging from $10$ to $10^4$. Clearly, the gap between PCS and APCS decreases with the growth of simulation budget. APCS reflects the general trend of PCS even for small budgets, which is much more important than the value of APCS itself for deriving an efficient allocation rule because the derivative information plays key roles in  the optimality conditions \citep{gao2017new}. Therefore, APCS is a good approximation of PCS for deriving a budget-adaptive rule.

\begin{figure}[t]
     \centering
     \begin{subfigure}[b]{0.495\textwidth}
         \centering
         \includegraphics[width=\textwidth]{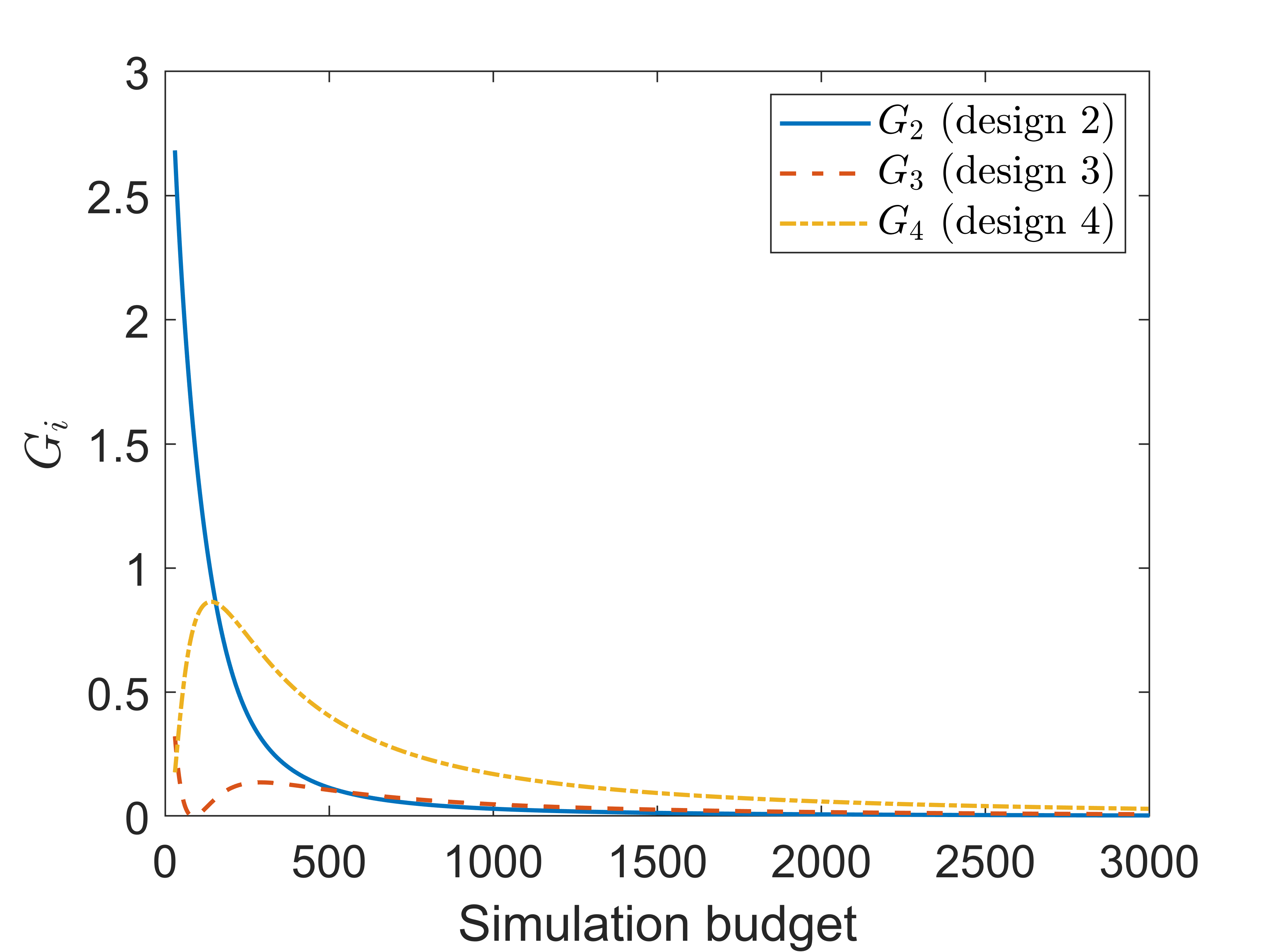}
         \caption{$k=10$}
         \label{sfig:Illustration of G_i with k = 10}
     \end{subfigure}
     \hfill
     \begin{subfigure}[b]{0.495\textwidth}
         \centering
         \includegraphics[width=\textwidth]{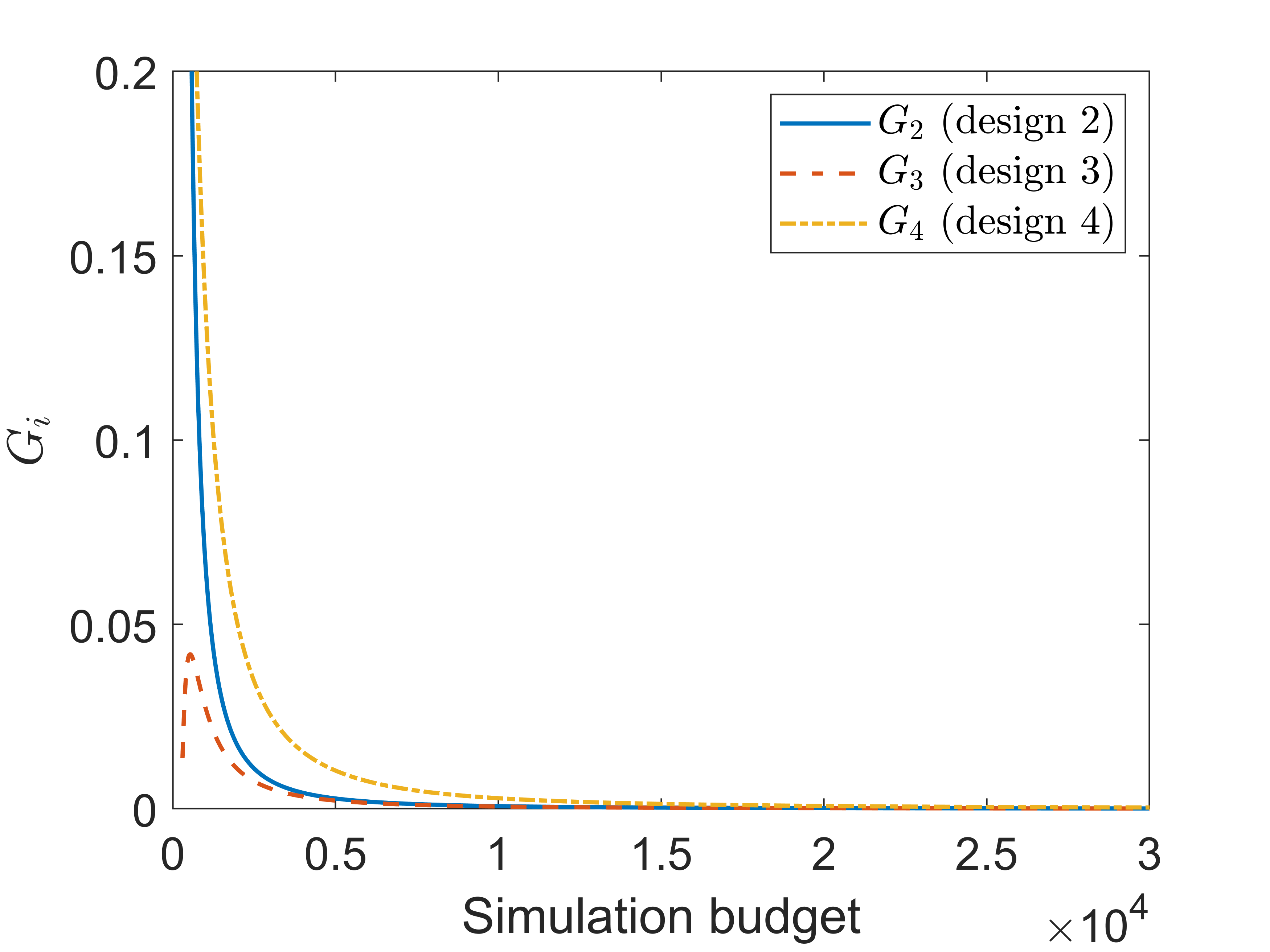}
         \caption{$k = 10^2$}
         \label{sfig:Illustration of G_i with k = 10^2}
     \end{subfigure}
     \hfill
     \begin{subfigure}[b]{0.495\textwidth}
         \centering
         \includegraphics[width=\textwidth]{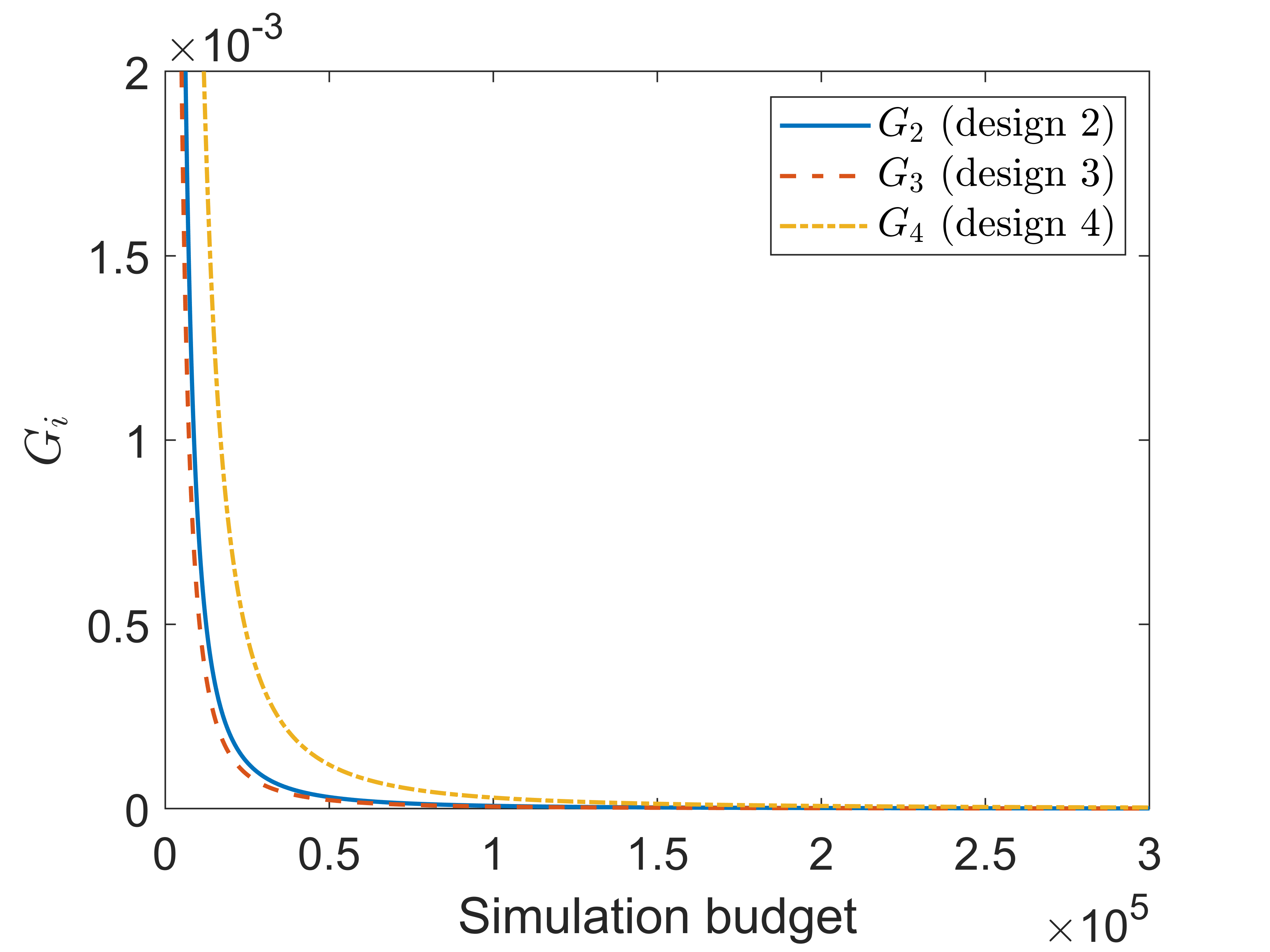}
         \caption{$k = 10^3$}
         \label{sfig:Illustration of G_i with k = 10^3}
     \end{subfigure}
     \hfill
     \begin{subfigure}[b]{0.495\textwidth}
         \centering
         \includegraphics[width=\textwidth]{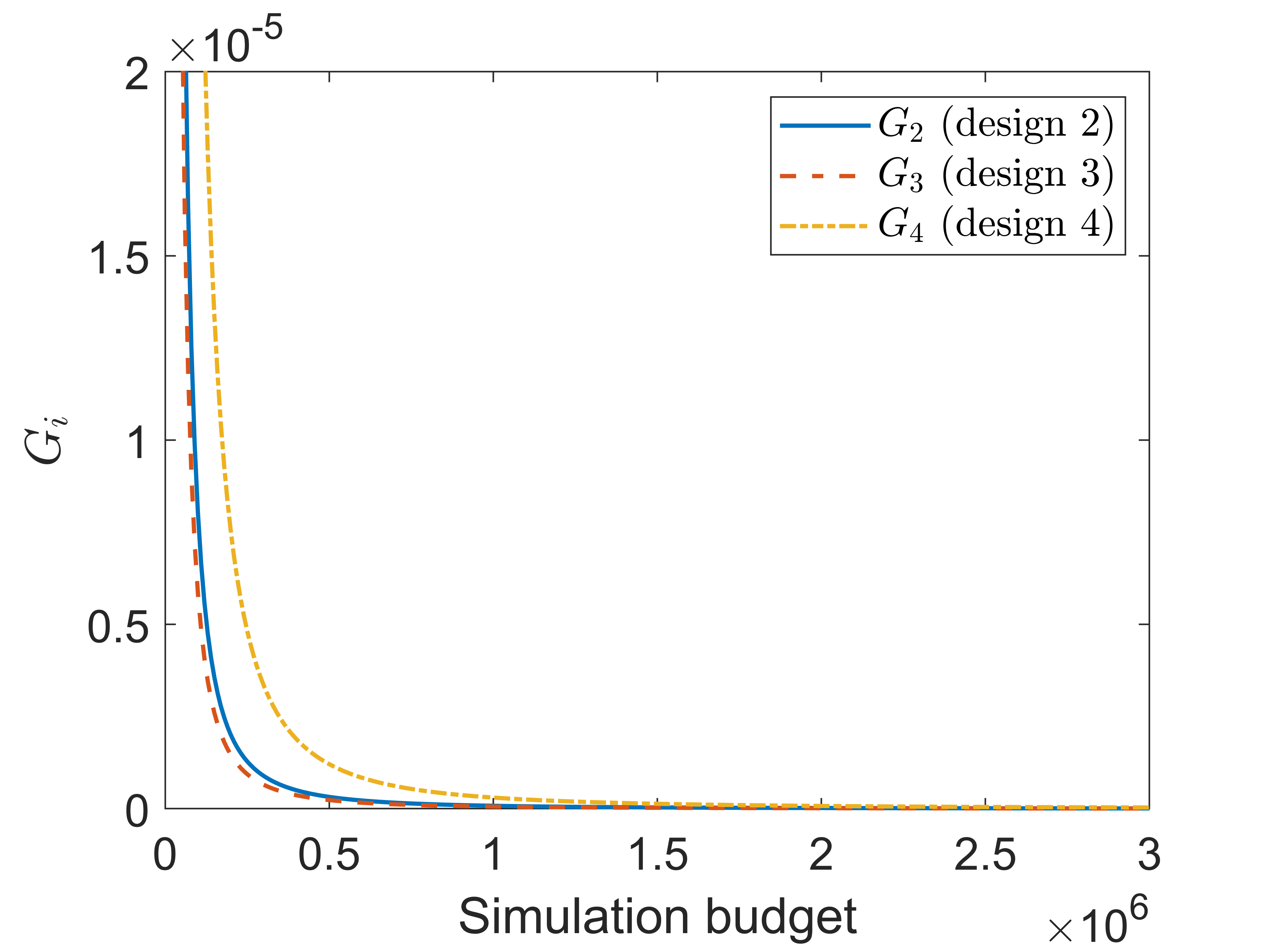}
         \caption{$k = 10^4$}
         \label{sfig:Illustration of G_i with k = 10^4}
     \end{subfigure}
     \caption{(Color online) Illustration of $G_i$}
    \label{fig:Illustration of G_i}
\end{figure}

To investigate the accuracy of the first-order Taylor series expansion, we define
\begin{equation*}
    G_i = \left|\log w_i  - \log w_i^* - \left(w_i - w_i^*\right)/w_i^* \right|, \ \text{for} \ i \in \mathcal{K}^\prime,
\end{equation*}
where $w = (w_1, w_2, \dots, w_k)$ is the exact solution to $\mathcal{C}_1$, $\mathcal{C}_3$, $\mathcal{C}_4$, and \eqref{simplified condition C2-1 (section 4.2)}, and $w^* = (w_1^*, w_2^*, \dots, w_k^*)$ is the asymptotic OCBA rule. For $i \in \mathcal{K}^\prime$, $G_i$ measures the gap between $\log w_i$ and its first-order Taylor series expansion at $w_i^*$. The $G_i$ value is calculated based on known distribution parameters. In calculating $G_i$, $w$ is determined by a non-linear optimization solver in MATLAB. Figure \ref{fig:Illustration of G_i} shows the value of $G_i$ under different budgets with $k$ ranging from $10$ to $10^4$; for better visualization, only the value of $G_2$, $G_3$, and $G_4$ are shown. Clearly, $G_i$, for $i = 2,3,4$, decrease rapidly to 0 as simulation budget increases. Additionally, these gaps are extremely small even for small budgets when $k$ is large, e.g., when $k = 10^4$. Therefore, we conclude that the first-order Taylor series expansion of $\log w_i$ at $w_i^*$ is a good approximation of $\log w_i$, especially for large-scale problems.

\subsection{Efficiency test for FAA and DAA}

\subsubsection{Benchmark algorithms}

We use four algorithms as benchmarks for comparison.
\begin{itemize}
    \item Equal allocation (EA). The simulation budget is equally allocated to all designs, i.e., $N_i = T/k$ and $w_i = 1/k$, for $i \in \mathcal{K}$. The equal allocation is a good benchmark for performance comparison.
    \item OCBA allocation \citep{chen2000simulation}. OCBA is guided by the asymptotically optimal allocation rule defined in \eqref{OCBA allocation rule}. We implement a fully sequential OCBA procedure, which allocates a single replication in each iteration according to the “most starving” technique in \cite{chen2011stochastic}. Similarly, in each iteration, sample means and variances are used as plug-in estimates for the true means and variances to calculate the OCBA allocation ratios.
    \item AOAP allocation \citep{peng2018ranking}. AOAP is an efficient budget allocation procedure that achieves both one-step-head optimality and asymptotic optimality. It requires the variances of designs to be known, and again, we use sample variances as plug-in estimates for the true variances. As iteration increases, AOAP achieves the asymptotically optimal budget allocation ratios defined in \eqref{large deviation ratio}.
    \item FBKT allocation \citep{hong2022solving}. FBKT is an efficient fixed-budget allocation algorithm specifically designed for solving large-scale problems. Different from the other benchmarks and our proposed algorithms, FBKT does not aim to achieve a certain static allocation rule. We implement a sequential FBKT procedure (with $\phi = 3$) instead of its variant for parallel computing environment.
\end{itemize}

\subsubsection{Test problems}

To demonstrate the efficiency of the proposed FAA and DAA, we consider six problem settings, in which both small- and large-scale problems are included.

\textit{Example 1 (small-scale problem)}: There are 10 alternative designs with sampling distributions $N(i,6^2)$, for $i = 1,2,\dots,10$. The goal is to identify the best design via simulation samples, where the best is $b=1$ in this example.

\textit{Example 2 (small-scale problem)}: This is a variant of Example 1. All settings are the same except that the variance is decreasing with respect to the indices. In this example, better designs are with larger variances. The designs' sampling distributions are $N(i,(11-i)^2)$, for $i = 1,2,\dots,10$. Again, the best design is $b = 1$.

\textit{Example 3 (small-scale problem)}: This is another variant of Example 1 with larger number of designs and variances. The designs' sampling distributions are $N(i,10^2)$, for $i = 1,2,\dots,50$. Again, the best design is $b = 1$.

\textit{Example 4 (medium-scale problem)}: There are 500 normal alternative designs. The sampling distribution of the best design, design $1$, is $N(0,6^2)$. As for non-best designs $i$, for $i = 2,3,\dots,500$, their sampling distributions are $N(\mu_i, \sigma_i^2)$, where $\mu_i$ and $\sigma_i$ are generated from two uniform distributions $U(1,16)$ and $U(3,9)$, respectively.

\textit{Example 5 (large-scale problem)}: This is a variant of Example 4. The identity of the best designs is fixed, but the number of designs is increased to $10^4$.

\textit{Example 6 (facility location problem)}:
The facility location problem is a practical test problem provided by the Simulation Optimization Library (\url{https://github.com/simopt-admin/simopt}) and has also been studied in \cite{gao2016new}. There is a company selling one product that will never be out of stock in a city. Without loss of generality, the city is assumed to be a unit square, i.e., $[0,1]^2$, and distances are measured in units of 30 km. Two warehouses are located in the city and each of them has 10 trucks delivering orders individually. Orders are generated from 8 AM to 5 PM by a stationary Poisson process with a rate parameter 1/3 per minutes and are located in the city according to a density function
\begin{equation*}
    f(x,y) = 1.6 - ( |x-0.8| + |y - 0.8| ), \quad x, y \in [0,1].
\end{equation*}
When order arrives, it is dispatched to the nearest warehouse with available trucks. Otherwise, it is placed into a queue and satisfied by following the first-in-first-out pattern when trucks become idle. Then, the trucks pick the order up, travel to the delivery point, deliver the products and return to their assigned warehouses waiting for the next order, where the pick-up and deliver time are exponentially distributed with mean 5 and 10, respectively. All trucks travel in Manhattan fashion at a constant speed 30 km/hour, and orders must be delivered on the day when it is received. The objective is to find the locations of the two warehouses that can maximize the proportion of orders which are delivered within 60 minutes. Let $(z_{i,1},z_{i,2})$ and $(z_{i,3},z_{i,4})$ be the two locations, respectively. We consider 10 alternatives $(z_{i,1},z_{i,2},z_{i,3},z_{i,4}) = (0.49+0.01i,0.59+0.01i,0.59+0.01i,0.79+0.01i)$, for $i = 1,2,\dots,10$. In this experiment, we run 30 days of simulation in each replication, and the proportion of orders satisfied within 60 minutes is the average proportion of satisfied orders during the 30 days. Thus, the proportion of orders satisfied within 60 minutes is approximately normally distributed. By comparing $100,000$ simulation samples of each design, the best design is determined, i.e., $(z_{1,1},z_{1,2},z_{1,3},z_{1,4}) = (0.5,0.6,0.6,0.8)$. 

\begin{figure}[t]
     \centering
     \begin{subfigure}[b]{0.267\textwidth}
         \centering
         \includegraphics[width=\textwidth]{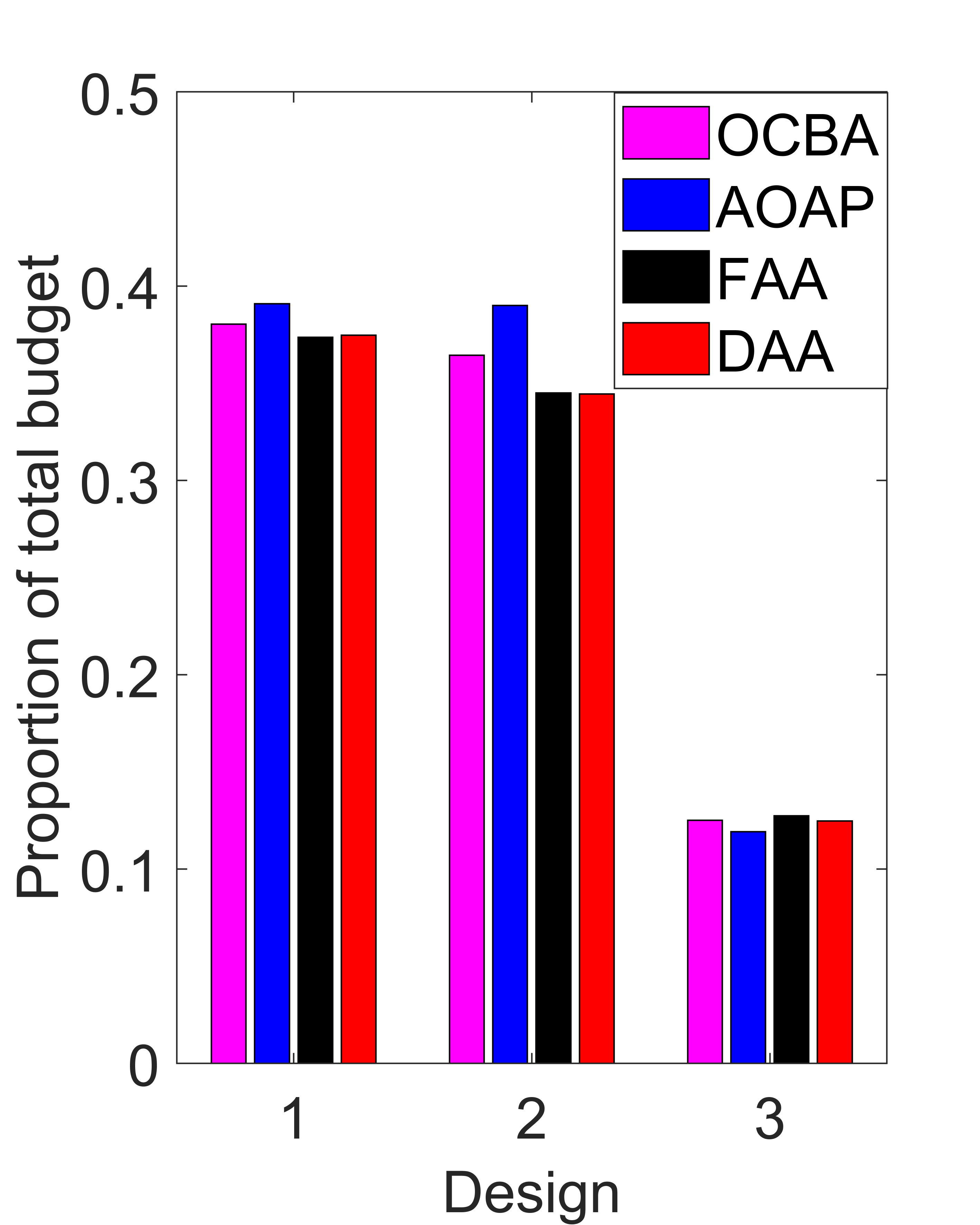}
         \caption{Budget allocation ratios of top 3 designs}
         \label{sfig:Example 1 ratios of top 3 designs}
     \end{subfigure}
     \hfill
     \begin{subfigure}[b]{0.267\textwidth}
         \centering
         \includegraphics[width=\textwidth]{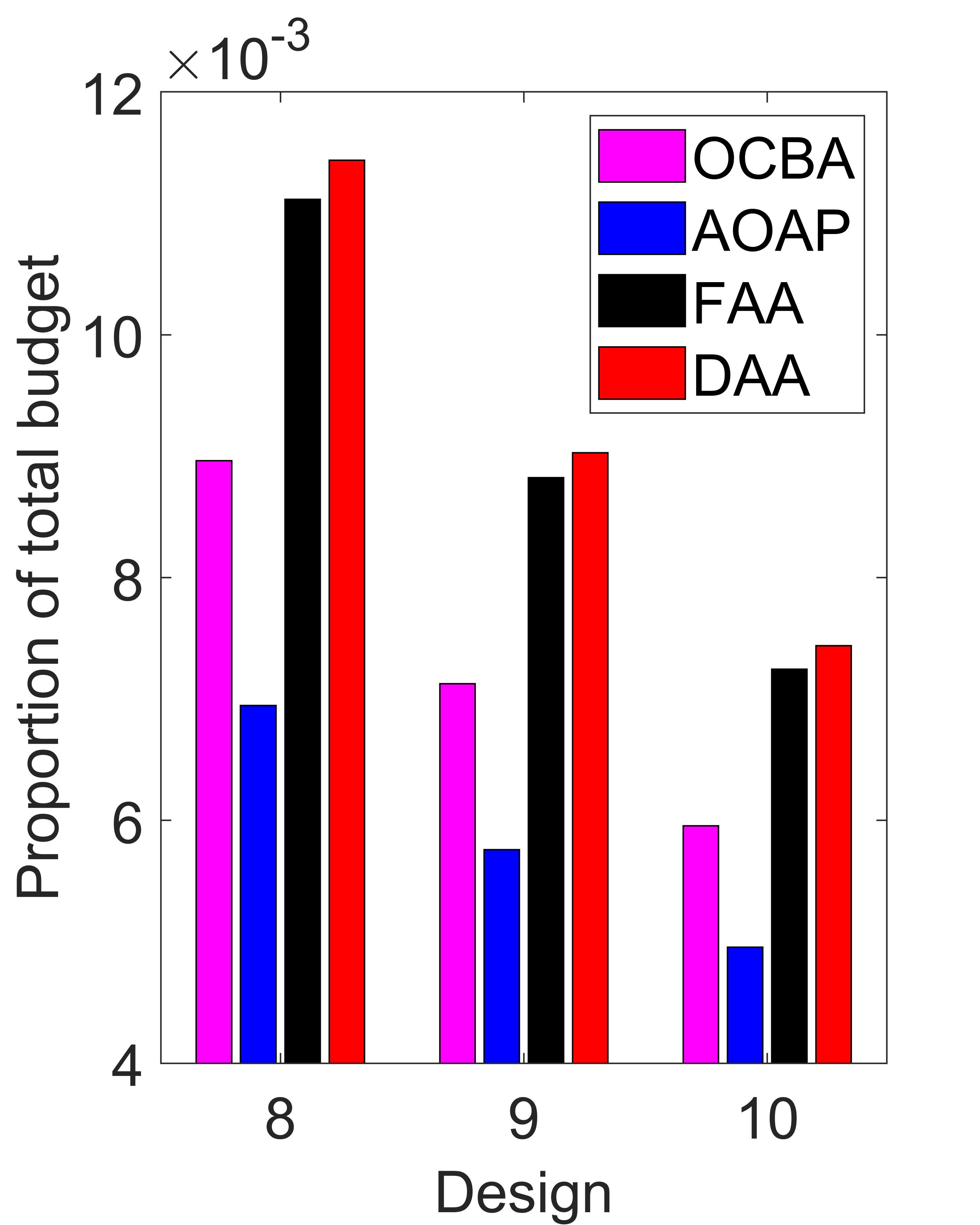}
         \caption{Budget allocation ratios of last 3 designs}
         \label{sfig:Example 1 ratios of last 3 designs}
     \end{subfigure}
     \hfill
     \begin{subfigure}[b]{0.45\textwidth}
         \centering
         \includegraphics[width=\textwidth]{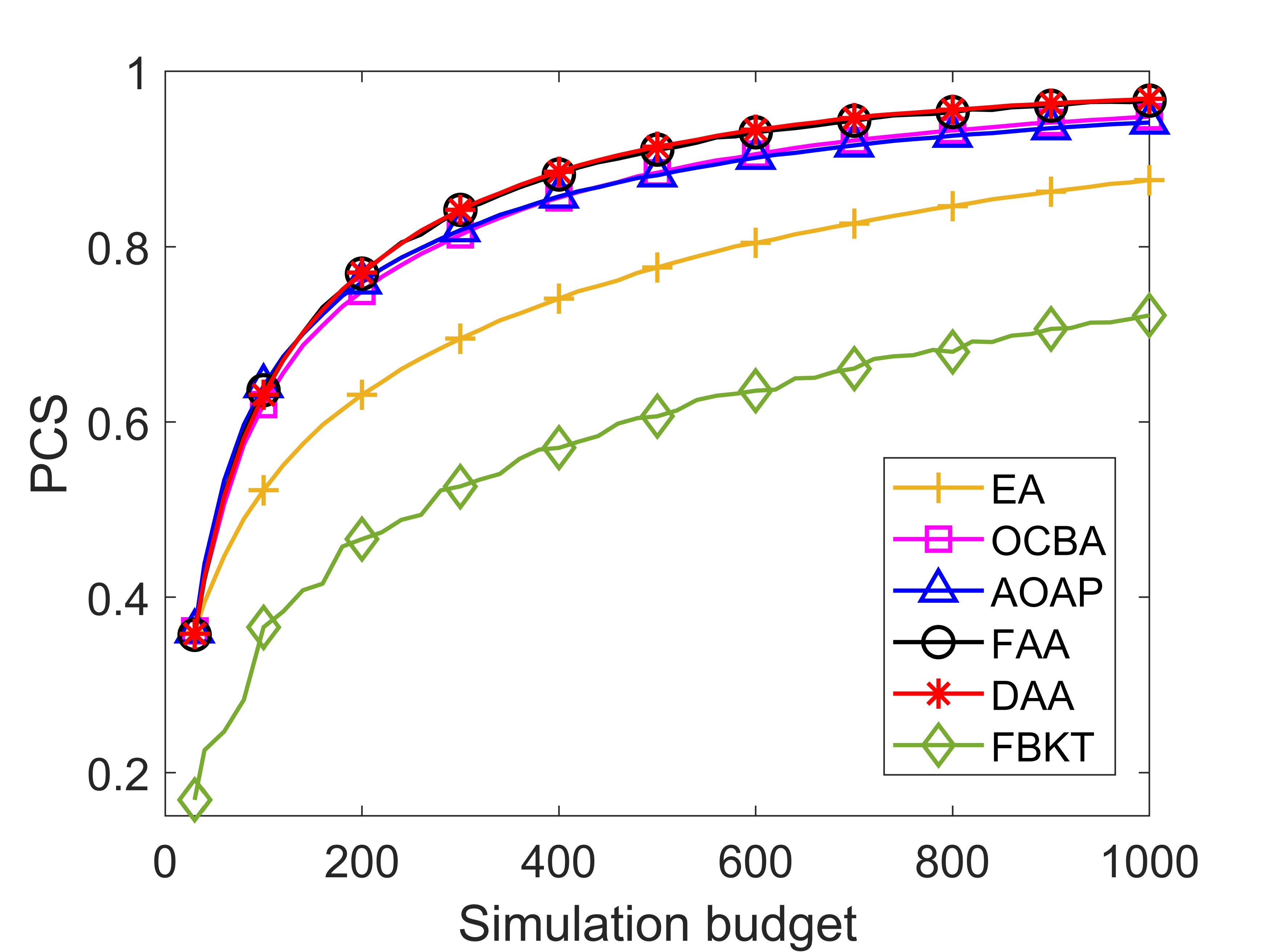}
         \caption{Comparison of PCS for the five competing procedures on Example 1}
         \label{sfig:Example 1 PCS of the five procedures}
     \end{subfigure}
     \caption{(Color online) Illustration of Example 1}
    \label{fig:Illustration of Example 1}
\end{figure}

\begin{figure}[t]
     \centering
     \begin{subfigure}[b]{0.267\textwidth}
         \centering
         \includegraphics[width=\textwidth]{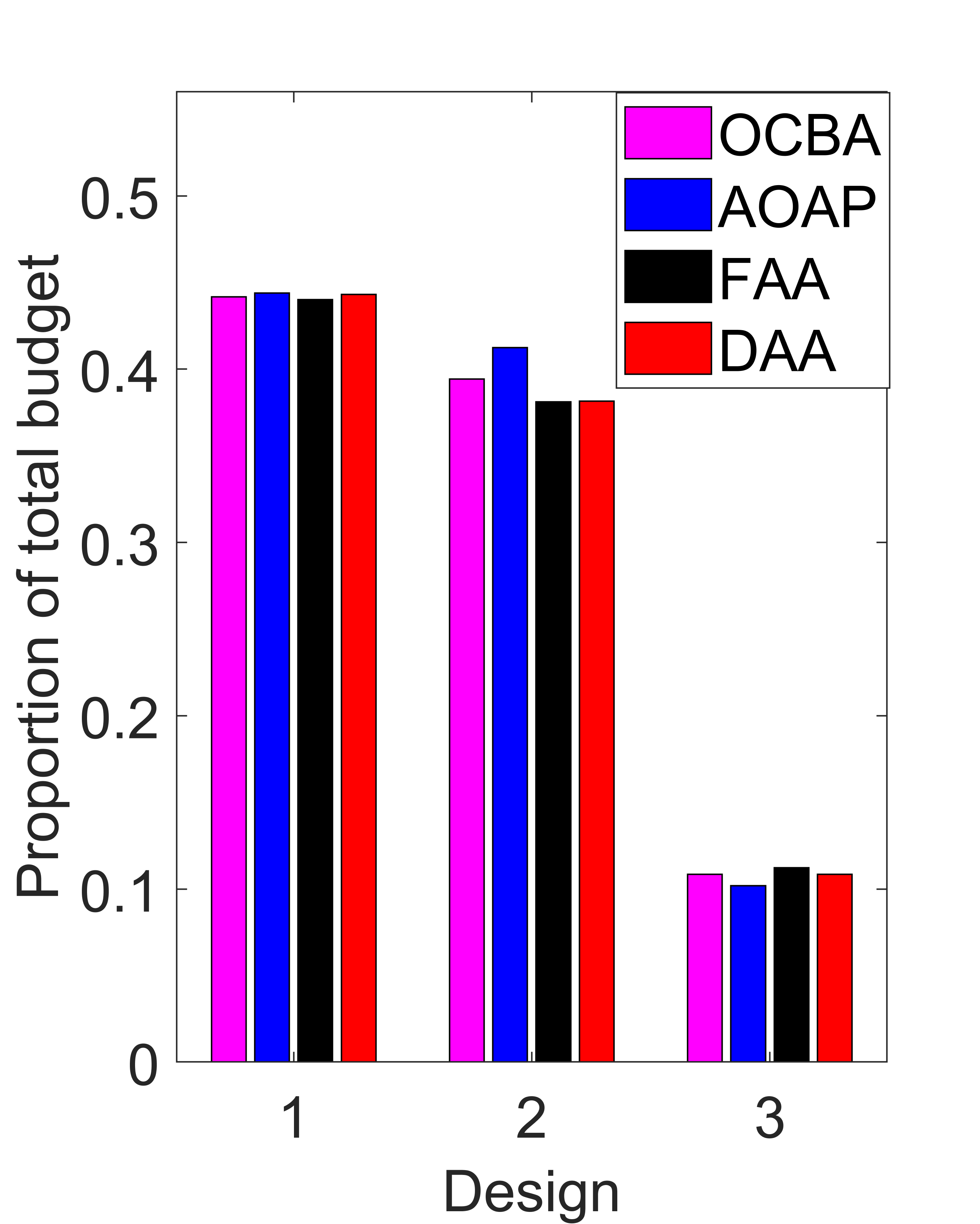}
         \caption{Budget allocation ratios of top 3 designs}
         \label{sfig:Example 2 ratios of top 3 designs}
     \end{subfigure}
     \hfill
     \begin{subfigure}[b]{0.267\textwidth}
         \centering
         \includegraphics[width=\textwidth]{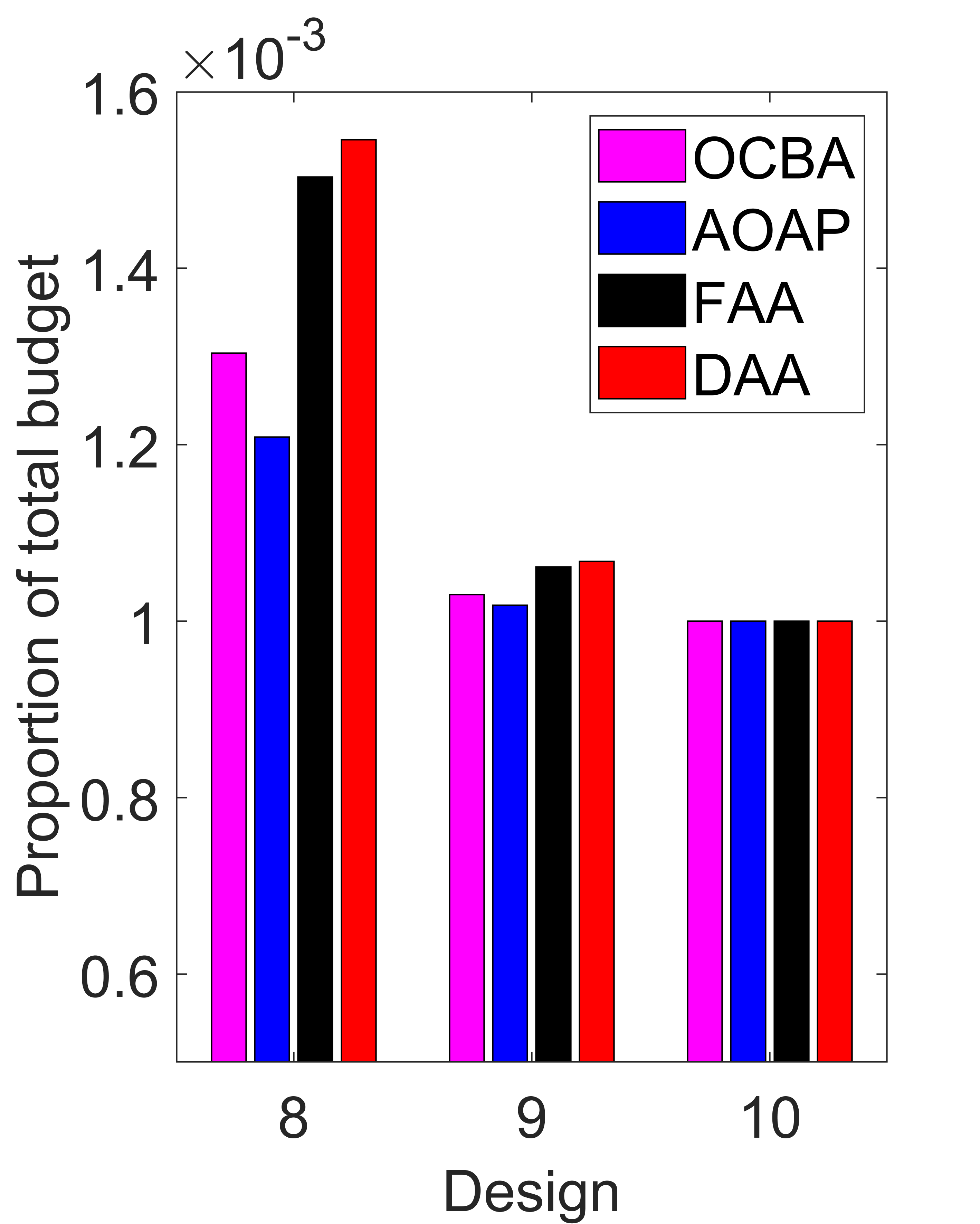}
         \caption{Budget allocation ratios of last 3 designs}
         \label{sfig:Example 2 ratios of last 3 designs}
     \end{subfigure}
     \hfill
     \begin{subfigure}[b]{0.45\textwidth}
         \centering
         \includegraphics[width=\textwidth]{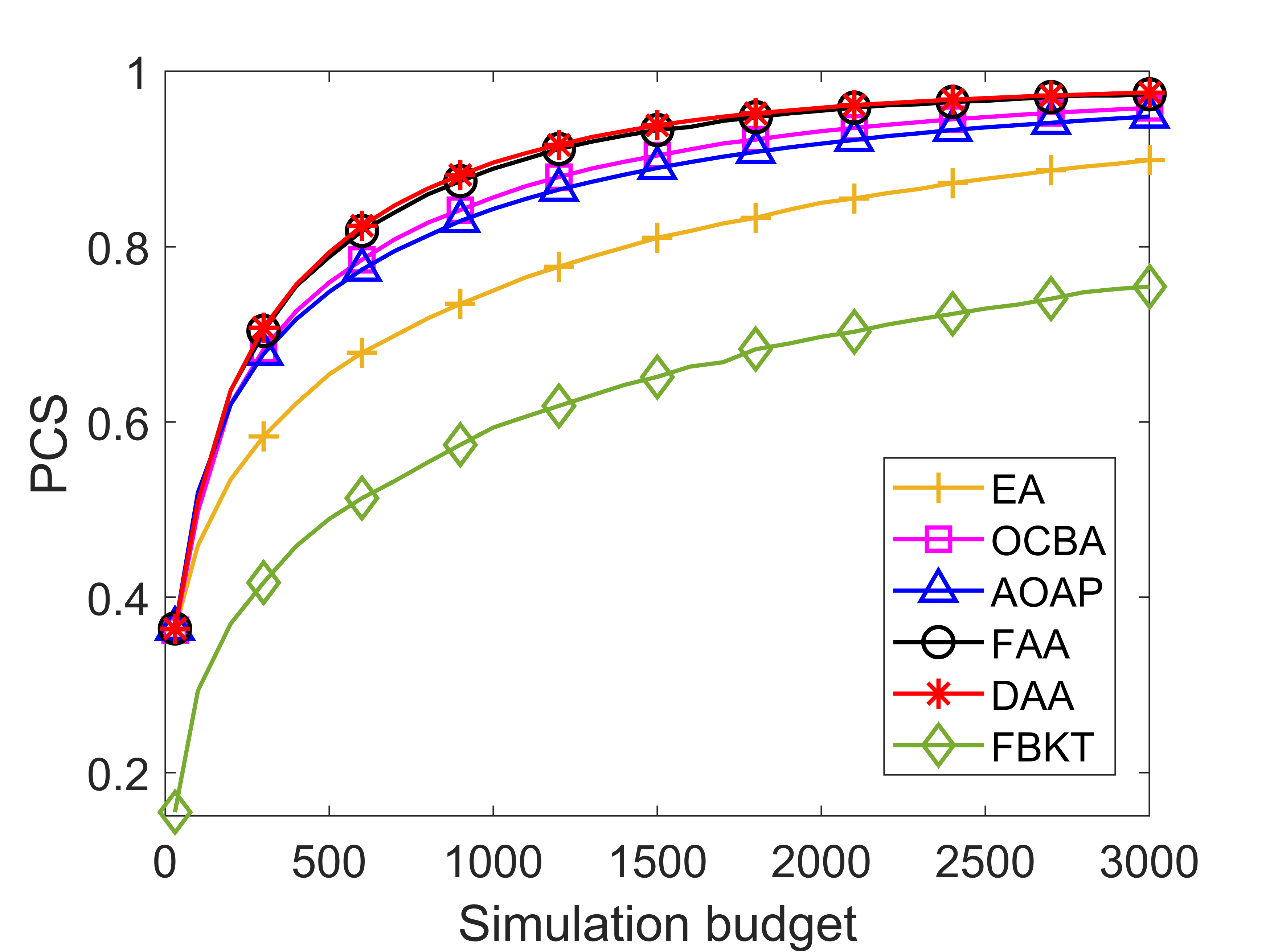}
         \caption{Comparison of PCS for the five competing procedures on Example 2}
         \label{sfig:Example 2 PCS of the five procedures}
     \end{subfigure}
     \caption{(Color online) Illustration of Example 2}
    \label{fig:Illustration of Example 2}
\end{figure}

\begin{figure}[t]
     \centering
     \begin{subfigure}[b]{0.267\textwidth}
         \centering
         \includegraphics[width=\textwidth]{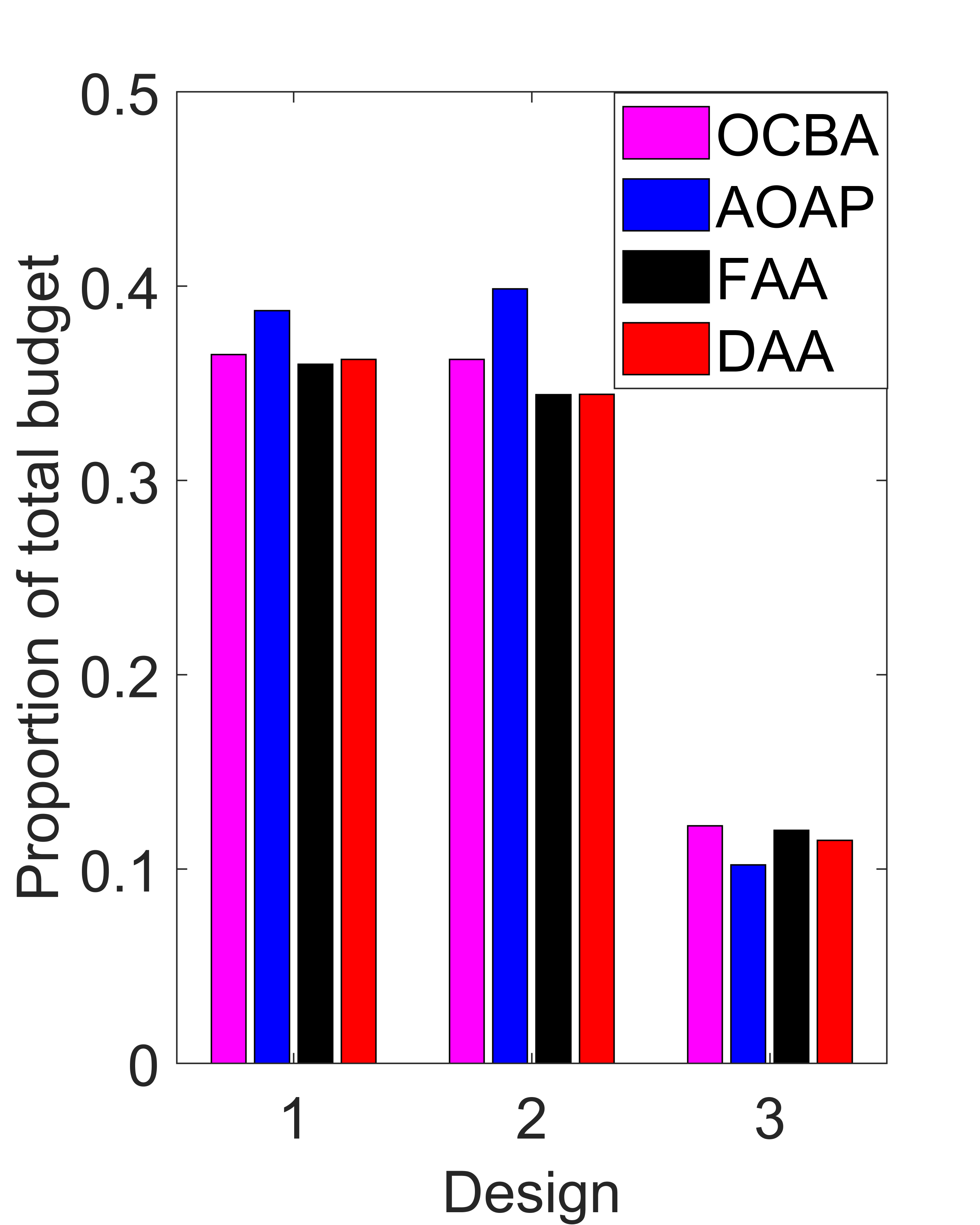}
         \caption{Budget allocation ratios of top 3 designs}
         \label{sfig:Example 3 ratios of top 3 designs}
     \end{subfigure}
     \hfill
     \begin{subfigure}[b]{0.267\textwidth}
         \centering
         \includegraphics[width=\textwidth]{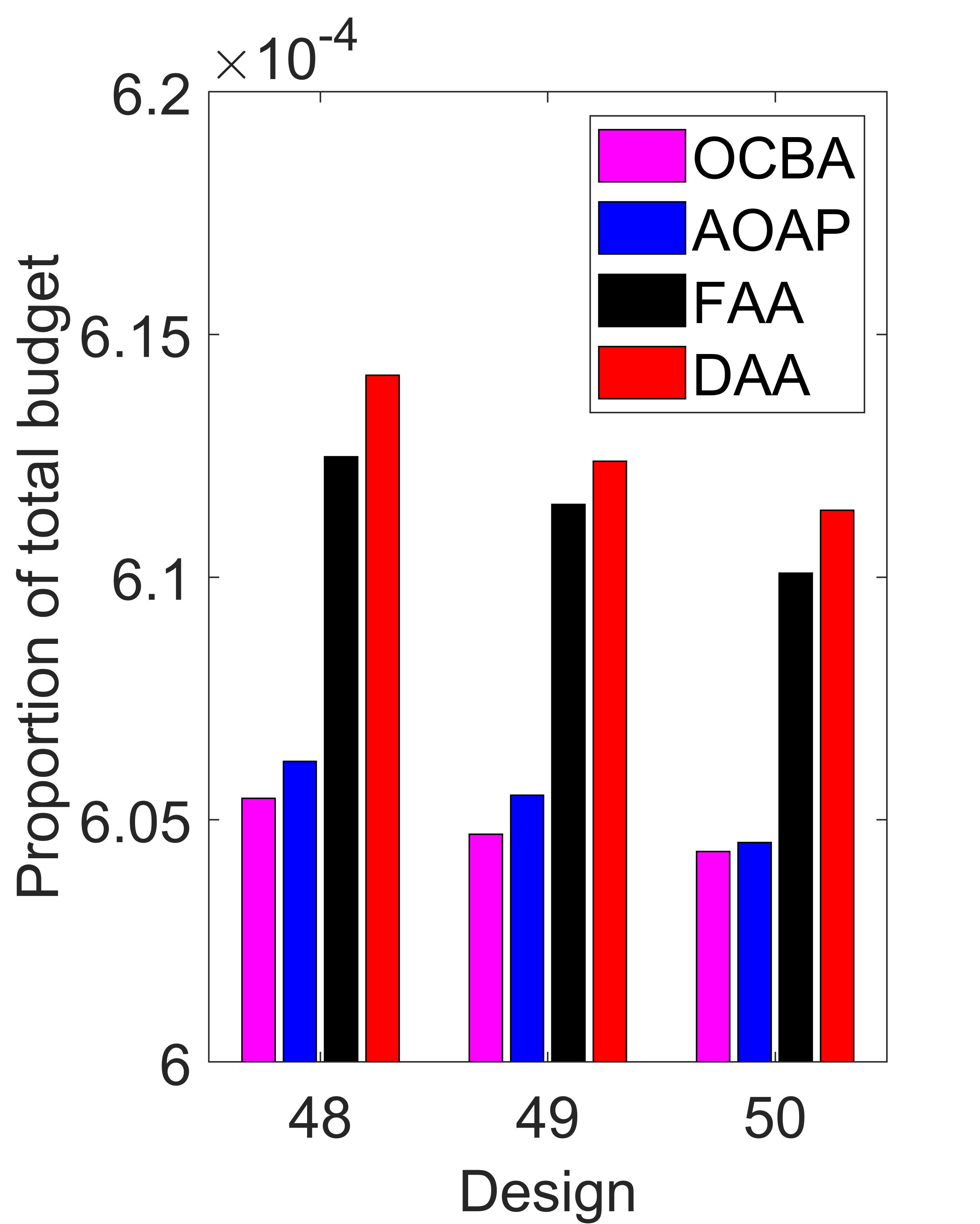}
         \caption{Budget allocation ratios of last 3 designs}
         \label{sfig:Example 3 ratios of last 3 designs}
     \end{subfigure}
     \hfill
     \begin{subfigure}[b]{0.45\textwidth}
         \centering
         \includegraphics[width=\textwidth]{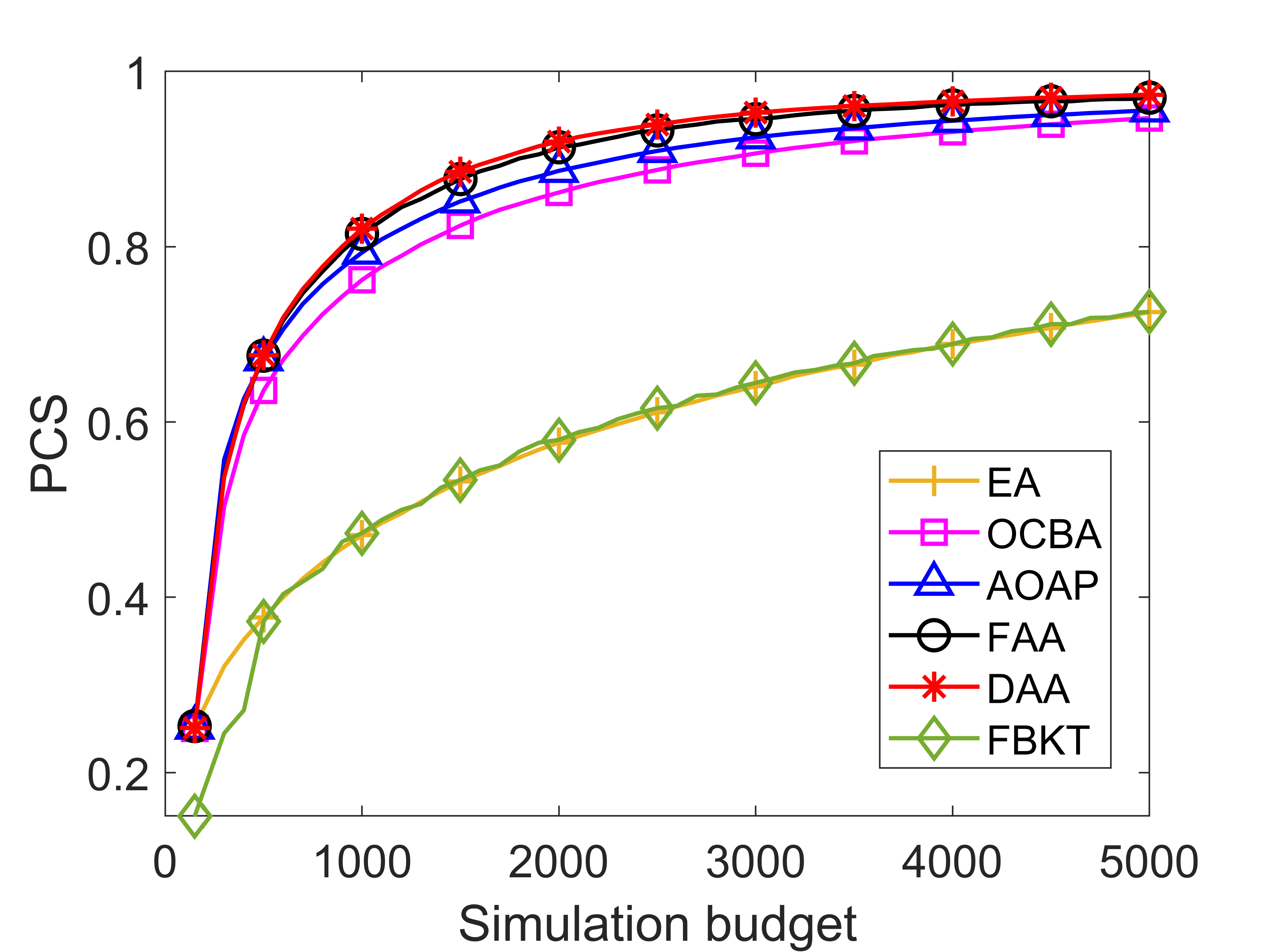}
         \caption{Comparison of PCS for the five competitive procedures on Example 3}
         \label{sfig:Example 3 PCS of the five procedures}
     \end{subfigure}
     \caption{(Color online) Illustration of Example 3}
    \label{fig:Illustration of Example 3}
\end{figure}

\begin{figure}[t]
     \centering
     \begin{subfigure}[b]{0.267\textwidth}
         \centering
         \includegraphics[width=\textwidth]{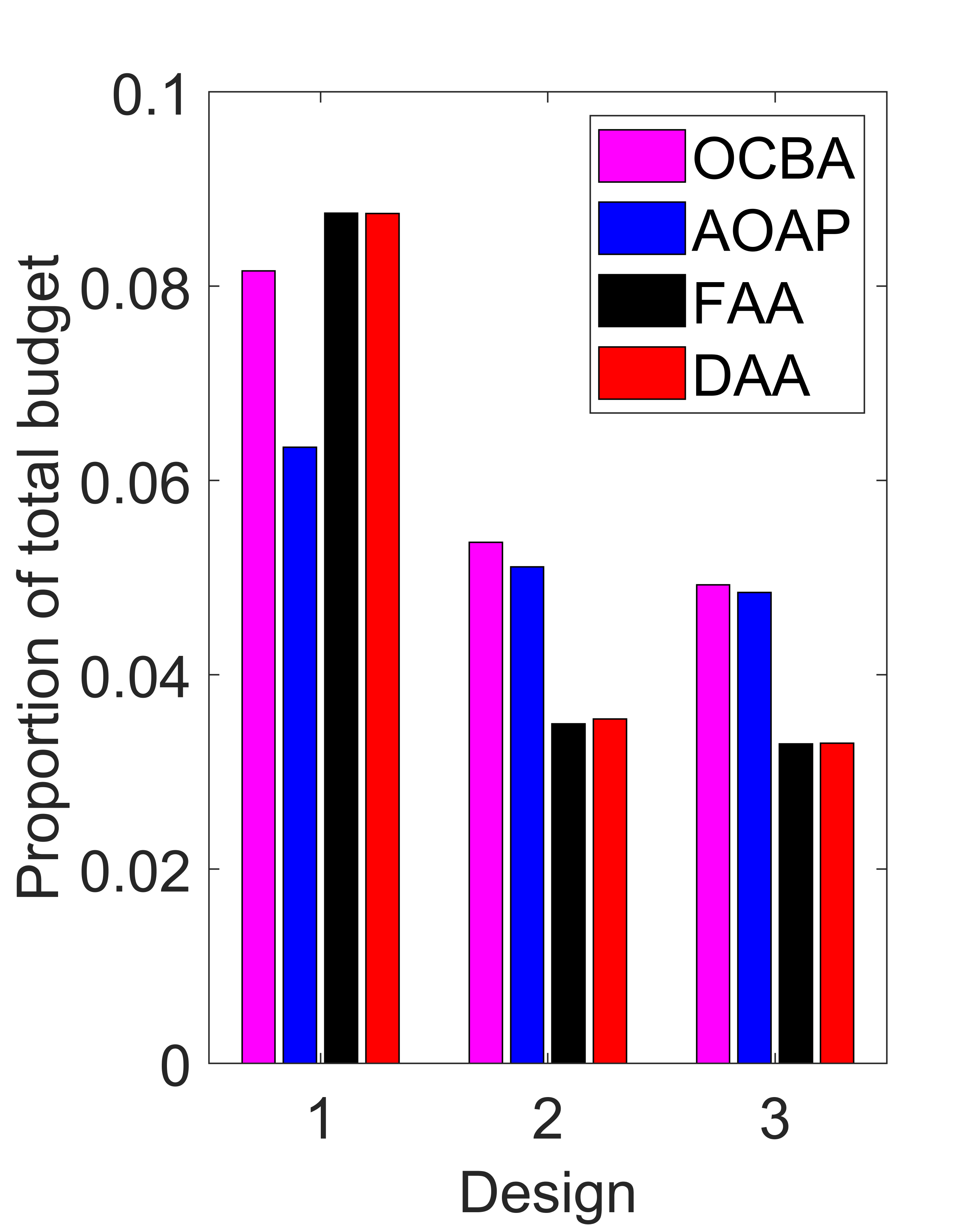}
         \caption{Budget allocation ratios of top 3 designs}
         \label{sfig:Example 4 ratios of top 3 designs}
     \end{subfigure}
     \hfill
     \begin{subfigure}[b]{0.267\textwidth}
         \centering
         \includegraphics[width=\textwidth]{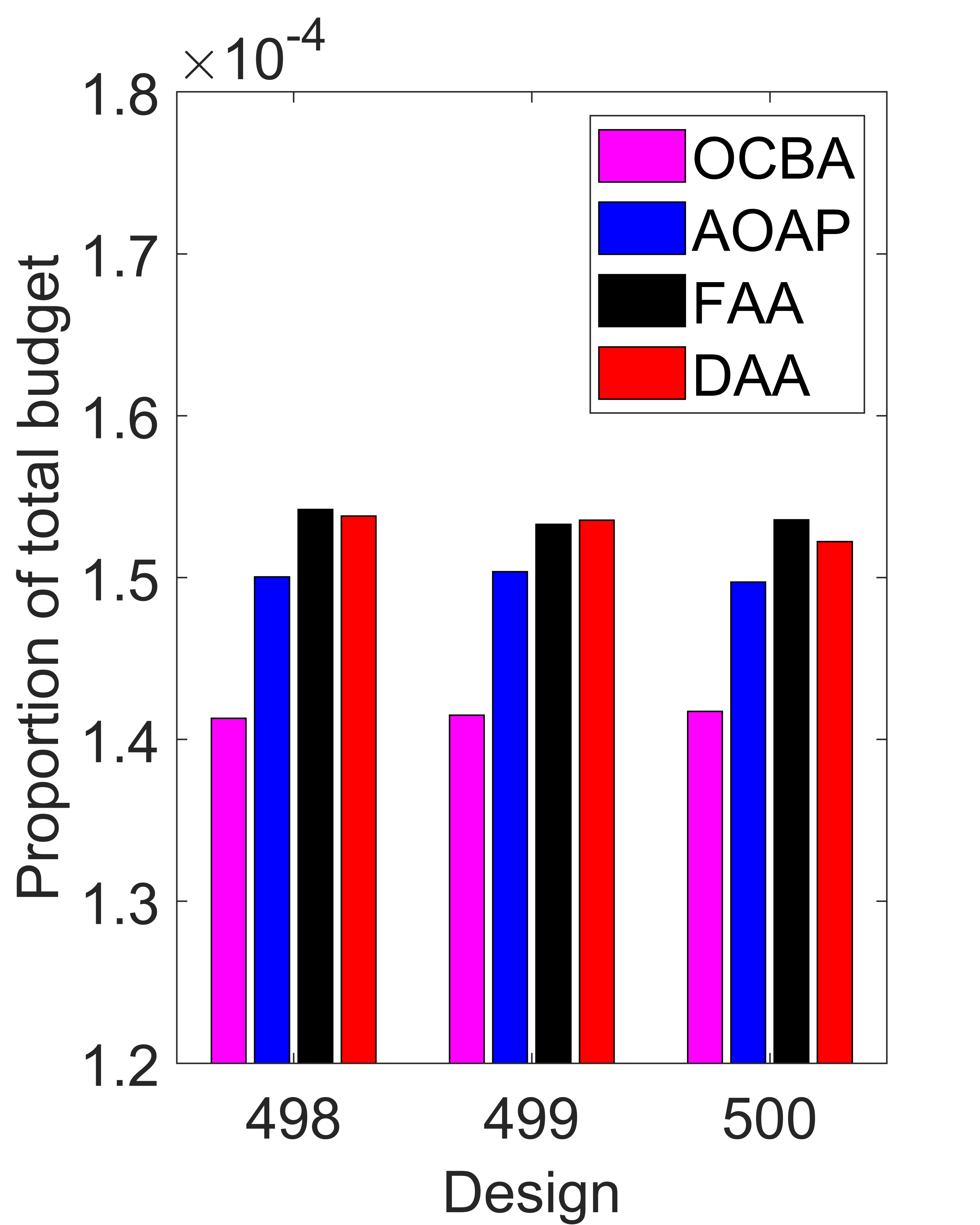}
         \caption{Budget allocation ratios of last 3 designs}
         \label{sfig:Example 4 ratios of last 3 designs}
     \end{subfigure}
     \hfill
     \begin{subfigure}[b]{0.45\textwidth}
         \centering
         \includegraphics[width=\textwidth]{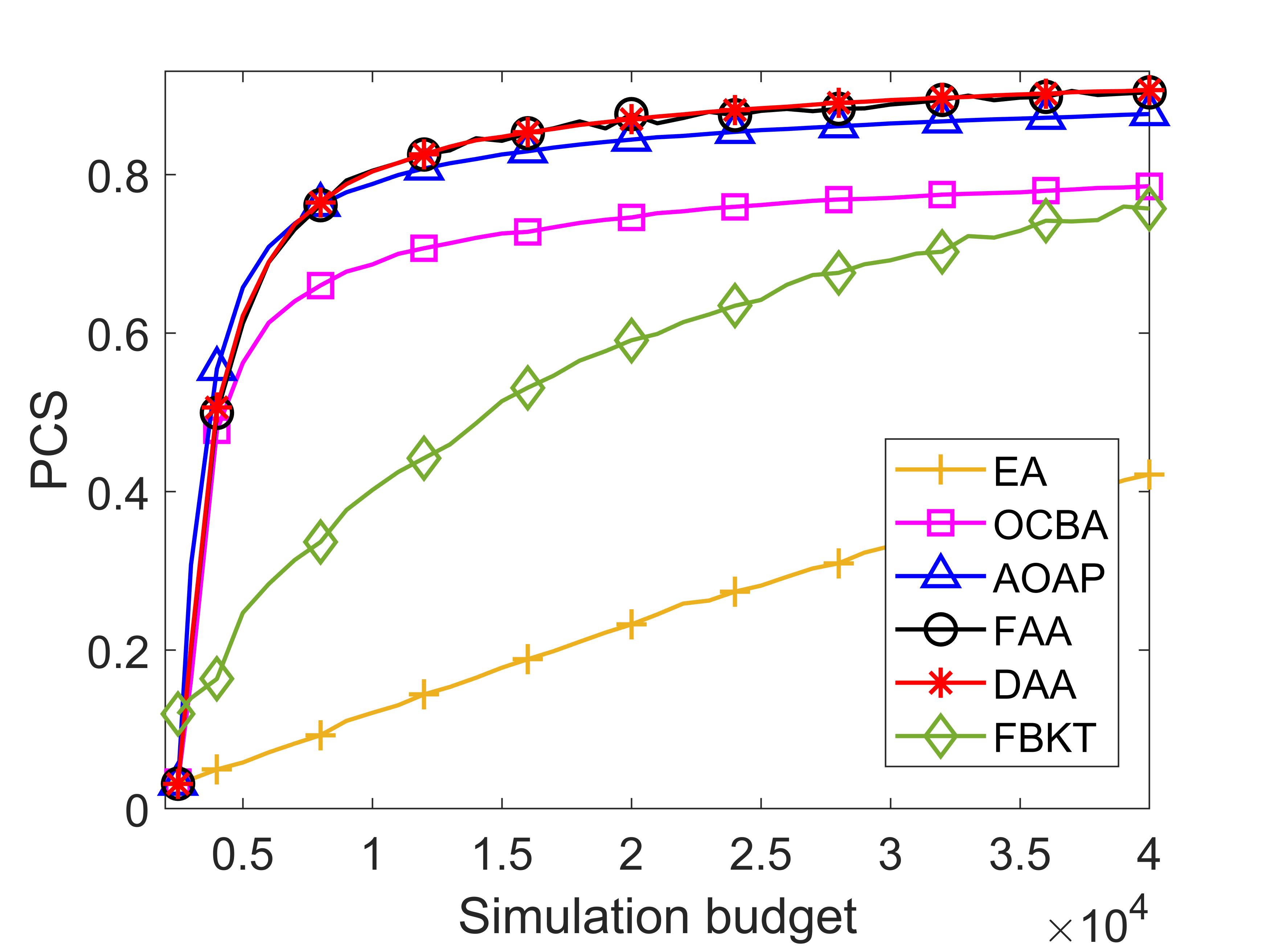}
         \caption{Comparison of PCS for the five competitive procedures on Example 4}
         \label{sfig:Example 4 PCS of the five procedures}
     \end{subfigure}
     \caption{(Color online) Illustration of Example 4}
    \label{fig:Illustration of Example 4}
\end{figure}

\begin{figure}[t]
     \centering
     \begin{subfigure}[b]{0.21\textwidth}
         \centering
         \includegraphics[width=\textwidth]{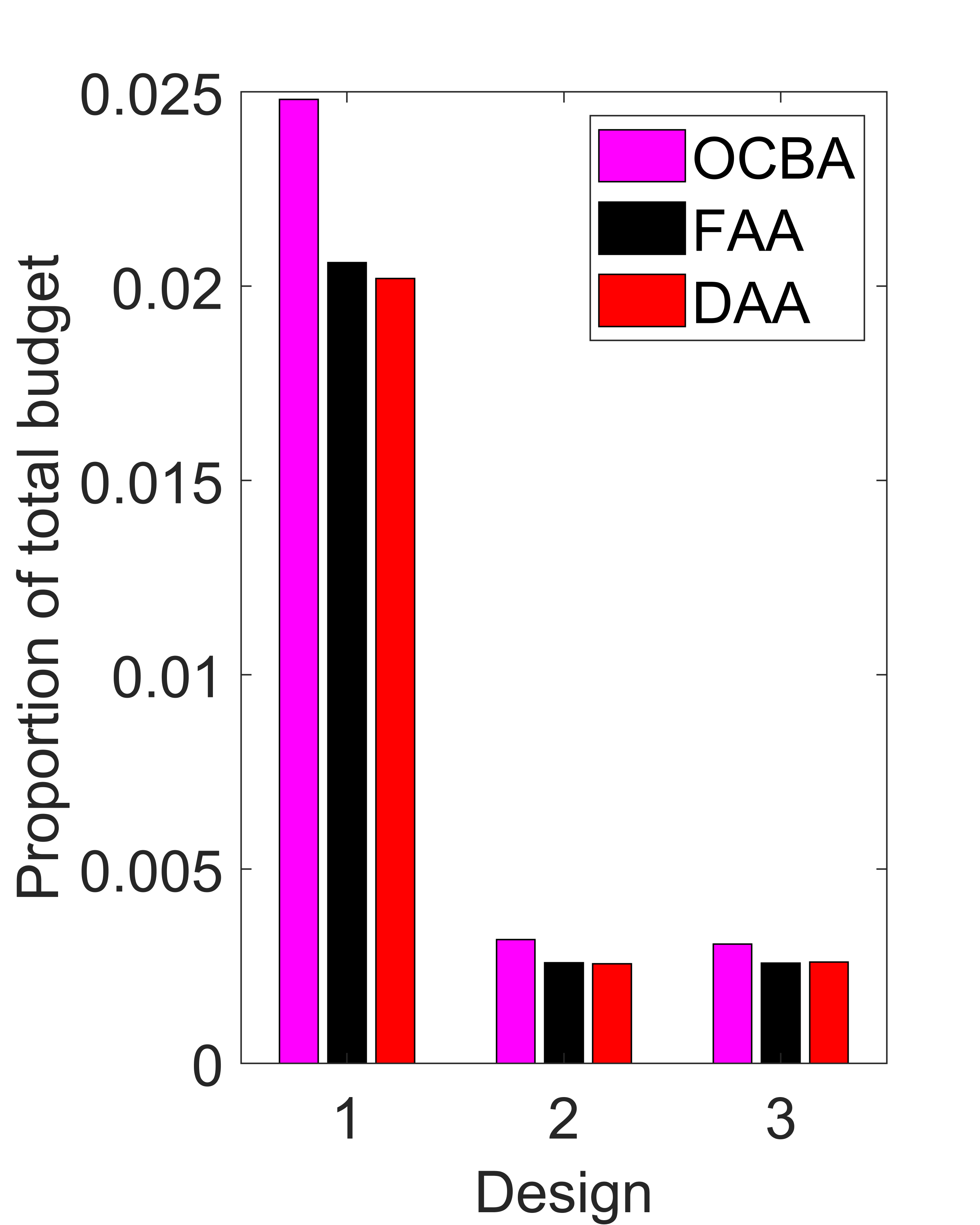}
         \caption{Budget allocation ratios of top 3 designs}
         \label{sfig:Example 5 ratios of top 3 designs}
     \end{subfigure}
     \hfill
     \begin{subfigure}[b]{0.21\textwidth}
         \centering
         \includegraphics[width=\textwidth]{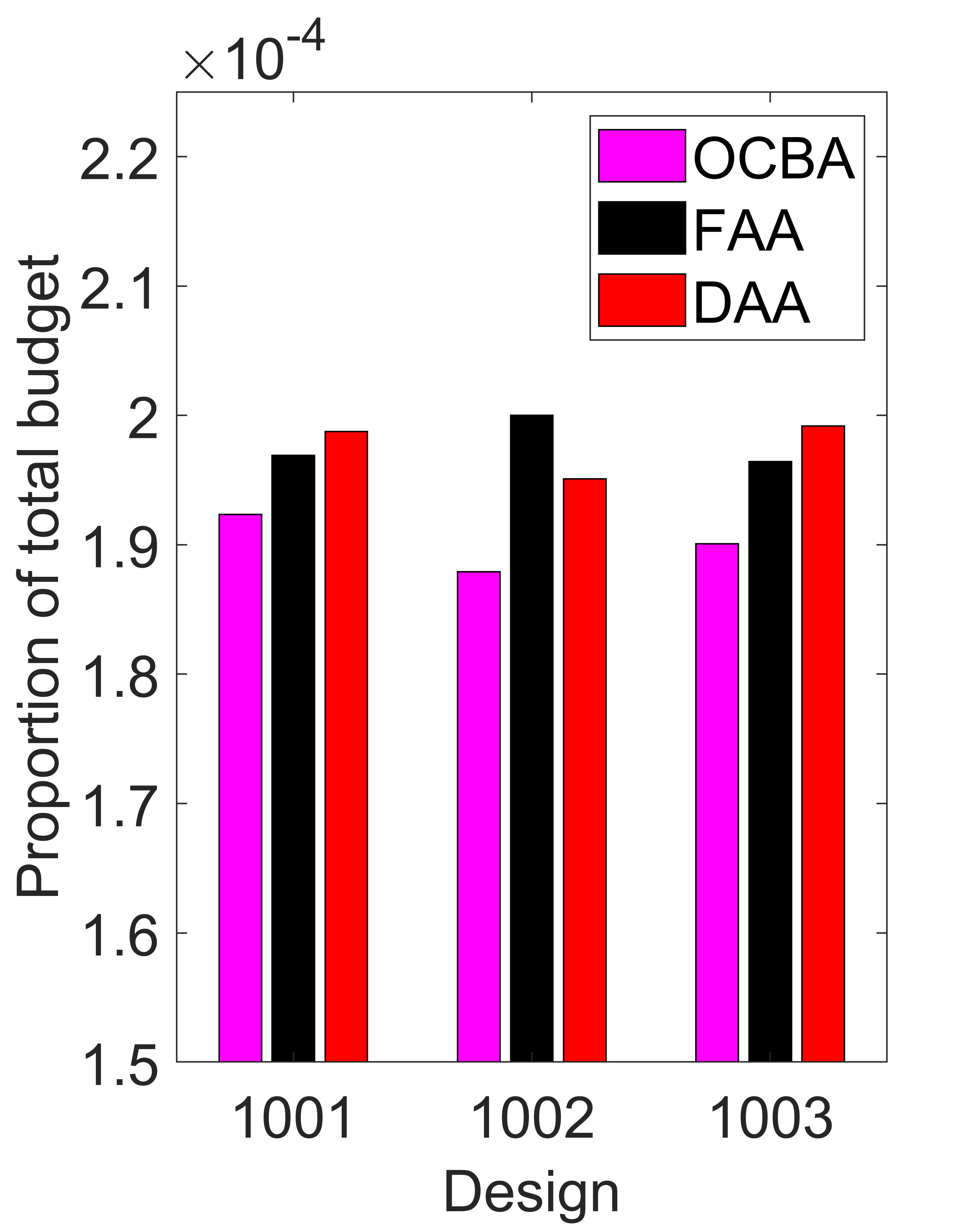}
         \caption{Budget allocation ratios of design 998, 999, and 1000}
         \label{sfig:Example 5 ratios of mid 3 designs}
     \end{subfigure}
     \hfill
     \begin{subfigure}[b]{0.21\textwidth}
         \centering
         \includegraphics[width=\textwidth]{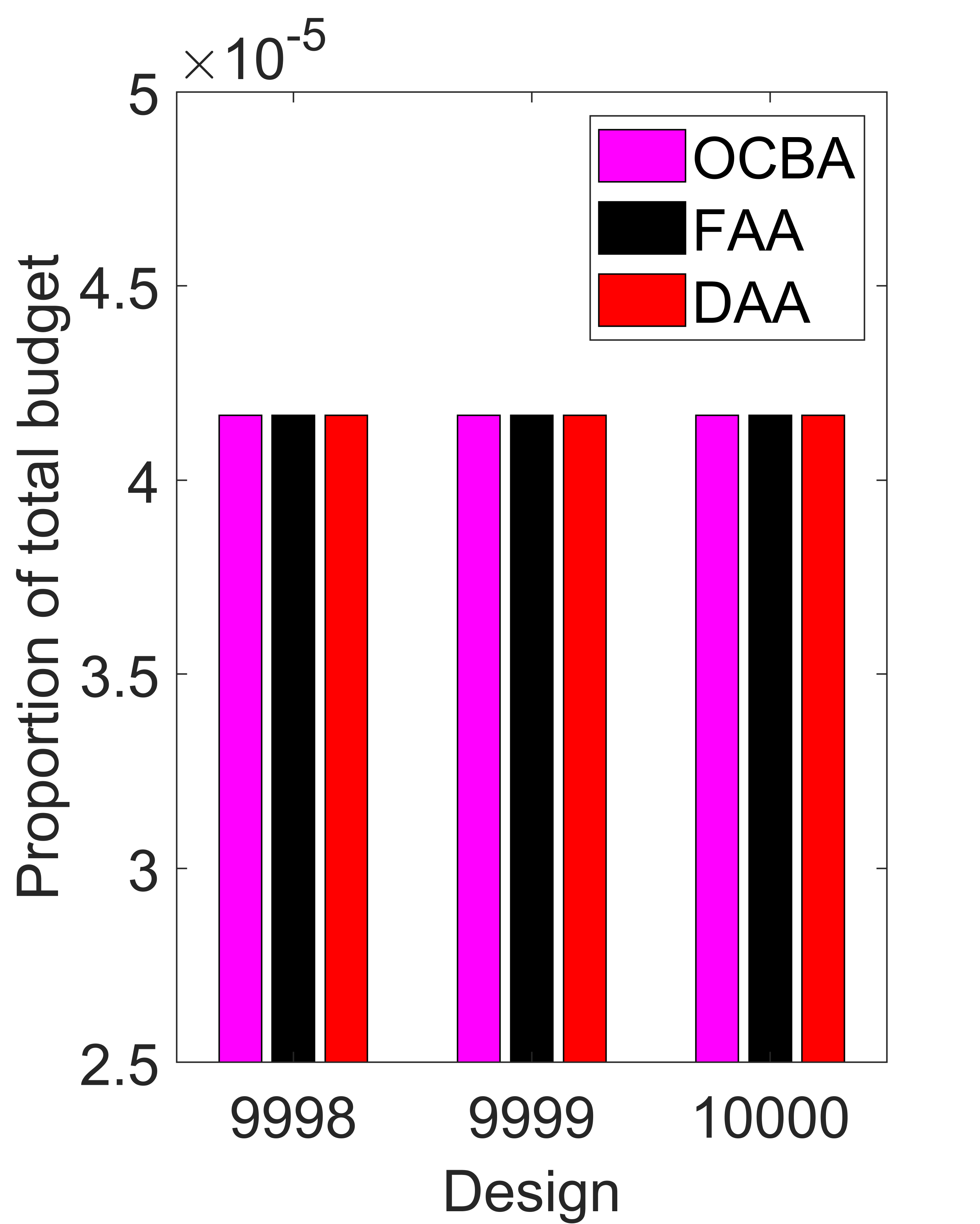}
         \caption{Budget allocation ratios of last 3 designs}
         \label{sfig:Example 5 ratios of last 3 designs}
     \end{subfigure}
     \hfill
     \begin{subfigure}[b]{0.35\textwidth}
         \centering
         \includegraphics[width=\textwidth]{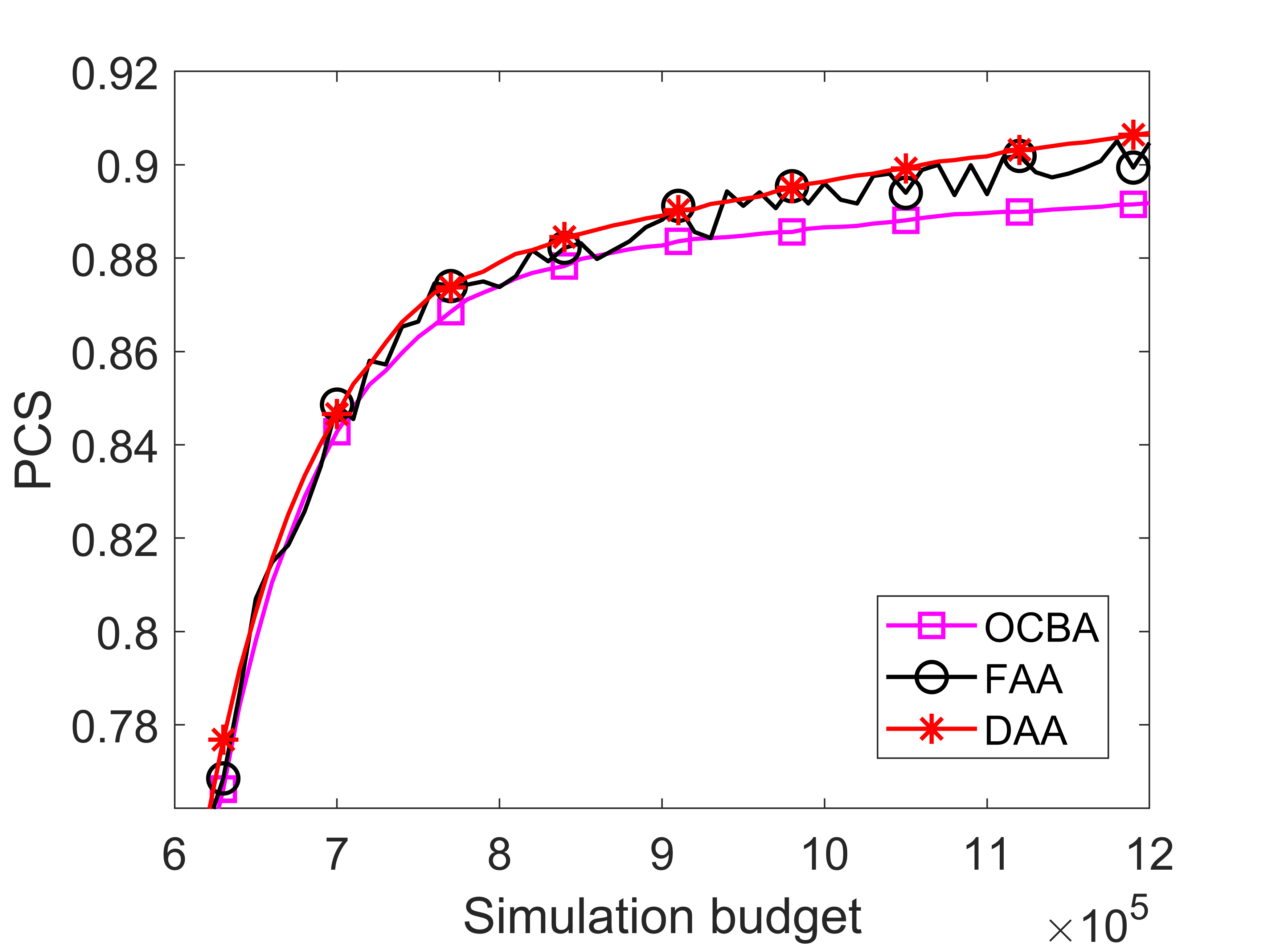}
         \caption{Comparison of PCS for OCBA, FAA, and DAA on Example 5}
         \label{sfig:Example 5 PCS of the five procedures}
     \end{subfigure}
     \caption{(Color online) Illustration of Example 5}
    \label{fig:Illustration of Example 5}
\end{figure}

\begin{figure}[t]
     \centering
     \begin{subfigure}[b]{0.267\textwidth}
         \centering
         \includegraphics[width=\textwidth]{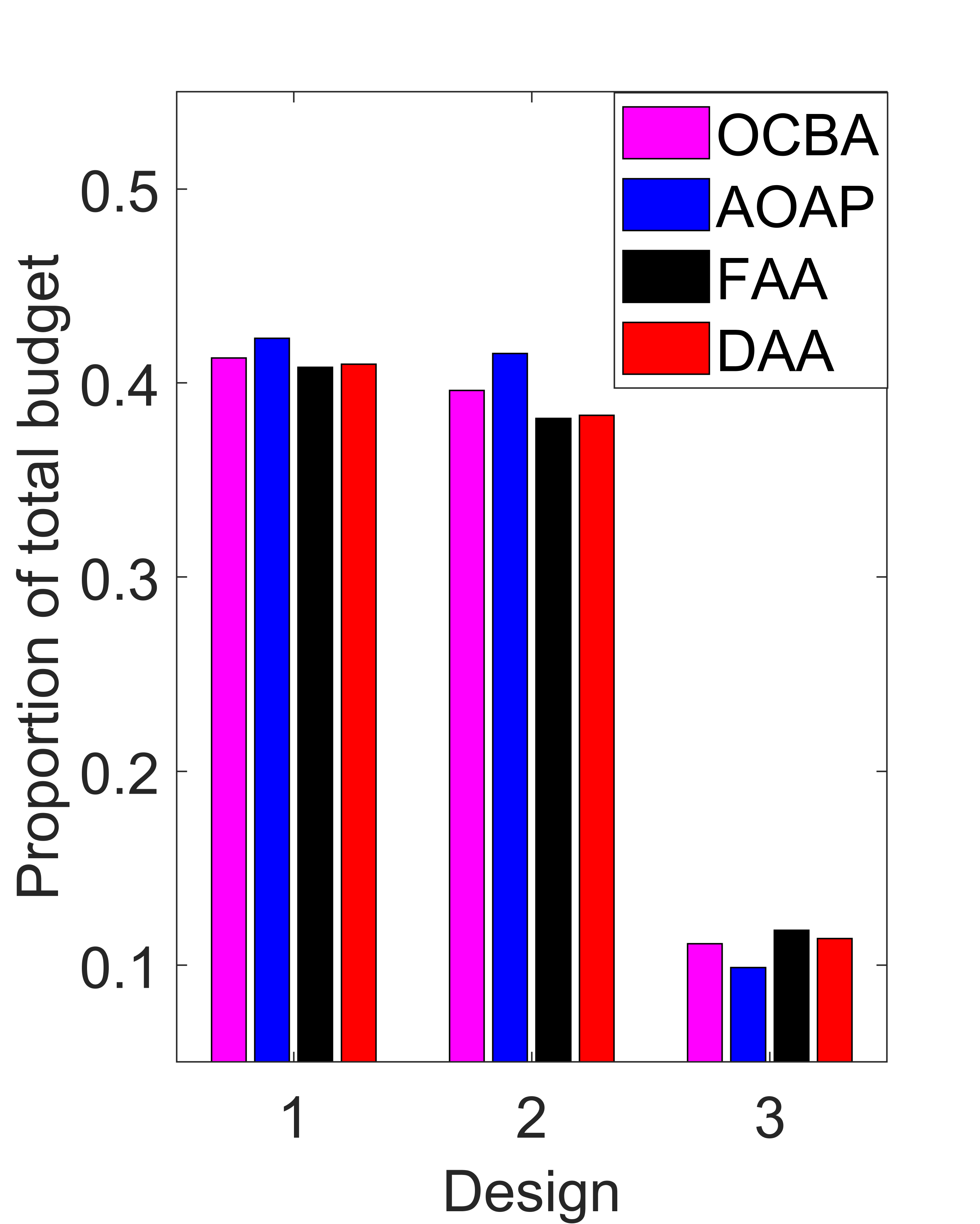}
         \caption{Budget allocation ratios of top 3 designs}
         \label{sfig:FACLOC ratios of top 3 designs}
     \end{subfigure}
     \hfill
     \begin{subfigure}[b]{0.267\textwidth}
         \centering
         \includegraphics[width=\textwidth]{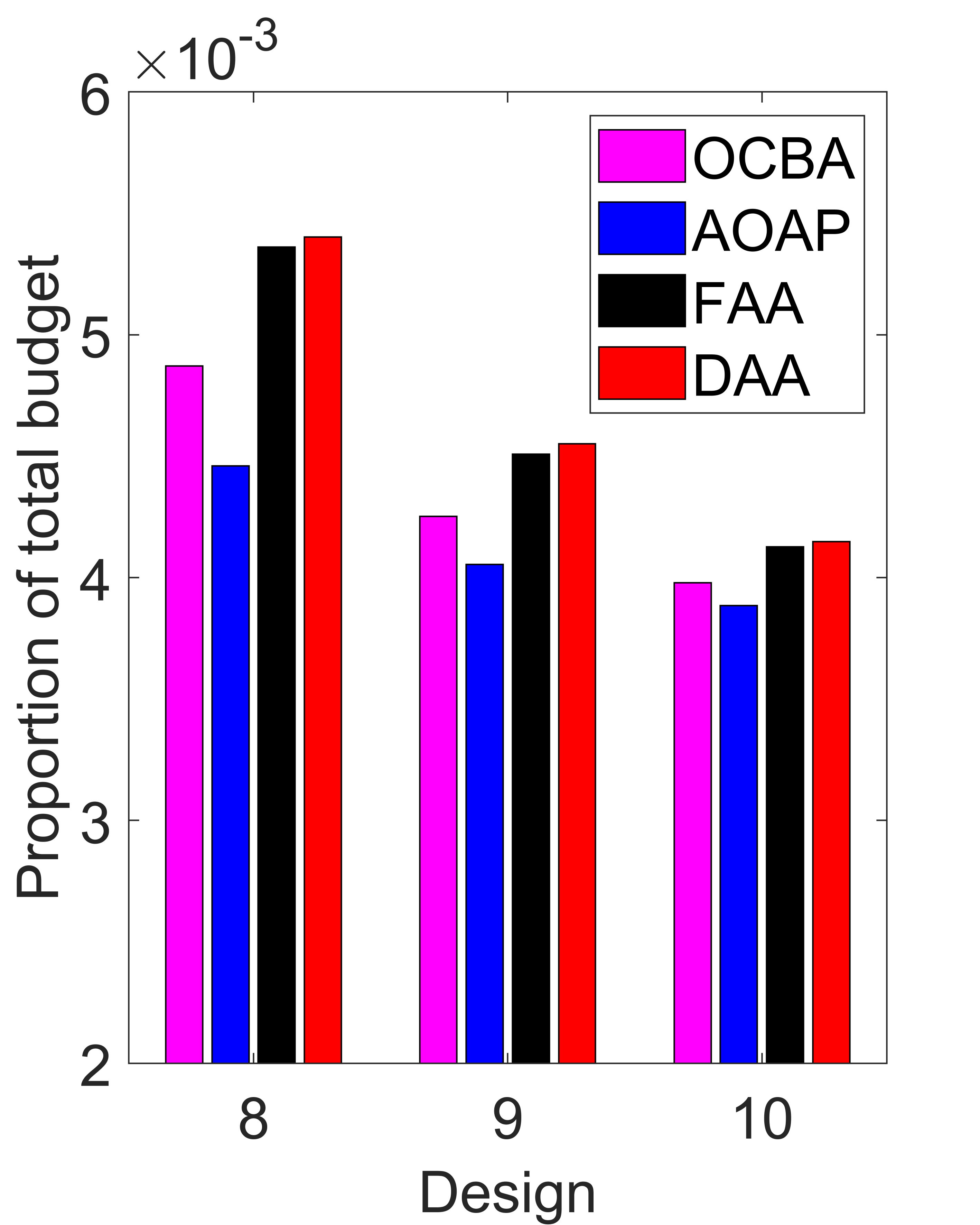}
         \caption{Budget allocation ratios of last 3 designs}
         \label{sfig:FACLOC ratios of last 3 designs}
     \end{subfigure}
     \hfill
     \begin{subfigure}[b]{0.45\textwidth}
         \centering
         \includegraphics[width=\textwidth]{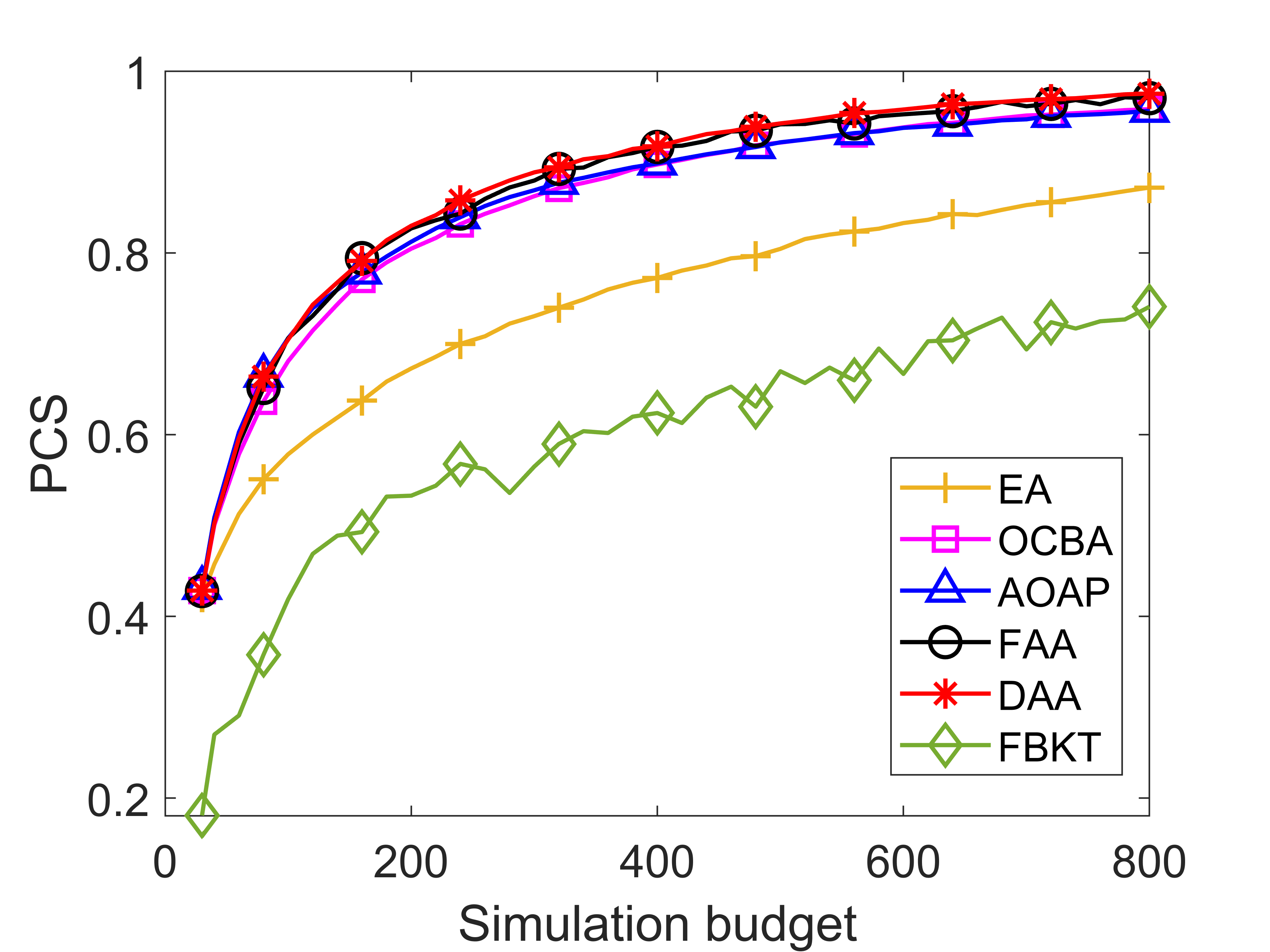}
         \caption{Comparison of PCS for the five competitive procedures on the facility location problem}
         \label{sfig:FACLOC PCS of the five procedures}
     \end{subfigure}
     \caption{(Color online) Illustration of Example 6 (facility location problem)}
    \label{fig:Illustration of the facility location problem}
\end{figure}

The initial number of simulation replications per design for Example 1-3 is set to be 3, i.e., $n_0 = 3$; for Example 4, $n_0 = 5$; for Example 5, $n_0 = 50$; and for Example 6, $n_0 = 3$. The total simulation budgets for Example 1-6 are $10^3$, $3 \times 10^3$, $5 \times 10^3$, $4 \times 10^4$, $1.2 \times 10^6$, and $8 \times 10^2$, respectively. The number of independent macro replications to evaluate empirical PCS for Example 1-3 and 4-6 is $10^5$ and $10^4$, respectively. Figure \ref{fig:Illustration of Example 1}-\ref{fig:Illustration of the facility location problem} and Table \ref{table: performance comparison} show the proportional allocations to designs made by algorithms when budgets are depleted and compare the empirical PCS achieved by all algorithms under different budgets. Allocations made by FBKT are not illustrated because it differs from the others that will converge to certain allocation ratios in the long run. For better visualization, the empirical PCS achieved by EA and FBKT for Example 5 is not shown in Figure \ref{fig:Illustration of Example 5} but reported in Table \ref{table: performance comparison}. Note that, AOAP is not tested for comparison in Example 5 due to its explosively increasing computational burden with the growth of the number of designs (as will be discussed in Section \ref{subsection: discussion on experiment results}).

\subsubsection{Discussion on experiment results} 
\label{subsection: discussion on experiment results}
From Figure~\ref{fig:Illustration of Example 1}-\ref{fig:Illustration of the facility location problem} and Table \ref{table: performance comparison}, both FAA and DAA have better performances than OCBA across all tested examples; and the difference in performance between FAA and DAA is subtle. Both FAA and DAA, compared with OCBA, discount the proportions of total budget allocated to competitive designs (non-best designs that are hard to distinguish from the best), e.g., design 2 in Example 1, while they increase the proportions of total budget allocated to non-competitive designs, e.g., design 10 in Example 1. In Online Figure D.1, the same phenomenon appears in the comparison of proportional allocations to designs made by FAA and OCBA in Example 1 even with different budgets, e.g., between a budget of 100 and 110.  These observations are consistent with our analyses on the potential impact of budget on the optimal budget allocation strategy: compared with the OCBA allocation rule, the budget allocation ratios of competitive designs and non-competitive designs should be discounted and increased based on the budget, respectively. As for the best design, its proportional allocations made by both FAA and DAA are larger than that made by OCBA in Example 2, 4, and 5 but are smaller than that made by OCBA in Example 1, 3, and 6. This verifies that whether the proportion of total budget allocated to the best designs should be discounted or increased really depends on the specific problem structure. In Figure \ref{sfig:Example 5 ratios of last 3 designs}, the proportional allocations to last 3 designs made by OCBA, FAA, and DAA are equal. This is because in the experiment with Example 5, the last 3 designs are clearly inferior to the best design, and they do not receive any sample after initialization. As shown in Table \ref{table: Average runtime of the five allocation procedures (in seconds)}, the average runtimes of FAA and DAA are longer than OCBA as both algorithms require additional computational time to take the impact of budget on allocation ratios into consideration. However, the average runtimes of both FAA and DAA are on the same magnitude as OCBA. Therefore, the numerical results and preceding analyses verify that both FAA and DAA possess desirable budget-adaptive property, and clearly, this desirable property improves their small-budget efficiency of correct selection compared with OCBA.

\begin{table}[thp]
\footnotesize
\setlength\tabcolsep{12pt}
\caption{Performance comparison of the six allocation procedures with different simulation budgets for synthetic examples and the facility location problem}
\label{table: performance comparison}
\renewcommand{\arraystretch}{0.667}
\centering
\begin{tabular}{@{}lccccccc@{}}
\midrule
\multirow{2}{*}{Example 1} & \multicolumn{7}{c}{Simulation budget}                               \\ \cmidrule(l){2-8} 
                           & 50      & 100     & 200     & 400     & 600     & 800     & 1000    \\ \midrule
EA                         & 0.425   & 0.523   & 0.631   & 0.744   & 0.805   & 0.846   & 0.876   \\
OCBA                       & 0.466   & 0.623   & 0.749   & 0.856   & 0.906   & 0.934   & 0.950   \\ 
AOAP                       & 0.492   & 0.643   & 0.760   & 0.857   & 0.902   & 0.928   & 0.943   \\ 
FAA                        & 0.474   & 0.631   & 0.771   & 0.881   & 0.930   & 0.954   & 0.967   \\ 
DAA                        & 0.473   & 0.631   & 0.771   & 0.886   & 0.934   & 0.957   & 0.969   \\ 
FBKT                       & 0.245   & 0.366   & 0.466   & 0.571   & 0.636   & 0.680   & 0.722   \\
\midrule
\multirow{2}{*}{Example 2} & \multicolumn{7}{c}{Simulation budget}                               \\ \cmidrule(l){2-8} 
                           & 50      & 150     & 500     & 1000    & 1500    & 2000    & 3000    \\ \midrule
EA                         & 0.389   & 0.505   & 0.654   & 0.753   & 0.811   & 0.850   & 0.900   \\ 
OCBA                       & 0.388   & 0.571   & 0.760   & 0.858   & 0.906   & 0.933   & 0.959   \\ 
AOAP                       & 0.404   & 0.583   & 0.751   & 0.844   & 0.892   & 0.919   & 0.949   \\ 
FAA                        & 0.398   & 0.589   & 0.789   & 0.890   & 0.935   & 0.955   & 0.974   \\ 
DAA                        & 0.396   & 0.586   & 0.792   & 0.895   & 0.938   & 0.958   & 0.976   \\ 
FBKT                       & 0.215   & 0.331   & 0.489   & 0.594   & 0.652   & 0.697   & 0.755   \\
\midrule
\multirow{2}{*}{Example 3} & \multicolumn{7}{c}{Simulation budget}                               \\ \cmidrule(l){2-8} 
                           & 200     & 500     & 800     & 1000    & 2000    & 3000    & 5000    \\ \midrule
EA                         & 0.281   & 0.382   & 0.443   & 0.375   & 0.581   & 0.643   & 0.725   \\ 
OCBA                       & 0.356   & 0.635   & 0.724   & 0.762   & 0.864   & 0.907   & 0.947   \\ 
AOAP                       & 0.429   & 0.672   & 0.755   & 0.791   & 0.886   & 0.924   & 0.955   \\ 
FAA                        & 0.383   & 0.677   & 0.775   & 0.814   & 0.912   & 0.945   & 0.970   \\ 
DAA                        & 0.382   & 0.679   & 0.782   & 0.822   & 0.920   & 0.953   & 0.974   \\
FBKT                       & 0.221   & 0.372   & 0.432   & 0.473   & 0.580   & 0.645   & 0.726   \\
\midrule
\multirow{2}{*}{Example 4} & \multicolumn{7}{c}{Simulation budget ($ \times 10^3$)}             \\ \cmidrule(l){2-8} 
                           & 6       & 8       & 10      & 15      & 20      & 30      & 40        \\ \midrule
EA                         & 0.071   & 0.093   & 0.121   & 0.178   & 0.233   & 0.331   & 0.422     \\ 
OCBA                       & 0.613   & 0.660   & 0.687   & 0.726   & 0.746   & 0.771   & 0.785     \\ 
AOAP                       & 0.709   & 0.762   & 0.788   & 0.825   & 0.844   & 0.864   & 0.876     \\ 
FAA                        & 0.689   & 0.761   & 0.805   & 0.843   & 0.876   & 0.888   & 0.904     \\ 
DAA                        & 0.690   & 0.765   & 0.804   & 0.848   & 0.870   & 0.894   & 0.907     \\ 
FBKT                       & 0.284   & 0.337   & 0.402   & 0.514   & 0.591   & 0.692   & 0.757   \\
\midrule
\multirow{2}{*}{Example 5} & \multicolumn{7}{c}{Simulation budget ($ \times 10^5$)}              \\ \cmidrule(l){2-8} 
                           & 6       & 7       & 8       & 9       & 10      & 11      & 12        \\ \midrule
EA                         & 0.035   & 0.046   & 0.058   & 0.070   & 0.091   & 0.105   & 0.124     \\ 
OCBA                       & 0.706   & 0.843   & 0.874   & 0.883   & 0.887   & 0.890   & 0.892     \\ 
AOAP                       & -       & -       & -       & -       & -       & -       & -         \\ 
FAA                        & 0.719   & 0.849   & 0.874   & 0.888   & 0.896   & 0.894   & 0.905     \\ 
DAA                        & 0.722   & 0.847   & 0.879   & 0.889   & 0.896   & 0.902   & 0.907     \\
FBKT                       & 0.685   & 0.724   & 0.749   & 0.789   & 0.806   & 0.830   & 0.837   \\
\midrule
\multirow{2}{*}{\makecell[l]{Example 6 \\(Facility location problem)}} & \multicolumn{7}{c}{Simulation budget}   \\ \cmidrule(l){2-8} 
                           & 40      & 120     & 200     & 300     & 500     &700      & 800       \\ \midrule
EA                         & 0.458   & 0.600   & 0.673   & 0.731   & 0.805   & 0.853   & 0.872     \\ 
OCBA                       & 0.501   & 0.715   & 0.805   & 0.863   & 0.922   & 0.952   & 0.959     \\ 
AOAP                       & 0.509   & 0.740   & 0.813   & 0.869   & 0.922   & 0.947   & 0.957     \\ 
FAA                        & 0.504   & 0.731   & 0.827   & 0.880   & 0.942   & 0.962   & 0.971     \\ 
DAA                        & 0.502   & 0.743   & 0.830   & 0.889   & 0.943   & 0.968   & 0.976     \\ 
FBKT                       & 0.270   & 0.469   & 0.533   & 0.565   & 0.670   & 0.694   & 0.741     \\
\bottomrule
\end{tabular}
\end{table}

\begin{table}[thp]
\footnotesize
\setlength\tabcolsep{14pt}
\caption{Average runtime of the six allocation procedures (in seconds)}
\label{table: Average runtime of the five allocation procedures (in seconds)}
\centering
\begin{tabular}{@{}lcccccc@{}}
\midrule
                        & EA      & OCBA     & AOAP     & FAA      & DAA       & FBKT   \\ \hline
Example 1               & 0.001   & 0.006    & 0.015    & 0.012    & 0.012     & 0.001   \\ 
Example 2               & 0.001   & 0.003    & 0.007    & 0.006    & 0.006     & 0.001   \\ 
Example 3               & 0.005   & 0.030    & 0.335    & 0.057    & 0.058     & 0.004   \\ 
Example 4               & 0.061   & 1.409    & 59.270   & 3.352    & 3.436     & 0.038   \\
Example 5               & 3.371   & 250.087  & -        & 497.449  & 499.945   & 1.218   \\
\makecell[l]{Example 6} & 147.132 & 151.365  & 150.298  & 149.878  & 150.663   & 66.661 \\
\bottomrule
\end{tabular}
\end{table}

For small- to medium-scale problems, i.e., Example 1-4, FAA and DAA performs slightly worse than AOAP at the beginning, but they then surpass AOAP when budget is relatively large. One possible reason why it is the case is that AOAP is a myopic procedure that aims to maximize one-step-ahead improvement, and this myopic approach leads to good performance at the very beginning. Since AOAP ignores the impact of budget on allocation strategy, and it is reasonable for FAA and DAA to gradually outperform AOAP. As shown in Table \ref{table: Average runtime of the five allocation procedures (in seconds)}, the average runtimes of AOAP are compatible with FAA and DAA for small-scale problems, e.g., in Example 1-2. However, as the scale of problem becomes large, e.g., in Example 3-4, the average runtimes of AOAP drastically increase and are much longer than FAA and DAA. This indicates that both FAA and DAA are computationally more efficient than AOAP for large-scale problems. But for Example 6, the average runtime of every algorithm except FBKT is almost the same. FBKT often terminates before the budget is exhausted, and thus its average runtimes are shorter than the other algorithms for Example 6. These results imply that when the complexity of the simulated system is high and the simulation time is relatively long, the additional runtimes for calculating $\alpha_i(T)$ in both FAA and DAA, as well as the computational concerns on AOAP, can be negligible. Note that, the real industrial systems can be much more complex than the logistic system considered in this paper. Therefore, FAA, DAA, and AOAP, are competitive and possess good applicability in real industrial applications. FBKT is dominated by FAA and DAA for all tested examples; it does not perform well for small-scale problems; but it tends to perform better with the increase of the scale of problem. This is reasonable because FBKT is specially designed for large-scale problems. Overall, FAA and DAA are efficient algorithms for both small- and large-scale problems.

\section{Conclusion}
\label{conclusion}
In this paper, we consider a simulation-based R\&S problem of identifying the best system design from a set of alternatives under a fixed budget setting. We propose a budget-adaptive rule under the setting where simulation budget is not large enough to simplify computations, the setting of which differs from the derivation of the asymptotic OCBA rule in \cite{chen2000simulation}. Based on the proposed budget-adaptive allocation rule, two heuristic algorithms FAA and DAA are developed. As iteration increases, both algorithms can learn the best design with probability 1 and converge to the asymptotic OCBA rule. In the numerical experiments, various settings are tested to demonstrate the budget-adaptive property of FAA and DAA improves their small-budget performances compared with the asymptotic OCBA rule. 

We highlight the most important implication of our contributions: simulation budget significantly impacts the budget allocation strategy, and a desirable budget allocation rule should be adaptive to the simulation budget. The proposed budget-adaptive allocation rule indicates that, compared with the asymptotic OCBA rule, the budget allocation ratios of non-best designs that are hard to distinguish from the best design should be discounted, while the budget allocation ratios of those that are easy to distinguish from the best should be increased. These adjustments are based on the simulation budget, which is often limited in practice. Therefore, the budget-adaptive rule highlights the significant impacts of budget size on budget allocation strategy. We believe these findings can help and motivate researchers to develop more efficient budget allocation rules in future studies.

Many opportunities for future work remain. First, deriving an allocation rule that achieves the true optimality remains significant future directions. Second, the budget allocation problem can be essentially formulated as a stochastic dynamic program (DP) problem. Both FAA and DAA are efficient procedures that can adapt to the simulation budget, but they ignore the dynamic feedback of the final step while sampling at the current step. Recently, \cite{qin2022non} in their preliminary version of work formulate the problem as a DP and investigate a non-myopic knowledge gradient (KG) procedure, which can dynamically look multiple steps ahead and take the dynamic feedback mechanism into consideration. However, exactly solving the DP is intractable due to the extremely high computational cost caused by “curse of dimensionality”. As a result, how to derive a computationally tractable allocation rule that can incorporate the dynamic feedback mechanism, remains a critical future direction.  

\section*{Acknowledgments}
This research paper has been made possible by the funding support from the Singapore Maritime Institute \& the Center of Excellence in Modelling and Simulation for Next Generation Ports (Singapore Maritime Institute grant: SMI-2022-SP-02).

\appendix

\begingroup
\allowdisplaybreaks

\section{Definitions of Big-$\mathcal{O}$, Big-$\Theta$, and Little-$\omega$ notations}

According to \cite{cormen2022introduction}, we present the definitions of Big-$\Omega$, Big-$\mathcal{O}$, Big-$\Theta$, and Little-$\omega$ notations as follows.\\

\textit{Definition 1. For a given function $g(x)$, we denote by $\mathcal{O}(g(x))$ the set of functions
\begin{equation*}
\begin{split}
    \mathcal{O}(g(x)) = \{ f(x):\ &\text{there exists positive constants} \ M \ \text{and} \ x_0 \ \text{such that} \\ 
    &0 \leq f(x) \leq M \cdot g(x),\  \forall x \geq x_0 \}. 
\end{split}
\end{equation*}}

\textit{Definition 2. For any two functions $f(x)$ and $g(x)$, we have $f(x) = \Theta(g(x))$ if and only if $f(x) = \Omega(g(x))$ and $f(x) = \mathcal{O}(g(x))$.}

\textit{Definition 3. For a given function $g(x)$, we denote by $\omega(g(x))$ the set of functions
\begin{equation*}
\begin{split}
    \mathcal{O}(g(x)) = \{ f(x):\ &\text{there exists positive constants} \ M \ \text{and} \ x_0 \ \text{such that} \\ 
    &0 \leq M \cdot g(x) < f(x),\  \forall x \geq x_0 \}. 
\end{split}
\end{equation*}}
 
The notation $f(x) = \mathcal{O}(g(x))$ means that $g(x)$ can be viewed as the upper bound of the growth rate of $f(x)$. The notation $f(x) = \Theta(g(x))$ means that $f(x)$ grows in the exact same rate as that of $g(x)$. The notation $f(x) = \omega(g(x))$ means that $g(x)$ can be viewed as the lower bound of the growth rate of $f(x)$, but it is not a tight bound.

\section{Proof of Lemmas}

\subsection{Proof of Lemma 1}
\label{Proof of Lemma 1}
 
The constraints of Problem $\mathcal{P}1$ are all affine functions of $w_i$, for $i \in \mathcal{K}$. Furthermore, showing $\text{APCS}$ is concave is equivalent to showing
\begin{equation*}
    g(w) = \sum_{i \in \mathcal{K}^{\prime}}  \Phi \left( - \frac{\delta_{i,b}}{\sigma_{i,b}} \right),
\end{equation*}
is a convex function of $w$. To verify the convexity of $g(w)$, we need to show its Hessian matrix is positive semi-definite. The Hessian matrix for $g(w)$ is
\begin{equation*}
\nabla^2 g(w) = 
\begin{pmatrix}
\frac{\partial^{2} g(w)}{\partial^2 w_1} & \frac{\partial^{2} g(w)}{\partial w_1 \partial w_2} & \cdots & \frac{\partial^{2} g_(w)}{\partial w_1 \partial w_k} \\
\frac{\partial^{2} g(w)}{\partial w_2 \partial w_1} & \frac{\partial^{2} g(w)}{\partial^2 w_2} & \cdots & \frac{\partial^{2} g(w)}{\partial w_2 \partial w_k} \\
\vdots  & \vdots  & \ddots & \vdots  \\
\frac{\partial^{2} g(w)}{\partial w_k \partial w_1} & \frac{\partial^{2} g(w)}{\partial w_k \partial w_2} & \cdots & \frac{\partial^{2} g(w)}{\partial^2 w_k}
\end{pmatrix},
\end{equation*}
where, for $i,j \in \mathcal{K}^{\prime}$ and $i \neq j$, and $b$
\begin{align*}
    \frac{\partial^{2} g(w)}{\partial w_i^{2}}&= \frac{1}{2 \sqrt{2 \pi}} \exp \left(-\frac{\delta_{i,b}^{2}}{2 \sigma_{i,b}^{2}}\right) \times\left[\frac{\delta_{i,b}}{\sigma_{i,b}^{5} T^{2}} \frac{\sigma_{i}^{4}}{w_i^{4}}\left(\frac{1}{2} \frac{\delta_{i,b}^{2}}{\sigma_{i,b}^{2}}-\frac{3}{2}\right)+2 \frac{\delta_{i,b}}{\sigma_{i,b}^{3} T} \frac{\sigma_{i}^{2}}{w_i^{3}}\right], \\
    \frac{\partial^{2} g(w)}{\partial w_b^{2}}&= \frac{1}{2 \sqrt{2 \pi}} \sum\limits_{i \in \mathcal{K}^{\prime}} \exp \left(-\frac{\delta_{i,b}^{2}}{2 \sigma_{i,b}^{2}}\right) \times\left[\frac{\delta_{i,b}}{\sigma_{i,b}^{5} T^{2}} \frac{\sigma_{b}^{4}}{w_b^{4}}\left(\frac{1}{2} \frac{\delta_{i,b}^{2}}{\sigma_{i,b}^{2}}-\frac{3}{2}\right)+2 \frac{\delta_{i,b}}{\sigma_{i,b}^{3} T} \frac{\sigma_{b}^{2}}{w_b^{3}}\right], \\
    \frac{\partial^{2} g(w)}{\partial w_i \partial w_b}&= \frac{1}{2 \sqrt{2 \pi}} \exp \left(-\frac{\delta_{i,b}^{2}}{2 \sigma_{i,b}^{2}}\right) \times\left[\frac{\delta_{i,b}}{\sigma_{i,b}^{5} T^{2}} \frac{\sigma_{i}^{2}}{w_i^{2}} \frac{\sigma_{b}^{2}}{w_b^{2}}\left(\frac{1}{2} \frac{\delta_{i,b}^{2}}{\sigma_{i,b}^{2}}-\frac{3}{2}\right)\right],\\
\frac{\partial^{2} g(w)}{\partial w_i \partial w_{j}} &= 0.
\end{align*}
Furthermore, for any non-zero vector $a\in \mathbb{R}^k$, we have
\begin{align*}
    a^{T} \nabla^2 g(w) a &= \sum\limits_{i \in \mathcal{K}^{\prime}} \frac{\partial^{2} g(w)}{\partial w_i^2} a_i^2 +  \frac{\partial^{2} g(w)}{\partial w_b^{2}} a_b^2  + 2\sum\limits_{i \in \mathcal{K}^{\prime}} \frac{\partial^{2} g(w)}{\partial w_i \partial w_b} a_i a_b\\
    &= \frac{1}{2 \sqrt{2 \pi}} \sum\limits_{i \in \mathcal{K}^{\prime}} \exp \left(-\frac{\delta_{i,b}^{2}}{2 \sigma_{i,b}^{2}}\right)\times a_i^2 \left[\frac{\delta_{i,b}}{\sigma_{i,b}^{5} T^{2}} \frac{\sigma_{i}^{4}}{w_i^{4}}\left(\frac{1}{2} \frac{\delta_{i,b}^{2}}{\sigma_{i,b}^{2}}-\frac{3}{2}\right)+2 \frac{\delta_{i,b}}{\sigma_{i,b}^{3} T} \frac{\sigma_{i}^{2}}{w_i^{3}}\right]\\
    &\quad + \frac{1}{2 \sqrt{2 \pi}} \sum\limits_{i \in \mathcal{K}^{\prime}} \exp \left(-\frac{\delta_{i,b}^{2}}{2 \sigma_{i,b}^{2}}\right) \times a_b^2 \left[\frac{\delta_{i,b}}{\sigma_{i,b}^{5} T^{2}} \frac{\sigma_{b}^{4}}{w_b^{4}}\left(\frac{1}{2} \frac{\delta_{i,b}^{2}}{\sigma_{i,b}^{2}}-\frac{3}{2}\right)+2 \frac{\delta_{i,b}}{\sigma_{i,b}^{3} T} \frac{\sigma_{b}^{2}}{w_b^{3}}\right]\\
    &\quad + \frac{1}{2 \sqrt{2 \pi}} \sum\limits_{i \in \mathcal{K}^{\prime}} \exp \left(-\frac{\delta_{i,b}^{2}}{2 \sigma_{i,b}^{2}}\right) \times 2 a_i a_b \left[\frac{\delta_{i,b}}{\sigma_{i,b}^{5} T^{2}} \frac{\sigma_{i}^{2}}{w_i^{2}} \frac{\sigma_{b}^{2}}{w_b^{2}}\left(\frac{1}{2} \frac{\delta_{i,b}^{2}}{\sigma_{i,b}^{2}}-\frac{3}{2}\right)\right] \\
    &= \frac{1}{2 \sqrt{2 \pi}} \sum\limits_{i \in \mathcal{K}^{\prime}} \exp \left(-\frac{\delta_{i,b}^{2}}{2 \sigma_{i,b}^{2}}\right) \\
    & \quad \times \frac{\delta_{i,b}}{\sigma_{i,b}^{3} T}  \left[\frac{1}{\sigma_{i,b}^{2} T} \left(\frac{1}{2} \frac{\delta_{i,b}^{2}}{\sigma_{i,b}^{2}}-\frac{3}{2}\right) \left(\frac{\sigma_{i}^{2}}{w_i^{2}} a_i + \frac{\sigma_{b}^{2}}{w_b^{2}} a_b \right)^2  + 2 \left( \frac{\sigma_{i}^{2}}{w_i^{3}} a_i^2 + \frac{\sigma_{b}^{2}}{w_b^{3}} a_b^2 \right) \right].  \\
\end{align*}
Since
\begin{align*}
    2 \left( \frac{\sigma_{i}^{2}}{w_i^{3}} a_i^2 + \frac{\sigma_{b}^{2}}{w_b^{3}} a_b^2 \right) & = 2 \frac{\sigma_{i,b}^2}{\sigma_{i,b}^2} \left( \frac{\sigma_{i}^{2}}{w_i^{3}} a_i^2 + \frac{\sigma_{b}^{2}}{w_b^{3}} a_b^2 \right)\\
    & = \frac{2}{\sigma_{i,b}^2 T}  \left(\frac{\sigma_i^2}{w_i} + \frac{\sigma_b^2}{w_b} \right) \left( \frac{\sigma_{i}^{2}}{w_i^{3}} a_i^2 + \frac{\sigma_{b}^{2}}{w_b^{3}} a_b^2 \right)\\
    & = \frac{2}{\sigma_{i,b}^2 T} \left[\left(\frac{\sigma_{i}^{2}}{w_i^{2}} a_i + \frac{\sigma_{b}^{2}}{w_b^{2}} a_b \right)^2 + \frac{\sigma_i^2 \sigma_b^2}{w_i w_b} \left( \frac{a_i}{w_i} - \frac{a_b}{w_b} \right)^2\right], \\
\end{align*}
we have
\begin{align*}
    &\frac{1}{\sigma_{i,b}^{2} T} \left(\frac{1}{2} \frac{\delta_{i,b}^{2}}{\sigma_{i,b}^{2}}-\frac{3}{2}\right) \left(\frac{\sigma_{i}^{2}}{w_i^{2}} a_i+ \frac{\sigma_{b}^{2}}{w_b^{2}} a_b \right)^2  + 2 \left( \frac{\sigma_{i}^{2}}{w_i^{3}} a_i^2 + \frac{\sigma_{b}^{2}}{w_b^{3}} a_b^2 \right)  \\
    = & \frac{1}{\sigma_{i,b}^{2} T} \left[ \left(\frac{1}{2} \frac{\delta_{i,b}^{2}}{\sigma_{i,b}^{2}}-\frac{3}{2}\right) \left(\frac{\sigma_{i}^{2}}{w_i^{2}} a_i + \frac{\sigma_{b}^{2}}{w_b^{2}} a_b \right)^2 + 2 \left(\frac{\sigma_{i}^{2}}{w_i^{2}} a_i +  \frac{\sigma_{b}^{2}}{w_b^{2}} a_b \right)^2 + 2 \frac{\sigma_i^2 \sigma_b^2}{w_i w_b} \left( \frac{a_i}{w_i} - \frac{a_b}{w_b} \right)^2\right] \\
    = & \frac{1}{\sigma_{i,b}^{2} T} \left[ \left(\frac{1}{2} \frac{\delta_{i,b}^{2}}{\sigma_{i,b}^{2}}+\frac{1}{2}\right) \left(\frac{\sigma_{i}^{2}}{w_i^{2}} a_i + \frac{\sigma_{b}^{2}}{w_b^{2}} a_b \right)^2 + 2 \frac{\sigma_i^2 \sigma_b^2}{w_i w_b} \left( \frac{a_i}{w_i} - \frac{a_b}{w_b} \right)^2\right]\\
    \geq & 0,
\end{align*}
and therefore, $\nabla^2 g(w) \succeq 0$, and $g(w)$ is a convex function of $w$. Due to the two constraints $\sum_{i \in \mathcal{K}} w_i = 1$ and $w_i \geq 0$, for $i \in \mathcal{K}$, forming a convex set, Problem $\mathcal{P}1$ is a convex optimization problem. This result concludes the proof.


\subsection{Proof of Lemma 3}

By $\widehat{\mathcal{C}}_2$, for $i \in \mathcal{K}^{\prime}$, we have
\begin{equation}
\label{proof p1 w_i}
\begin{split}
    w_i &= \frac{\lambda - 2 \log I_i}{\frac{T}{I_i} + \frac{1}{w^*_i}} \\
    &= \frac{I_i (\lambda - 2 \log I_i)}{T+S} ,
\end{split}
\end{equation}
in which, $S = \sum_{i \in \mathcal{K}} I_i$ and $w^*_i = I_i / S$. We substitute $w_i$ provided by \eqref{proof p1 w_i} into $\mathcal{C}_1$, and obtain
\begin{equation}
\label{proof p1 w_b}
    w_b = \sigma_b \sqrt{ \sum\limits_{i \in \mathcal{K}^{\prime}} \frac{I_i^2(\lambda - 2 \log I_i )^2}{\sigma_i^2 (T+S)^2} }.
\end{equation}
For $i \in \mathcal{K}^{\prime}$, by $\mathcal{C}_3$, \eqref{proof p1 w_i}, and \eqref{proof p1 w_b}
\begin{equation}
\label{verify  lambda feasibility}
    \sum\limits_{i \in \mathcal{K}^{\prime}} \frac{I_i (\lambda - 2 \log I_i)}{T+S} + \sigma_b \sqrt{ \sum\limits_{i \in \mathcal{K}^{\prime}} \frac{I_i^2(\lambda - 2 \log I_i )^2}{\sigma_i^2 \left(T+S\right)^2} } = 1.
\end{equation}
Then, we have
\begin{equation*}
    \lambda =\left\{\begin{array}{ll}
\frac{- q + \sqrt{q^2 - 4 p r}}{2 p}  & \quad \text{if} \  w^*_b \neq \frac{1}{2}  \\
\frac{4\sum_{i \in \mathcal{K}^{\prime}} I_i \log I_i + T + S}{2\sum_{i \in \mathcal{K}^{\prime}} I_i} & \quad \text{if} \ w^*_b = \frac{1}{2} \\
\end{array}\right.
\end{equation*}
where
\begin{equation*}
\begin{split}
    p  &=  S(2 I_b - S), \\
    q  &= - 4 \sigma_b^2  \sum\limits_{i \in \mathcal{K}^{\prime}} \frac{I_i^2 \log I_i }{\sigma_i^2 } + 2 (S- I_b)  \left(2 \sum\limits_{i \in \mathcal{K}^{\prime}}  I_i \log I_i  + T+S \right), \\
    r  &=  4 \sigma_b^2  \sum\limits_{i \in \mathcal{K}^{\prime}} \frac{I_i^2 \log^2 I_i }{\sigma_i^2 } -  \left( 2 \sum\limits_{i \in \mathcal{K}^{\prime}}  I_i \log I_i + T+S \right)^2. 
\end{split}
\end{equation*}
Let $\alpha_i(T) = \frac{ (\lambda - 2 \log I_i)}{1 + T/S} $, for $i \in \mathcal{K}^{\prime}$, and substituting $\alpha_i(T)$ into \eqref{proof p1 w_i}. Then, Lemma 3 is proved.


\subsection{Proof of Lemma 4}


We now consider $\mathcal{C}_4$. For the best design  $b$, its budget allocation ratio $W_b(T)$ is always non-negative. As for non-best designs $i \in \mathcal{K}^{\prime}$, let $W_i(T) \geq 0$, and we have
\begin{align}
    \lambda & \geq 2 \log I_i \quad \forall i \in \mathcal{K}^{\prime}, \notag \\
    & \geq 2 \log I_{\langle k - 1 \rangle}. \notag 
\end{align}

If $w^*_b = 1/2$
\begin{equation*}
    \lambda = \frac{4\sum_{i \in \mathcal{K}^{\prime}} I_i \log I_i + T + S}{2\sum_{i \in \mathcal{K}^{\prime}} I_i} \geq 2 \log I_{\langle k - 1 \rangle},
\end{equation*}
and it can be checked that $T \geq 4 \sum_{i \in \mathcal{K}^{\prime}} I_i \log(I_{\langle k - 1 \rangle}/I_i) - S $.

If $w^*_b \neq 1/2$
\begin{equation*}
    \lambda = \frac{- q + \sqrt{q^2 - 4 p r}}{2 p} \geq 2 \log I_{\langle k - 1 \rangle},
\end{equation*}
where $p$, $q$, and $r$ are given in Lemma 3. 

Case 1: If $ p = 2w_b^* - 1 > 0$, i.e. $w_b^* > 1/2$, we need to solve the following inequality
\begin{equation}
\label{proof p2 case1}
    \sqrt{q^2 - 4 p r} \geq 4 p \log I_{\langle k - 1 \rangle} + q.
\end{equation}
Additionally
\begin{equation*}
\begin{split}
    &4 p \log I_{\langle k - 1 \rangle} + q  \\
    =& 4 \sum\limits_{i \in \mathcal{K}^{\prime}} \left(\frac{\sigma^2_b (w_i^*)^2}{\sigma_i^2} - (1-w_b^*)w_i^* \right) \log \frac{ I_{\langle k - 1 \rangle}}{I_i} + 2(1-w_b^*)\frac{T+S}{S}.\\ 
\end{split}
\end{equation*}
Furthermore, we define
\begin{equation*}
    T_1 = 2 \sum\limits_{i \in \mathcal{K}^{\prime}} \left[\frac{\sigma^2_b I_i^2}{\sigma_i^2 (S - I_b)} - I_i \right] \log \frac{I_{\langle k - 1 \rangle}}{I_i}  - S.
\end{equation*}
When $T <  T_1 $, $4 p \log I_{\langle k - 1 \rangle} + q < 0$, then Inequality \eqref{proof p2 case1} always holds. Otherwise, when $T \geq  T_1 $, we take square on both sides of \eqref{proof p2 case1} and obtain
\begin{equation*}
    T \in (- \infty, T_3] \cup [T_2,+\infty),
\end{equation*}
in which
\begin{equation*}
\begin{split}
    T_2 &= 2 \left[\sum\limits_{i \in \mathcal{K}^{\prime}}  I_i \log \frac{I_{\langle k - 1 \rangle}}{I_i}  +  \sigma_b \sqrt{  \sum\limits_{i \in \mathcal{K}^{\prime}} \frac{  I_i^2 }{\sigma_i^2 } \left(  \log \frac{I_{\langle k - 1 \rangle}}{I_i} \right)^2 }  \right]- S,\\
    T_3 &= 2 \left[\sum\limits_{i \in \mathcal{K}^{\prime}}  I_i \log \frac{I_{\langle k - 1 \rangle}}{I_i}  -  \sigma_b \sqrt{  \sum\limits_{i \in \mathcal{K}^{\prime}} \frac{  I_i^2 }{\sigma_i^2 } \left(  \log \frac{I_{\langle k - 1 \rangle}}{I_i}   \right)^2 }  \right]- S.\\
\end{split}
\end{equation*}
It can be checked that $T_2$ is strictly positive. Due to the fact $\lim_{T \rightarrow \infty} W_i(T) \rightarrow w_i^* \geq 0$, for $i \in \mathcal{K}^{\prime}$, we expect that $T$ can be sufficiently large. Therefore, a sufficient condition for $W(T)$ being feasible, when $w_b^* > 1/2$, is $T\geq \max\{0,T_1,T_2\}$.

Case 2: If $ p = 2w_b^* - 1 < 0$, i.e. $w_b^* < 1/2$, similarly, we need to solve the following inequality
\begin{equation}
\label{proof p2 case2}
    \sqrt{q^2 - 4 p r} \leq 4 p \log I_{\langle k - 1 \rangle} + q.
\end{equation}
When $T <  T_1 $, Inequality \eqref{proof p2 case2} never holds. Otherwise, we take square on both sides of \eqref{proof p2 case2} and, similarly, obtain $T \in (- \infty, T_3] \cup [T_2,+\infty)$. Therefore, a sufficient condition for $W(T)$ being feasible, when $w_b^* < 1/2$, is $T\geq \max\{0,T_1,T_2\}$.

Hence, the solution $W(T)$ is always feasible if $T\geq T_0$, where
\begin{equation*}
    T_0 =\left\{\begin{array}{ll}
\max\{0,T_1,T_2\}  & \quad \text{if} \  w^*_b \neq \frac{1}{2}  \\
\max\{0,4 \sum_{i \in \mathcal{K}^{\prime}} I_i \log\frac{I_{\langle k - 1 \rangle}}{I_i} - S \} & \quad \text{if} \ w^*_b = \frac{1}{2} \\
\end{array}\right.
\end{equation*}
These results conclude the proof.

\section{Proof of Propositions}

\subsection{Proof of Proposition 1}

We first show that for any pair of non-best designs $i,j \in \mathcal{K}^{\prime}$ and $i \neq j$, there exists a positive constant $c_{i,j} > 0$ such that $w_i/w_j \leq c_{i,j}$. We prove this by contradiction. Assume that there exists a pair of designs $i,j \in \mathcal{K}^{\prime}$ and $i \neq j$ such that $w_i/w_j$ can not be upper bounded, i.e., $w_i/w_j = \infty$. Since $w_i, w_j \in [0,1]$, it can be checked that $w_j \rightarrow 0$. By $\mathcal{C}_2$ in Lemma 3
\begin{equation}
\label{proof of proposition verifying: eqn 1}
\begin{split}
    &\left[\frac{\delta_{j,b}^2}{ (\sigma_j^2/w_j + \sigma_b^2/w_b) }  - \frac{\delta_{i,b}^2}{ (\sigma_i^2/w_i + \sigma_b^2/w_b) }  \right] \frac{T}{2}  = \log \frac{\delta_{j,b}  (\sigma_j^2/w^2_j) (\sigma_i^2/w_i + \sigma_b^2/w_b)^{\frac{3}{2}}}{\delta_{i,b}  (\sigma_i^2/w^2_i) (\sigma_j^2/w_j + \sigma_b^2/w_b)^{\frac{3}{2}} }. 
\end{split} 
\end{equation}
As $w_j \rightarrow 0$, the term $\frac{\delta_{j,b}^2}{ (\sigma_j^2/w_j + \sigma_b^2/w_b) }$ in \eqref{proof of proposition verifying: eqn 1} will vanish, and the right-hand side in \eqref{proof of proposition verifying: eqn 1} will approach infinity. Then, as $w_j \rightarrow 0$, we have
\begin{equation*}
    - \frac{\delta_{i,b}^2}{ 2(\sigma_i^2/w_i + \sigma_b^2/w_b) } T
    =  + \infty,
\end{equation*}
which is true if $\sigma_i^2/w_i + \sigma_b^2/w_b \rightarrow 0^-$. However, this contradicts that $\sigma_i^2/w_i + \sigma_b^2/w_b > 0$. Thus, the assumption is false, and for any pair of designs $i,j \in \mathcal{K}^{\prime}$ and $i \neq j$, there must exist a positive constant $c_{i,j} > 0$, such that $w_i/w_j \leq c_{i,j}$.

Let $w_{\text{min}}$ be the minimum budget allocation ratio of non-best designs, i.e., $w_{\text{min}} = \min_{i \in \mathcal{K}^{\prime}} w_i$. Therefore, there exists a positive constant $c$ such that $w_i \leq c w_{\text{min}}$, for $i \in \mathcal{K}^{\prime}$. By $\mathcal{C}_1$ in Lemma 3
\begin{equation*}
    w_b =  \sigma_b \sqrt{\sum\limits_{i \in \mathcal{K}^{\prime}} \frac{w_i^{2}}{\sigma_i^2}} \geq  \frac{ \sigma_b\sqrt{k-1} }{\bar{\sigma}}  w_{\text{min}}.
\end{equation*}
Then, for any non-best design $i \in \mathcal{K}^{\prime}$, we have
\begin{equation*}
    \frac{w_i}{w_b} \leq \frac{w_i}{w_{\text{min}}} \frac{\bar{\sigma}}{\sigma_b \sqrt{k-1} } \leq \frac{\bar{\sigma}}{\underline{\sigma}} \frac{c}{\sqrt{k-1}}.
\end{equation*}
This result concludes the proof, and therefore, $w_i/w_b = \mathcal{O}(1/\sqrt{k})$, for $i \in \mathcal{K}^{\prime}$.

\subsection{Proof of Proposition 2}
Let $w_{\text{min}}$ and $w_{\text{max}}$ be the minimum and maximum budget allocation ratios of non-best designs, respectively, i.e., $w_{\text{min}} = \min_{i\in \mathcal{K}^\prime} w_i$ and $w_{\text{max}} = \max_{i\in \mathcal{K}^\prime} w_i$. According to the proof of Proposition 1, we have $w_{\text{max}}/ w_{\text{min}} \leq c_1$ for some positive constants $c_1$. First, we show $(\log I_i) / (w_i/I_i) = \Theta(k)$, for $i\in \mathcal{K}^\prime$. On the one hand, for each $i\in \mathcal{K}^\prime$, we have 
\begin{align}
    \frac{\log I_i}{w_i/I_i} &= \frac{I_i\log I_i}{w_i} \notag\\
    & = \frac{\left(\sum_{i\in \mathcal{K}^\prime} w_i + w_b\right) I_i\log I_i}{w_i} \notag\\
    &= \frac{\left(\sum_{i\in \mathcal{K}^\prime} w_i + \sigma_b \sqrt{\sum_{i \in \mathcal{K}^{\prime}}w_i^{2} /\sigma_i^2}\right) I_i\log I_i}{w_i}. \label{proof of proposition 2-1}
\end{align}
If there exists a non-best design such that $\log I_i < 0$,
\begin{align*}
    \eqref{proof of proposition 2-1} &\geq \frac{ \left(k-1 + \frac{\bar{\sigma}}{\underline{\sigma}} \sqrt{k-1}\right) 
\frac{\bar{\sigma}^2}{\underline{\delta}^2} \log \frac{\underline{\sigma}^2}{\bar{\delta}^2} w_{\text{max}}  }{ w_{\text{min}}} \\
&\geq c_1 \frac{\bar{\sigma}^2}{\underline{\delta}^2} \left(\log \frac{\underline{\sigma}^2}{\bar{\delta}^2} \right) \left(k-1 + \frac{\bar{\sigma}}{\underline{\sigma}} \sqrt{k-1} \right)
 = lb_{1,1}.
\end{align*}
Otherwise,
\begin{align*}
    \eqref{proof of proposition 2-1} &\geq \frac{ \left(k-1 + \frac{\underline{\sigma}}{\bar{\sigma}} \sqrt{k-1}\right) 
\frac{\underline{\sigma}^2}{\bar{\delta}^2} \log \frac{\underline{\sigma}^2}{\bar{\delta}^2} w_{\text{min}}  }{ w_{\text{max}}} \\
& \geq \frac{1}{c_1} \frac{\underline{\sigma}^2}{\bar{\delta}^2} \left(\log \frac{\underline{\sigma}^2}{\bar{\delta}^2} \right) \left( k-1 + \frac{\underline{\sigma}}{\bar{\sigma}} \sqrt{k-1} \right) = lb_{1,2}.
\end{align*}
On the other hand, if there exists a non-best design such that $\log I_i > 0$,
\begin{align*}
    \eqref{proof of proposition 2-1} \leq c_1 \frac{\bar{\sigma}^2}{\underline{\delta}^2} \left(\log \frac{\bar{\sigma}^2}{\underline{\delta}^2}\right)  \left( k-1 + \frac{\bar{\sigma}}{\underline{\sigma}} \sqrt{k-1} 
  \right) = ub_{1,1}.
\end{align*}
Otherwise,
\begin{align*}
    \eqref{proof of proposition 2-1} \leq \frac{1}{c_1} \frac{\underline{\sigma}^2}{\bar{\delta}^2} \left(\log \frac{\bar{\sigma}^2}{\underline{\delta}^2}\right) \left(k-1 + \frac{\underline{\sigma}}{\bar{\sigma}} \sqrt{k-1}
  \right) = ub_{1,2}.
\end{align*}
Based on preceding analyses, for $i \in \mathcal{K}^\prime$, we have
\begin{align*}
    \min\{ lb_{1,1}, lb_{1,2} \} \leq \frac{\log I_i}{ w_i / I_i} \leq \max\{ ub_{1,1}, ub_{1,2} \}.
\end{align*}
Since $lb_{1,1}$, $lb_{1,2}$, $ub_{1,1}$, and $ub_{1,2}$ are all $\Theta(k)$, $(\log I_i)/ (w_i/I_i)$, for $i \in \mathcal{K}^\prime$, are also $\Theta(k)$.

Second, we show $(\log w_i) / (w_i/I_i) = \Theta(k \log k)$, for $i\in \mathcal{K}^\prime$. For each $i\in \mathcal{K}^\prime$, we have 
\begin{align}
    \frac{\log w_i}{w_i/I_i} &= \frac{I_i\log w_i}{w_i} \notag\\
    & = \frac{\sum_{i\in \mathcal{K}^\prime} w_i + w_b}{w_i} 
 \cdot  I_i \log \frac{w_i}{\sum_{i\in \mathcal{K}^\prime} w_i + w_b}. \label{proof of proposition 2-2}
\end{align}
On the one hand, we have 
\begin{align*}
    \eqref{proof of proposition 2-2} &\geq \frac{\left(k-1 + \frac{\bar{\sigma}}{\underline{\sigma}} \sqrt{k-1}\right) w_{\text{max}}}{w_{\text{min}}} \cdot  \frac{\bar{\sigma}^2}{\underline{\delta}^2} \log \frac{w_{\text{min}}}{\left(k-1 + \frac{\bar{\sigma}}{\underline{\sigma}} \sqrt{k-1}\right) w_{\text{max}}} \\
    & \geq -  \frac{\bar{\sigma}^2}{\underline{\delta}^2}c_1 \left(k-1 + \frac{\bar{\sigma}}{\underline{\sigma}} \sqrt{k-1}\right) \cdot \log c_1  \left(k-1 + \frac{\bar{\sigma}}{\underline{\sigma}} \sqrt{k-1}\right)  = lb_2.
\end{align*}
On the other hand, we have 
\begin{align*}
    \eqref{proof of proposition 2-2} & \leq \frac{\left(k-1 + \frac{\underline{\sigma}}{\bar{\sigma}} \sqrt{k-1}\right) w_{\text{min}}}{w_{\text{max}}} \cdot  \frac{\underline{\sigma}^2}{\bar{\delta}^2} \log \frac{w_{\text{max}}}{\left(k-1 + \frac{\underline{\sigma}}{\bar{\sigma}} \sqrt{k-1}\right) w_{\text{min}}} \\
    & \leq  \frac{\underline{\sigma}^2}{\bar{\delta}^2} \frac{1}{c_1} \left(k-1 + \frac{\underline{\sigma}}{\bar{\sigma}} \sqrt{k-1}\right) \cdot \log c_1  \left(k-1 + \frac{\underline{\sigma}}{\bar{\sigma}} \sqrt{k-1}\right) = ub_2.
\end{align*}
Since both $lb_2$ and $ub_2$ are $\Theta(k \log k)$, $(\log w_i) / (w_i/I_i) = \Theta(k \log k)$, for $i\in \mathcal{K}^\prime$.

Therefore, as $k \to \infty$, if $T = \omega(k \log k)$, both $\log I_i$ and $\log w_i$ are negligible to $Tw_i/I_i$, because $ (\log I_i)/(Tw_i/I_i) \to 0$ and $ (\log w_i)/(Tw_i/I_i) \to 0$. These results complete the proof.

\subsection{Proof of Proposition 3}
According to Lemma 3, we have
\begin{equation*}
    \alpha_i(T) = \frac{ (\lambda - 2 \log I_i)}{1 + T/S}, \quad i \in \mathcal{K}^{\prime},
\end{equation*}
which is monotone decreasing with respect to $I_i$. Since $I_{\langle 1 \rangle} \leq I_{\langle 2 \rangle} \leq \dots \leq I_{\langle k - 1 \rangle}$, we have $\alpha_{\langle 1 \rangle}(T) \geq \alpha_{\langle 2 \rangle}(T) \geq \dots \geq \alpha_{\langle k - 1 \rangle}(T)$. Furthermore, if $I_{i}$, for $i \in \mathcal{K}^{\prime}$, are all equal, then $\lambda = 1 + T/S + 2 \log I_i$ and $\alpha_{\langle 1 \rangle}(T) = \alpha_{\langle 2 \rangle}(T) = \dots = \alpha_{\langle k - 1 \rangle}(T) = 1$. On the contrary, if $\alpha_{\langle 1 \rangle}(T) = \alpha_{\langle 2 \rangle}(T) = \dots = \alpha_{\langle k - 1 \rangle}(T)$, it can be checked that all $I_i$s are equal.

We show $\alpha_{\langle 1 \rangle}(T) \geq 1 $ by contradiction, and $\alpha_{\langle k - 1 \rangle}(T) \leq 1$ can be proved similarly. Without lose of generality, we assume $\alpha_{\langle 1 \rangle}(T) < 1$. Based on preceding analyses, we know $\alpha_{\langle k - 1 \rangle}(T) \leq \alpha_{[k-2]}(T) \leq \dots \leq \alpha_{\langle 1 \rangle}(T) < 1$. According to Lemma 3, we have $W_i(T) = w^*_i \alpha_{i}(T) < w^*_i $ for non-best designs $i \in \mathcal{K}^{\prime}$, and $W_b(T) = \sigma_b \sqrt{\sum_{i \in \mathcal{K}^{\prime}}\frac{(w_i^* \alpha_{i}(T) )^{2}  }{\sigma_i^2}} < w^*_b$ for the best design $b$. Then we have $\sum_{i \in \mathcal{K}} W_b(T) < \sum_{i \in \mathcal{K}} w^*_b = 1$, which contradicts that $\sum_{i \in \mathcal{K}} W_b(T) = 1$. Therefore, $\alpha_{\langle 1 \rangle}(T)$ must be larger or equal to 1. These results conclude the proof.
\endgroup

\subsection{Proof of Proposition 5}

Consider a probability space $(\varOmega,\mathcal{F},\mathbb{P})$ where all Gaussian random variables are well defined. $\exists \tilde{\varOmega} \in \varOmega$, $\tilde{\varOmega}$ is measurable and $\mathbb{P}(\tilde{\varOmega}) = 1$. By the strong Law of Large Numbers, for any sample path $\omega \in \tilde{\varOmega}, \forall i \in \mathcal{K}, \lim_{t \rightarrow T}\hat{\mu}_i^{(t)} \rightarrow \mu_i$ almost surely if $N_i^{(t)} \rightarrow \infty$ as $t \rightarrow T$. Additionally, $\lim_{t \rightarrow T} (\hat{\sigma}_i^{(t)})^2 \rightarrow (\hat{\sigma_i})^2 $ almost surely by the same argument. We fix a sample path $\omega \in \tilde{\varOmega}$. For notational simplicity, we omit the dependence of terms on $\omega$.

We first prove Proposition 5(1).  
Let $A = \{i| N^{(t)}_i \rightarrow \infty,\text{as} \ t \rightarrow T \}$. Suppose that there exists a design $i_0 \notin A$, that is, $\exists \xi > 0$, design $i_0$ will not receive samples when $t > \xi$. For any $i \in \mathcal{K}$, denote by
\begin{equation*}
\begin{split}
    \hat{w}_i &= \lim_{t \rightarrow T} \frac{1}{t} \cdot (t+1) \cdot \widetilde{W}_i^{(t)}(t+1) \\
    & = \left(1 + \frac{1}{T}\right) \lim_{t \rightarrow T} \hat{w}_i^{*,(t)} \cdot \hat{\alpha}_i^{(t)}(t+1) \\
    & = \lim_{t \rightarrow T} \hat{w}_i^{*,(t)} > 0,
\end{split}
\end{equation*}
where $\hat{\alpha}_i^{(t)}(t+1)$ is the estimated value of $\alpha_i(t+1)$ with consumed budget $t$, and the third equality holds because we consider the setting where $k \to \infty$ and $T = \omega(k \log k)$. Then, we have $\sum_{i \in \mathcal{K}} \hat{w}_i = 1$. Note that,
\begin{equation}
\label{proof p4-1}
    \sum_{i \in \mathcal{K}} \lim_{t \rightarrow T} \frac{N_i^{(t)}}{t} = \sum_{i \in A} \lim_{t \rightarrow T} \frac{N_i^{(t)}}{t} = 1.
\end{equation}
Otherwise, $\lim_{t \rightarrow T} N_{i_0}^{(t)}/t > 0$. It indicates that $\exists \xi > 0$, design $i_0$ will receive samples when $t > \xi$. This contradicts that $i_0 \notin A$. Additionally, it can be checked that for any $i \in A$,
\begin{equation}
\label{proof p4-2}
    \hat{w}_i \geq \lim_{t \rightarrow T} \frac{N_i^{(t)}}{t} + \hat{w}_{i_0}.
\end{equation}
Otherwise, we have,
\begin{equation*}
\begin{split}
    \hat{w}_i - \lim_{t \rightarrow T} \frac{N_i^{(t)}}{t} &< \hat{w}_{i_0} - \lim_{t \rightarrow T} \frac{N_{i_0}^{(t)}}{t} \\
    \implies \lim_{t \rightarrow T} \frac{(t+1) \cdot \widetilde{W}_i^{(t)}(t+1) - N_i^{(t)}}{t} &< \lim_{t \rightarrow T} \frac{(t+1) \cdot \widetilde{W}_{i_0}^{(t)}(t+1) - N_{i_0}^{(t)}}{t}.
\end{split}
\end{equation*}
This implies that $\exists \xi_1 > 0$, design $i$ will not receive samples when $t > \xi_1$, which contradicts that $i \in A$. Then, with \eqref{proof p4-1} and \eqref{proof p4-2}, it can be checked that,
\begin{equation*}
    1 = \sum_{i \in A} \hat{w}_i \geq 1 + \sum_{i \in A} \hat{w}_{i_0} > 1.
\end{equation*}
This contradicts that $\sum_{i \in A} \hat{w}_i = 1$. Therefore, $\forall i \in \mathcal{K}$, $N_i^{(t)} \rightarrow \infty$ and $\hat{b}^{(t)} = b$ as $t \rightarrow T$.

Before proving Proposition 4(2), we introduce a lemma.\\
\textbf{Lemma 5.} \citep{li2023convergence}
\textit{
Let $\{ N_i^{(t)}| i \in \mathcal{K}, t = 1,2,\dots \}$ be a sequence of positive integers that satisfies $N_i^{(t)} \rightarrow \infty$ as $t \rightarrow \infty$. Denote by $w_i^{(t)} = N_i^{(t)}/\sum_{j \in \mathcal{K}} N_j^{(t)}$ and $M_i^{(t)} = M_i^{(t-1)} + \mathbbm{1}\{ 
N_i^{(t)} > N_i^{(t-1)}\}$, $\forall i \in \mathcal{K}$, where $\mathbbm{1}\{\cdot\}$ is an indicator function and $M_i^{(0)} = 0$, $\forall i \in \mathcal{K}$. If each subsequence $\{ w_i^{(t_p)} | M_i^{(p)} \rightarrow \infty \ \text{as} \ p \rightarrow \infty, i \in \mathcal{K}, p = 1,2,\dots \}$ of $\{ w_i^{(t)} | i \in \mathcal{K} \}$ has a convergent subsequence and the convergent subsequence converges to $\{w_i|i\in \mathcal{K}\}$, then $\lim_{t \rightarrow \infty} w_i^{(t)} = w_i$, $\forall i \in \mathcal{K}$.
}

According to Bolzano–Weierstrass theorem, we know that each subsequence $\{ w_i^{(t_p)} | M_i^{(p)} \rightarrow \infty \ \text{as} \ p \rightarrow \infty, i \in \mathcal{K}, p = 1,2,\dots \}$ has a convergent subsequence. Take any convergent subsequence $\{ w_i^{(t_{p_q})} | M_i^{(q)} \rightarrow \infty \ \text{as} \ q \rightarrow \infty, i \in \mathcal{K}, q = 1,2,\dots \}$. Suppose that it converges to $\{ w_i^\prime| i \in \mathcal{K} \} \neq \{ w_i^*| i \in \mathcal{K} \}$. Then, $\exists i_1$, $w_{i_1}^* < w_{i_0}^\prime$. It implies that $\exists \xi_2 > 0$, design $i_0$ will not receive samples when $q > \xi_2$, which contradicts that $M_{i_0}^{(q)} \rightarrow \infty$ as $q \rightarrow \infty$. Then, we must have $\{ w_i^\prime| i \in \mathcal{K} \} = \{ w_i^*| i \in \mathcal{K} \}$, and therefore, $\lim_{t \rightarrow T} w^{(t)}_i = w_i^* $, $\forall i \in \mathcal{K}$ (because Lemma 5 applies).  

\section{Illustrations on proportional allocations to designs for Example 1 with different budgets}

\begin{figure}[H]
     \centering
     \begin{subfigure}[b]{0.32\textwidth}
         \centering
         \includegraphics[width=\textwidth]{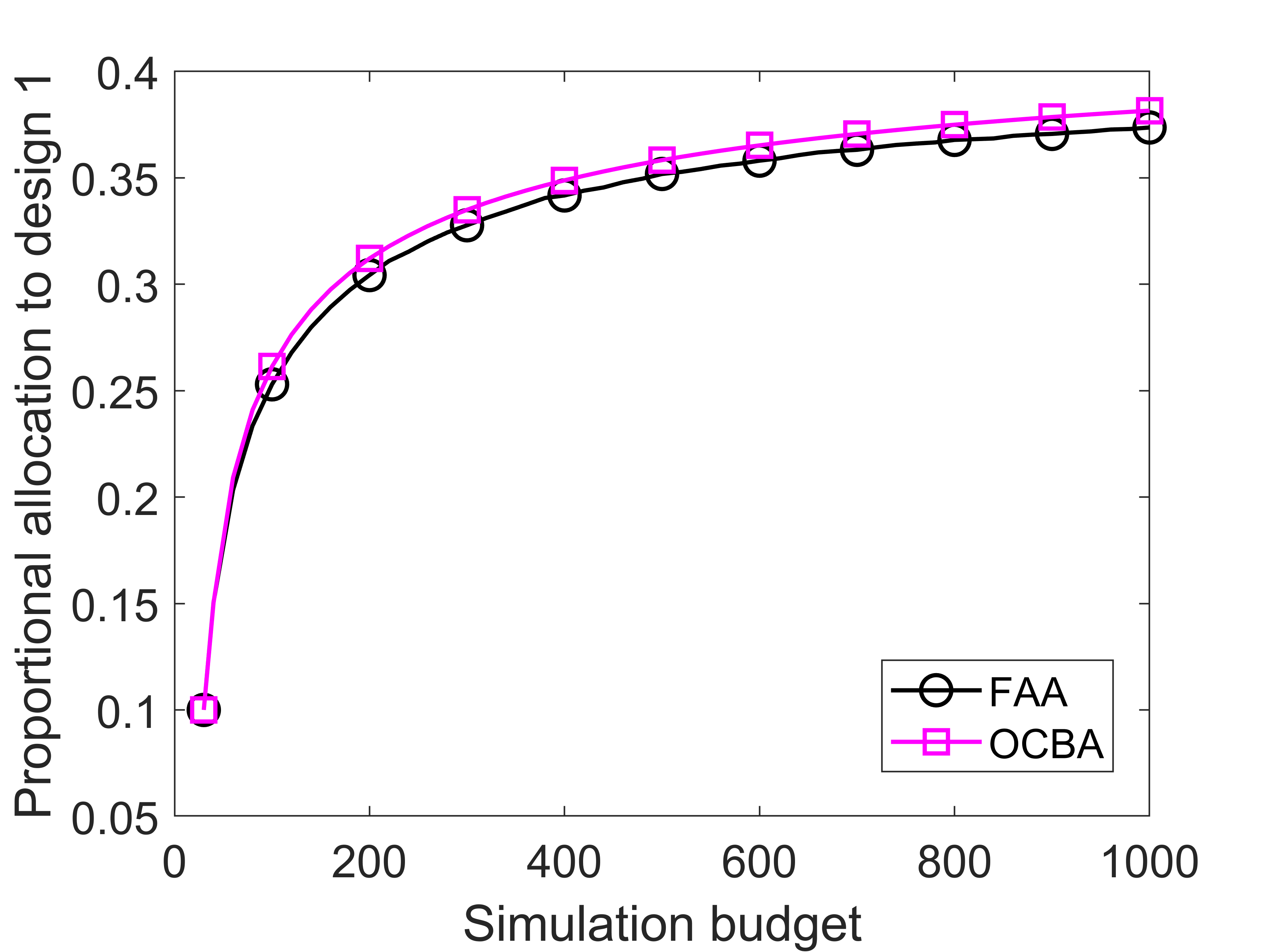}
         \caption{Proportional allocations to design 1}
         \label{sfig:Proportional allocations to design 1}
     \end{subfigure}
     \hfill
     \begin{subfigure}[b]{0.32\textwidth}
         \centering
         \includegraphics[width=\textwidth]{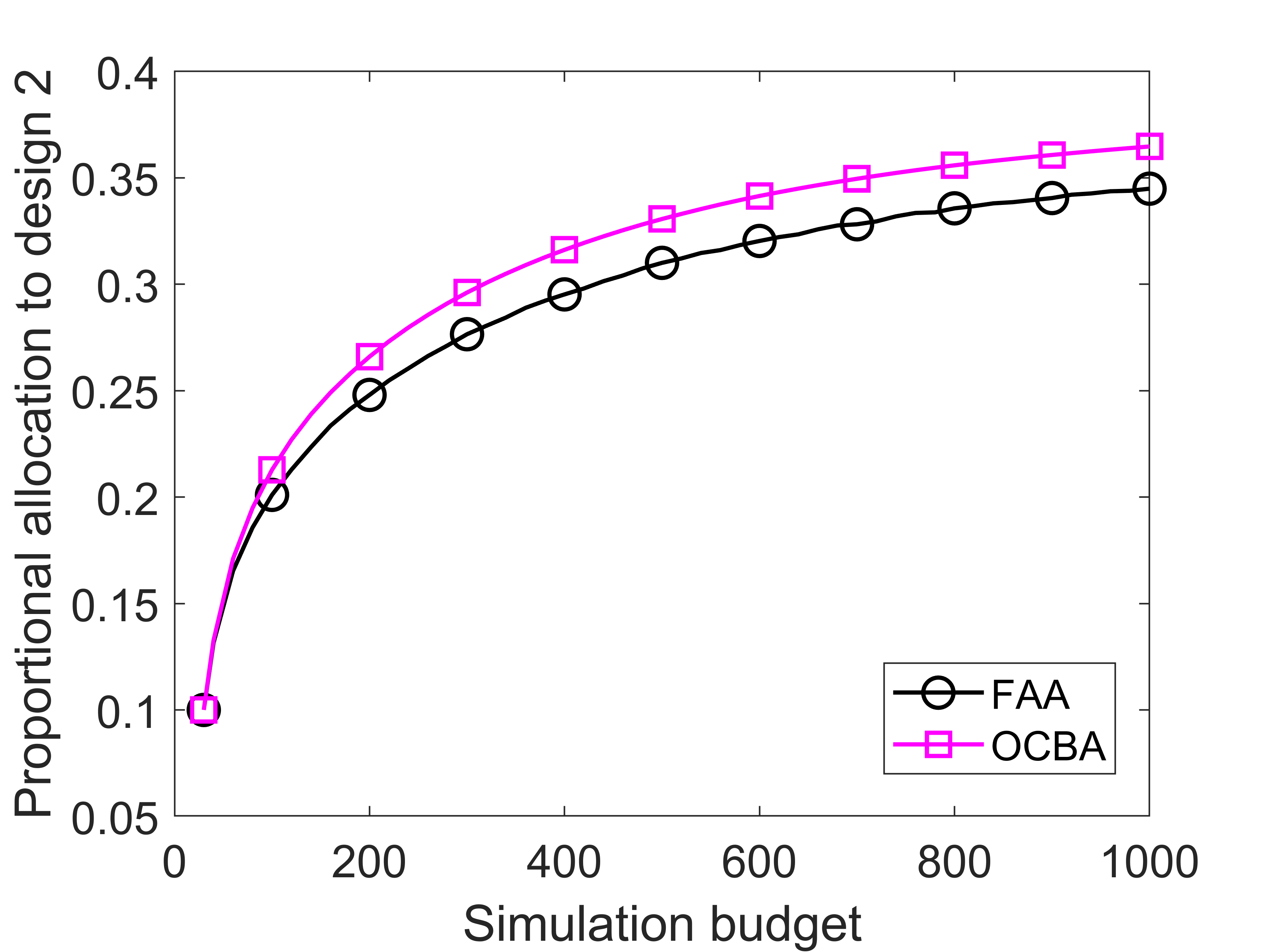}
         \caption{Proportional allocations to design 2}
         \label{sfig:Proportional allocations to design 2}
     \end{subfigure}
     \hfill
     \begin{subfigure}[b]{0.32\textwidth}
         \centering
         \includegraphics[width=\textwidth]{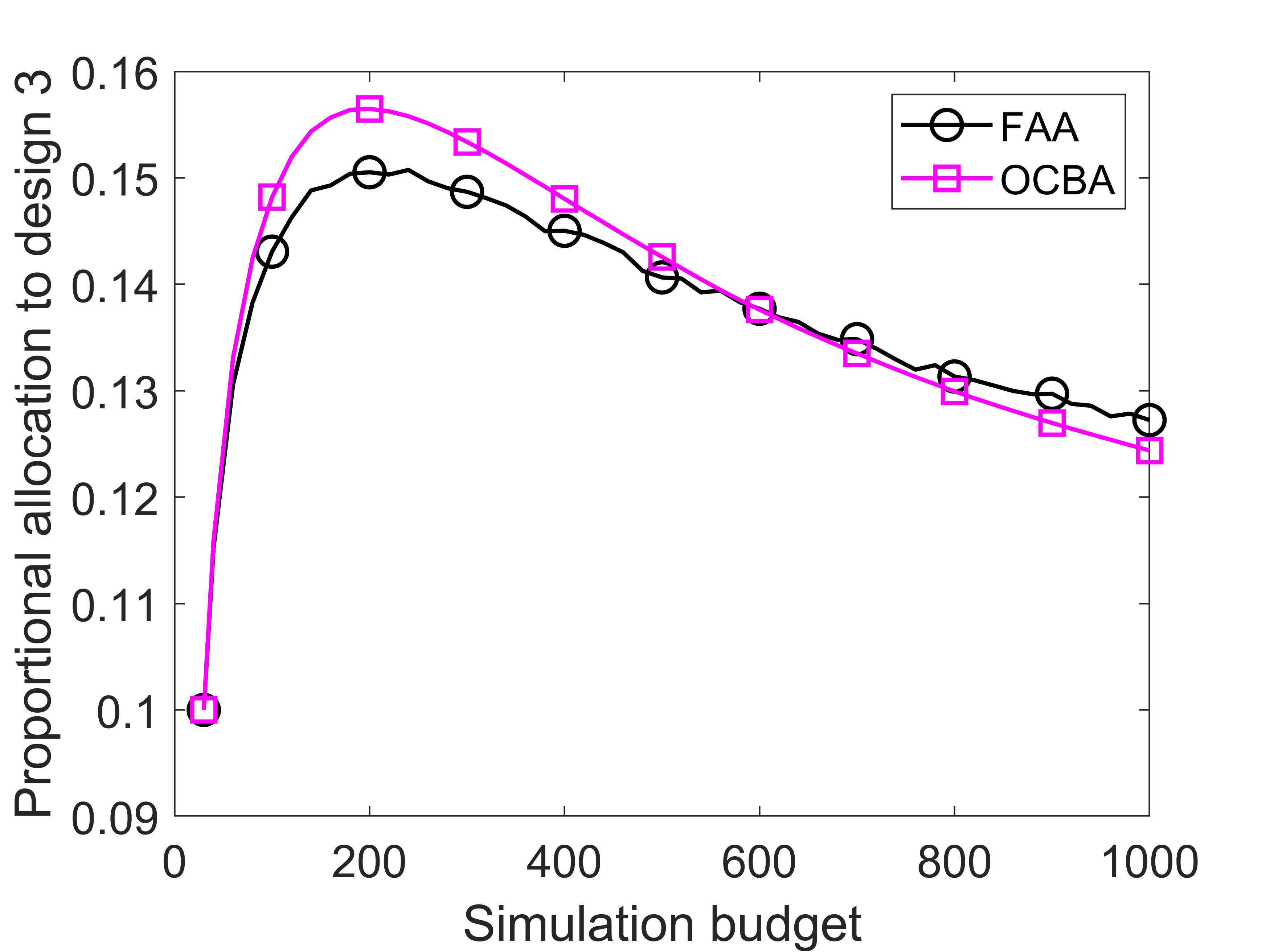}
         \caption{Proportional allocations to design 3}
         \label{sfig:Proportional allocations to design 3}
     \end{subfigure}
     \hfill
     \begin{subfigure}[b]{0.32\textwidth}
         \centering
         \includegraphics[width=\textwidth]{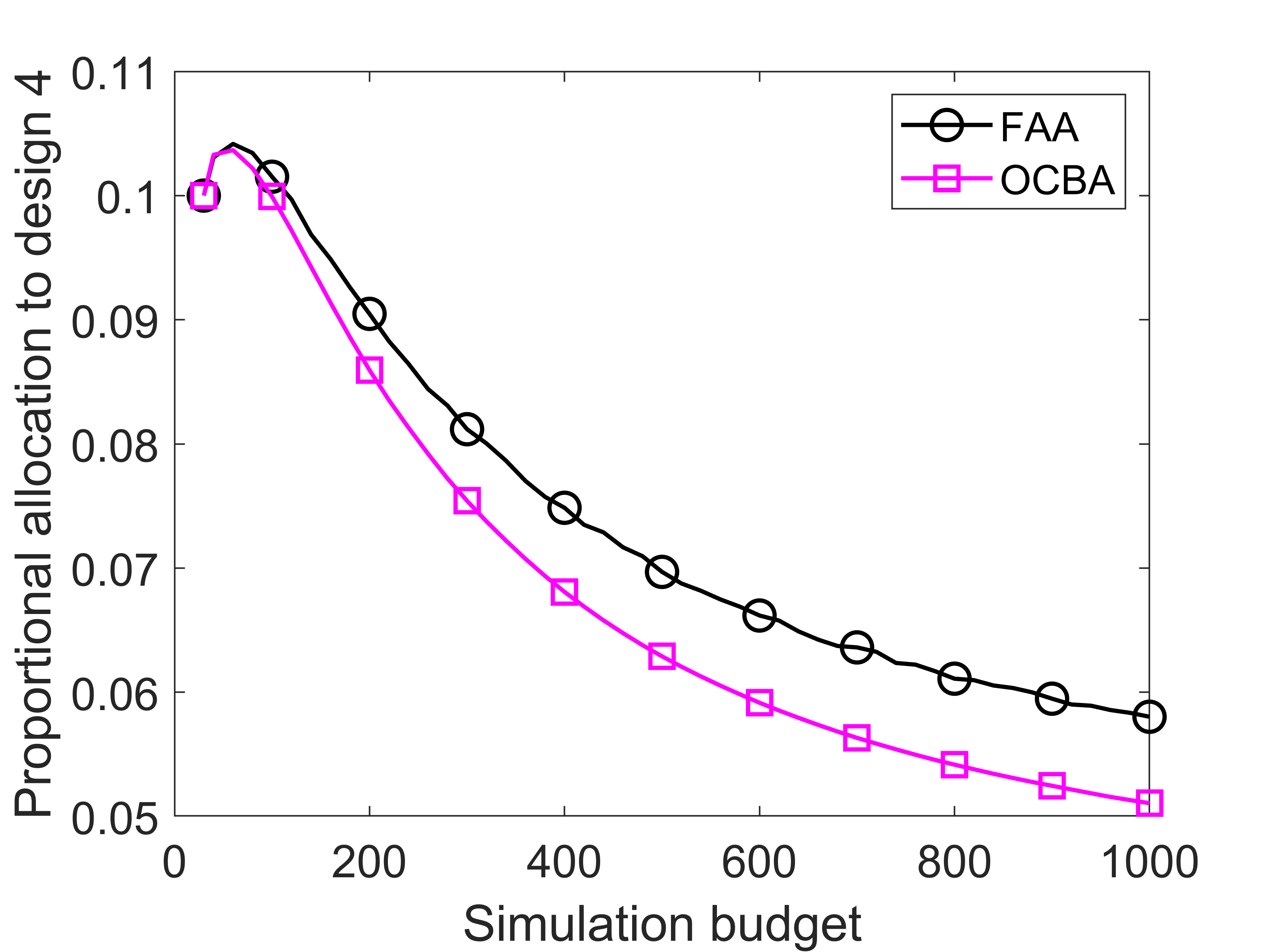}
         \caption{Proportional allocations to design 4}
         \label{sfig:Proportional allocations to design 4}
     \end{subfigure}
     \hfill
     \begin{subfigure}[b]{0.32\textwidth}
         \centering
         \includegraphics[width=\textwidth]{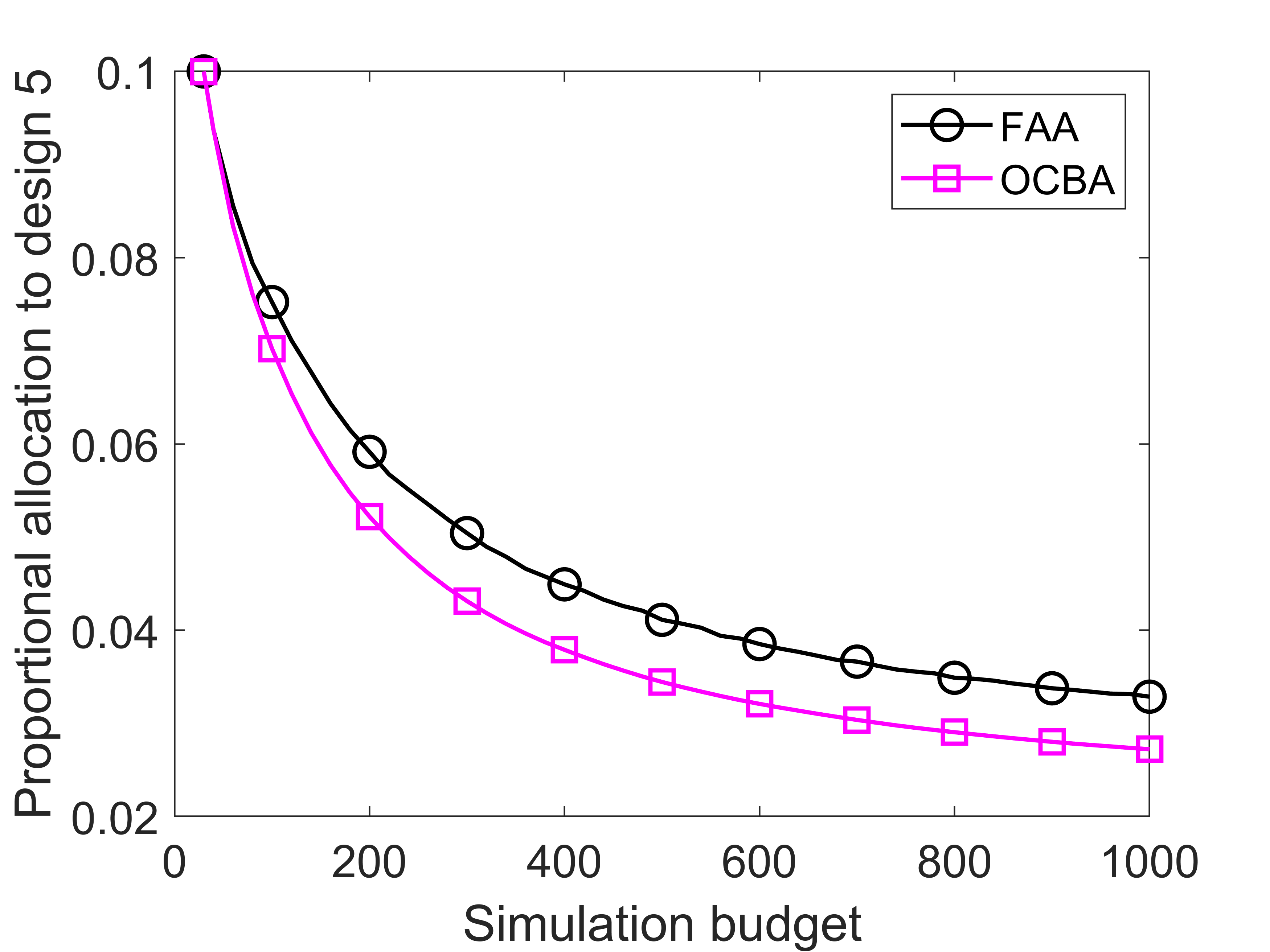}
         \caption{Proportional allocations to design 5}
         \label{sfig:Proportional allocations to design 5}
     \end{subfigure}
     \hfill
     \begin{subfigure}[b]{0.32\textwidth}
         \centering
         \includegraphics[width=\textwidth]{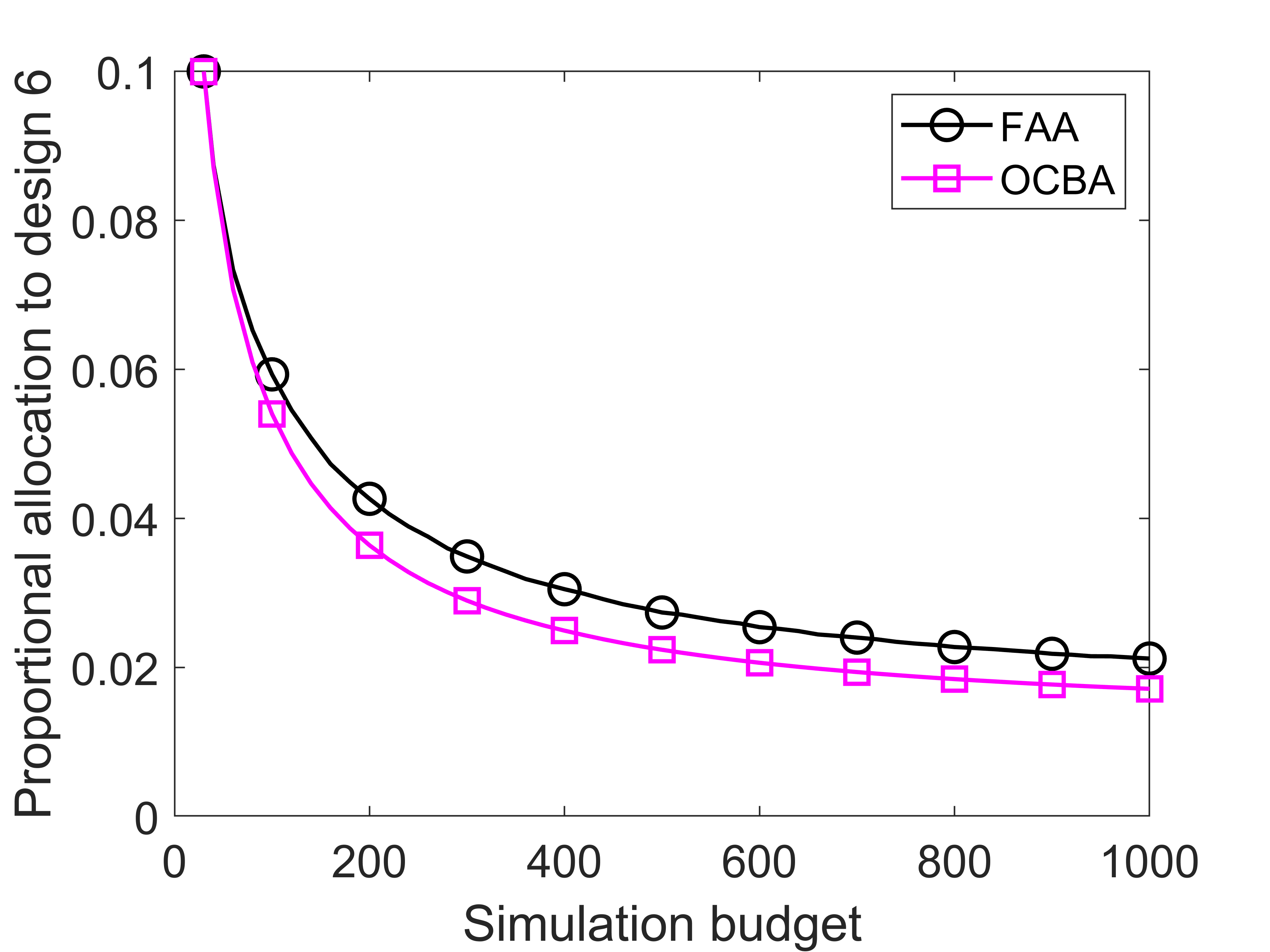}
         \caption{Proportional allocations to design 6}
         \label{sfig:Proportional allocations to design 6}
     \end{subfigure}
     \hfill
     \begin{subfigure}[b]{0.32\textwidth}
         \centering
         \includegraphics[width=\textwidth]{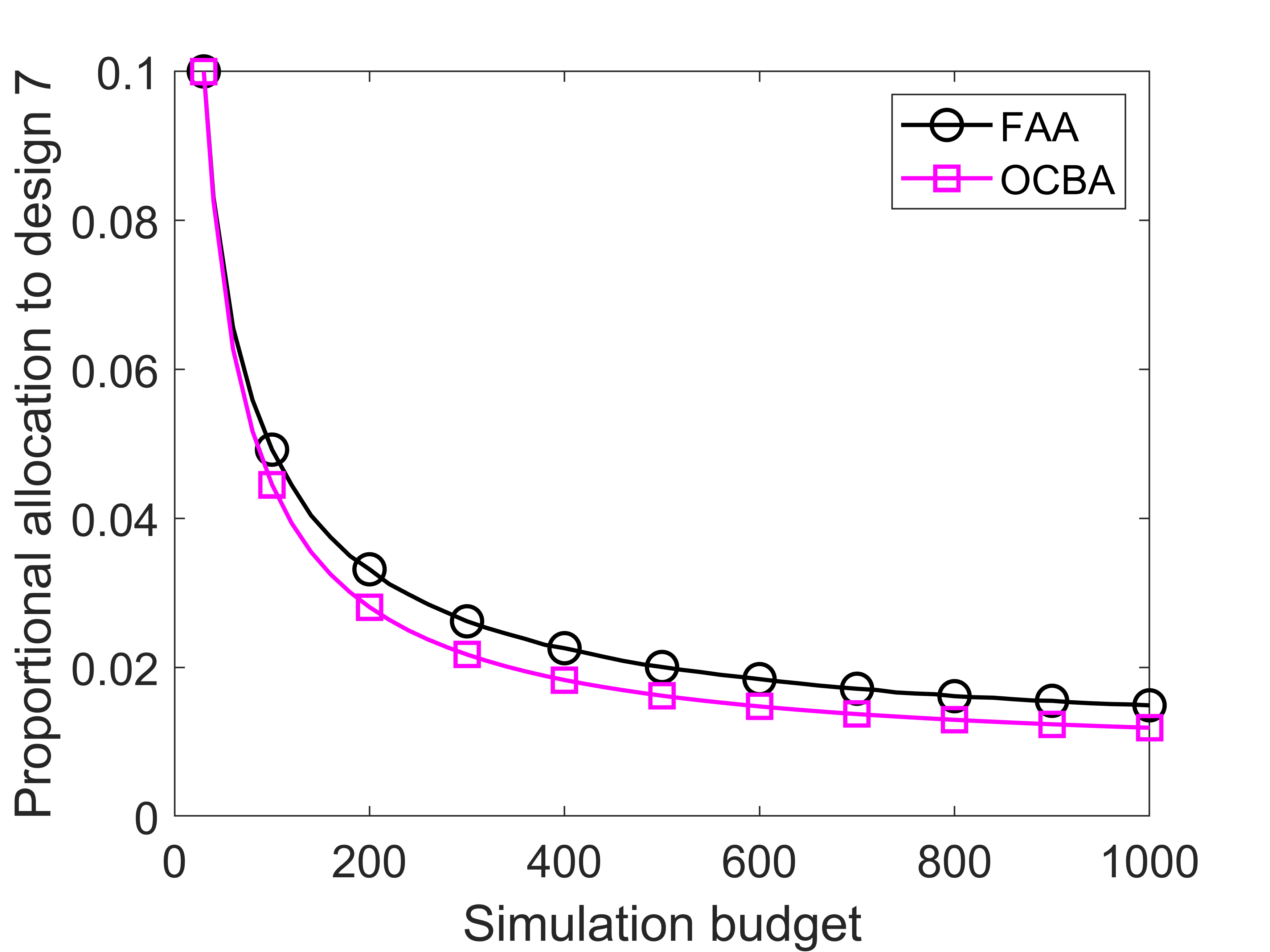}
         \caption{Proportional allocations to design 7}
         \label{sfig:Proportional allocations to design 7}
     \end{subfigure}
     \hfill
     \begin{subfigure}[b]{0.32\textwidth}
         \centering
         \includegraphics[width=\textwidth]{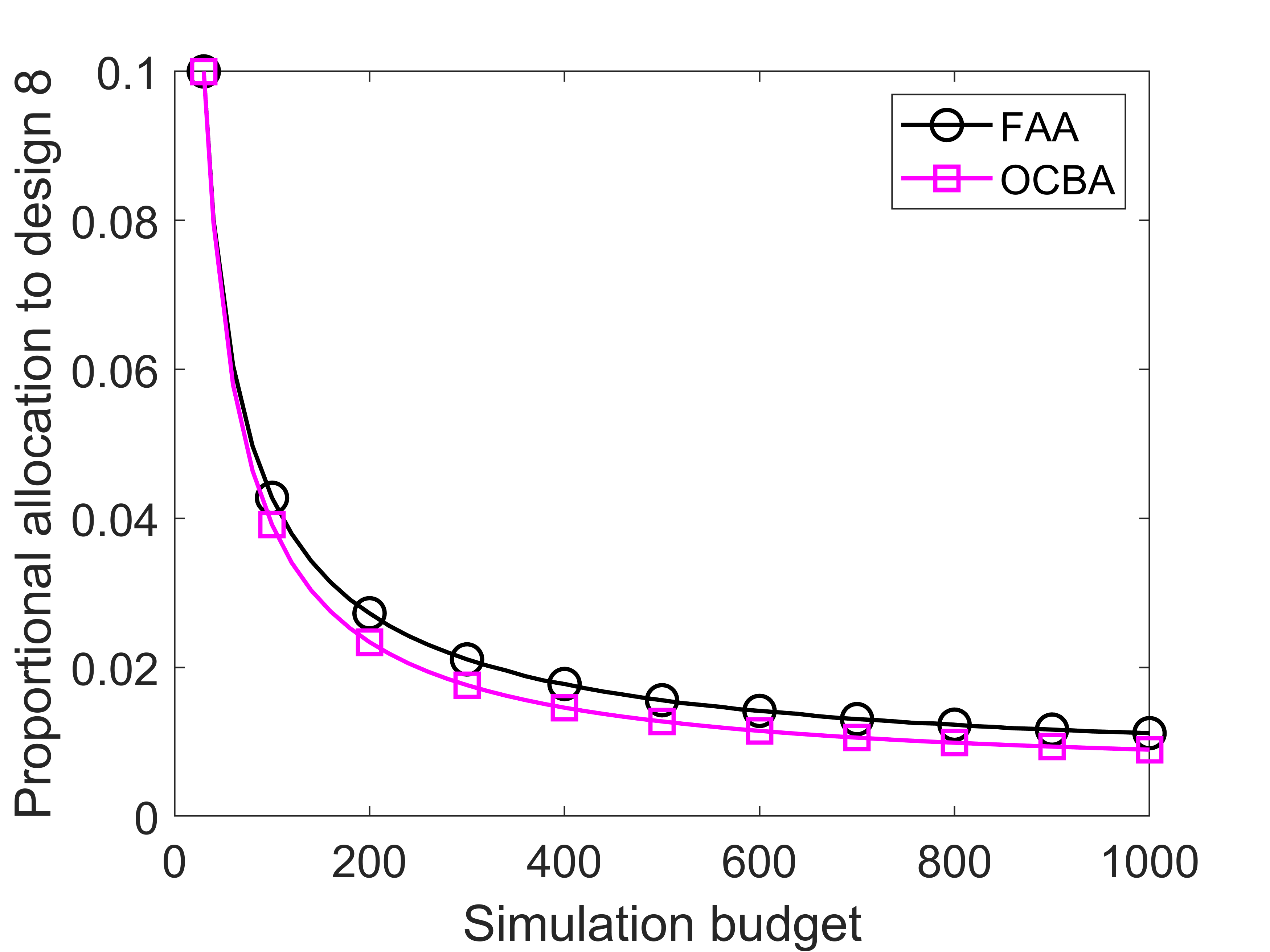}
         \caption{Proportional allocations to design 8}
         \label{sfig:Proportional allocations to design 8}
     \end{subfigure}
     \hfill
     \begin{subfigure}[b]{0.32\textwidth}
         \centering
         \includegraphics[width=\textwidth]{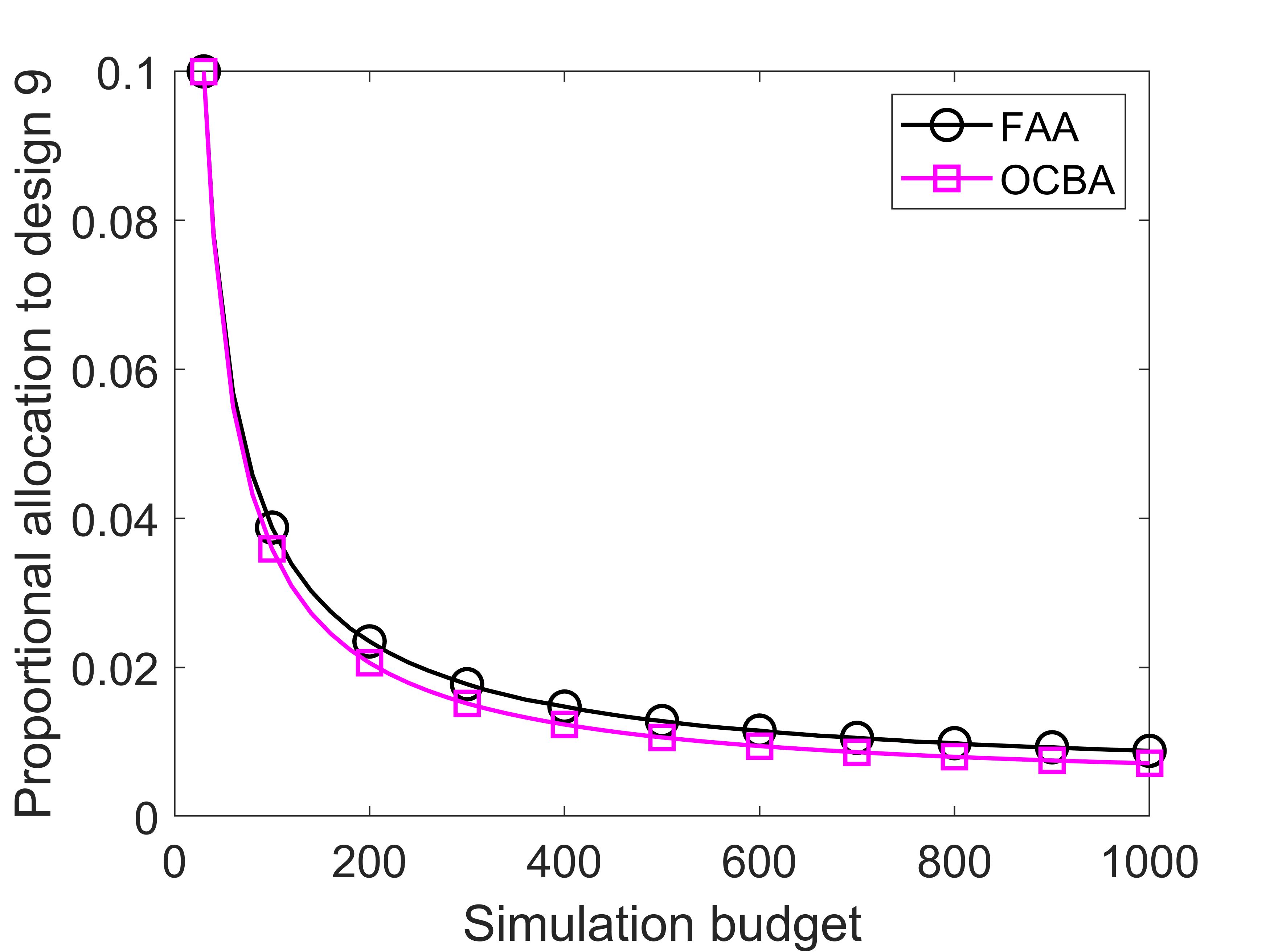}
         \caption{Proportional allocations to design 9}
         \label{sfig:Proportional allocations to design 9}
     \end{subfigure}
     \hfill
     \begin{subfigure}[b]{0.32\textwidth}
         \centering
         \includegraphics[width=\textwidth]{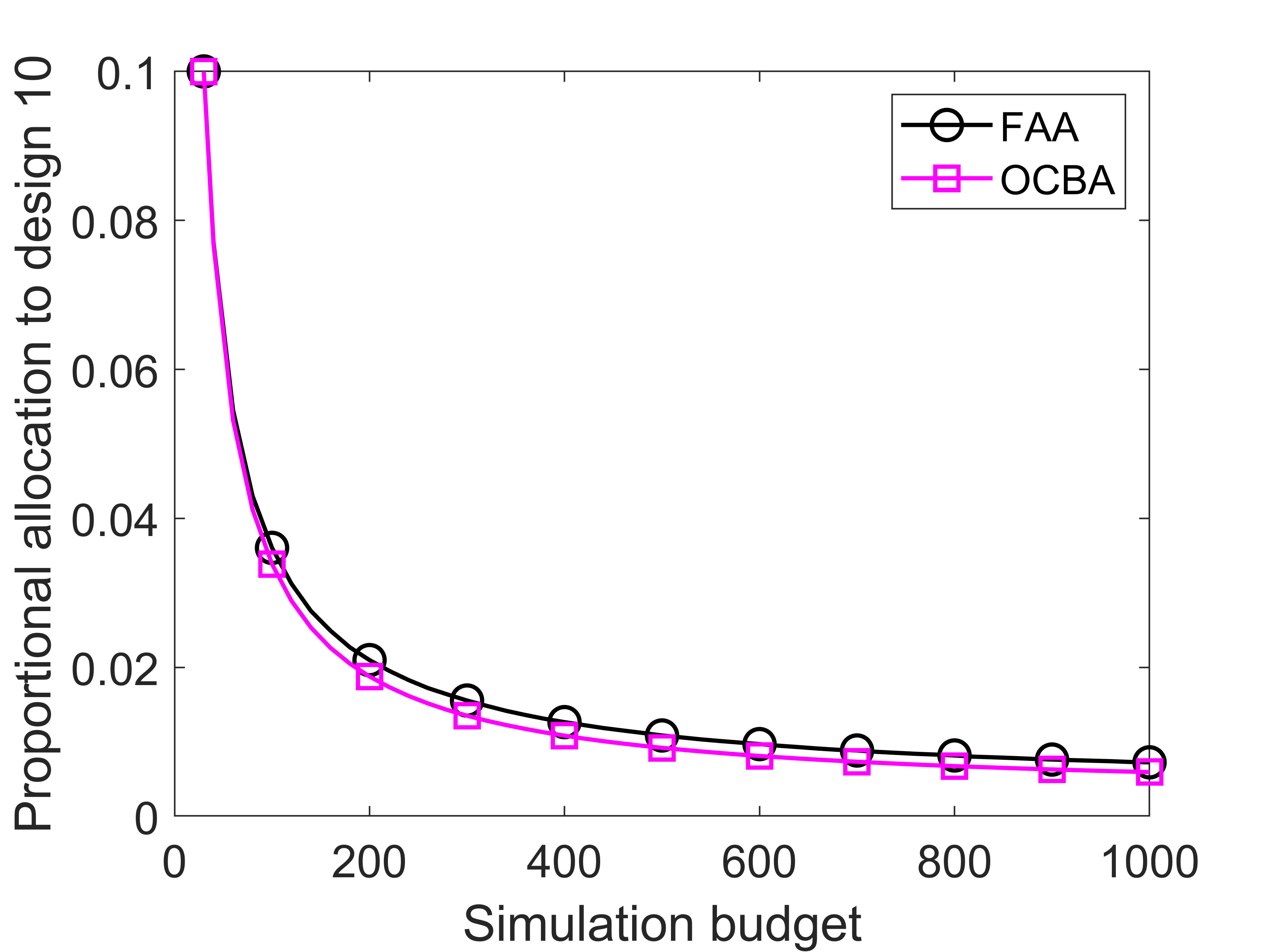}
         \caption{Proportional allocations to design 10}
         \label{sfig:Proportional allocations to design 10}
     \end{subfigure}
     \caption{(Color online) Illustration on proportional allocations to designs for Example 1 with different budgets}
    \label{fig:Illustration on proportional allocations to designs for Example 1 with different budgets}
\end{figure}

Figure~\ref{fig:Illustration on proportional allocations to designs for Example 1 with different budgets} illustrates the proportional allocation to designs made by FAA and OCBA for Example 1 with different budgets. Only the proportional allocations to designs made by FAA are illustrated for comparison with OCBA. Because FAA anchors the final budget, and its allocation policy does not change during the allocation procedure. Note that, for each budget size $T$, we run FAA with $T$ (anchored by FAA) and obtain the corresponding proportional allocations to designs made by FAA (with anchored budget $T$).

Clearly, the proposed budget-adaptive allocation rule behaves differently from OCBA with different budgets. Specifically, compared with OCBA, it discounts the proportions of total budget allocated to competitive designs (e.g., design 2) while increases the proportions of total budget allocated to non-competitive designs (e.g., design 10). Note that, these adjustments to proportional allocations to designs are effective across different budget sizes, e.g., between a budget size of 100 and 110. It is interesting to see that the proportion of total budget allocated to design 3 made by FAA is initially higher than that of OCBA, but later it is surpassed by OCBA’s allocation. This observation indicates that $\alpha_i(T)$ is not monotonic with respect to $T$.

\bibliographystyle{apalike}

\bibliography{reference}

\end{document}